\documentclass[10pt]{article}%
\usepackage{amsmath,amssymb}
\usepackage{amsfonts}
\usepackage{graphicx}
\usepackage{amsmath}
\usepackage{amssymb}%
\numberwithin{equation}{section} 
\def\beq{\begin{equation}}
\def\eeq{\end{equation}}

\def\cB{ {{\cal B}^n}}

\def\cF{ {{\cal F}^n}}

\def\cP{ {{\cal P}}}

\def\bE{ {{\mathbb{E}}}}
\def\tr{ {{\rm{tr}}} }

\renewcommand{\min}[2]{#1 \wedge #2}
\newtheorem{defn}{{\bf Definition}}[section]
\newtheorem{thm}[defn]{{\bf Theorem}}
\newtheorem{cor}[defn]{{\bf Corollary}}
\newtheorem{lem}[defn]{{\bf Lemma}}
\newtheorem{prop}[defn]{{\bf Proposition}}
\newtheorem{rem}[defn]{{\bf Remark}}

\newtheorem{notation}[defn]{Notation}
\newenvironment{proof}[1][Proof]{\textbf{#1.} }{\hfill \rule{0.5em}{0.5em}}
\makeatletter \@ifundefined{pdfoutput} {
\DeclareGraphicsRule{.wmf}{bmp}{}{}
\DeclareGraphicsRule{.jpg}{bmp}{}{}
\DeclareGraphicsRule{.png}{bmp}{}{}
\DeclareGraphicsRule{.cdr}{bmp}{}{}
\DeclareGraphicsRule{.gif}{bmp}{}{} } { \ifnum\pdfoutput=0\relax
\DeclareGraphicsRule{.wmf}{bmp}{}{}
\DeclareGraphicsRule{.jpg}{bmp}{}{}
\DeclareGraphicsRule{.png}{bmp}{}{}
\DeclareGraphicsRule{.cdr}{bmp}{}{}
\DeclareGraphicsRule{.gif}{bmp}{}{} \fi \ifnum\pdfoutput=1\relax
\DeclareGraphicsRule{.eps}{pdf}{}{}
\DeclareGraphicsRule{.wmf}{jpg}{}{} \fi } \makeatother
\begin{document}

\title{Path Integrals on a Compact Manifold with Non-negative Curvature}
\author{Adrian P. C. Lim \\
Department of Mathematics, Cornell University \\
Email: al382@cornell.edu} \date{} \maketitle

\begin{abstract}
A typical path integral on a manifold, $M$ is an informal
expression of the form \begin{equation} \frac{1}{Z}\int_{\sigma
\in H(M)} f(\sigma) e^{-E(\sigma)}\mathcal{D}\sigma, \nonumber
\end{equation} where $H(M)$ is a Hilbert manifold of paths with energy $E(\sigma) < \infty$,
$f$ is a real valued function on $H(M)$, $\mathcal{D}\sigma$ is a
\textquotedblleft Lebesgue measure \textquotedblright\ and $Z$ is
a normalization constant. For a compact Riemannian manifold $M$,
we wish to interpret $\mathcal{D}\sigma$ as a Riemannian
\textquotedblleft volume form \textquotedblright\ over $H(M)$,
equipped with its natural $G^{1}$ metric. Given an equally spaced
partition,
${\mathcal{P}}$ of $[0,1],$ let $H_{{\mathcal{P}}%
}(M)$ be the finite dimensional Riemannian submanifold of $H\left(
M\right) $ consisting of piecewise geodesic paths adapted to
$\mathcal{P.}$ Under certain curvature restrictions on $M,$ it is
shown that
\[
\frac{1}{Z_{{\mathcal{P}}}}e^{-\frac{1}{2}E(\sigma)}dVol_{H_{{\mathcal{P}}}%
}(\sigma)\rightarrow\rho(\sigma)d\nu(\sigma)\text{ as }\mathrm{mesh}%
({\mathcal{P}})\rightarrow0,
\]
where $Z_{{\mathcal{P}}}$ is a \textquotedblleft
normalization\textquotedblright\ constant, $E:H\left(  M\right)
\rightarrow\lbrack0,\infty)$ is the energy functional, $Vol_{H_{{\mathcal{P}}%
}}$ is the Riemannian volume measure on $H_{\mathcal{P}}\left(
M\right)  ,$ $\nu$ is Wiener measure on continuous paths in $M,$
and $\rho$ is a certain density determined by the curvature tensor
of $M.$

\end{abstract}



\section{Introduction}\label{s.intro}

Suppose we have a Riemannian manifold $(M,\ g)$ of dimension $d$
with metric $g$. We will only consider $M$ to be compact or
$\mathbb{R}^{d}$. Fix a point $o$ on the manifold $M$ and let $V :
M \rightarrow\mathbb{R}$ be a potential function. In classical
mechanics, the path $\sigma: [0,T] \rightarrow M,\ \sigma(0) = o$,
subject to the potential $V$, can be obtained by solving Newton's
equation of motion
\[
\frac{\nabla}{dt} \sigma^{\prime}\left(  t\right)
=-\text{grad~}V\left( \sigma\left(  t\right)  \right)  ,
\]
where given any vector field $X(s)$ on $\sigma(s)$, define
\begin{equation}
\frac{\nabla X(s)}{ds} := //_{s}(\sigma)
\frac{d}{ds}\left\{//_{s}^{-1} (\sigma)
X(s) \right\}, \label{eq.coderive}%
\end{equation}
and $//_{s}(\sigma): T_{o}M \rightarrow T_{\sigma(s)}M$ is
parallel translation along $\sigma$ relative to the Levi Civita
covariant derivative $\nabla$. The Hamiltonian of the system, $H$
is then given by
\begin{equation}
H(\sigma(t), \sigma^{\prime}(t)) = \frac{1}{2}\parallel
\sigma^{\prime}(t)\parallel^{2} + V(\sigma(t)),\nonumber
\end{equation}
where $\parallel v\parallel^{2} := g(v,v)$ and mass is set to be
1.

In Quantum Mechanics, observables are no longer functions, but
rather Hermitian operators on some Hilbert space. Let $q = (q_{1},
q_{2}, \ldots q_{d})$ be the cartesian coordinates on
$\mathbb{R}^{d}$ and $p_{i}$ be the momentum corresponding to
$q_{i}$. In canonical quantization on $\mathbb{R}^{d}$, the
quantum mechanical operator $\widehat{H}$ corresponding to the
classical Hamiltonian, $H = \frac{1}{2}\sum_{i} p_{i}^{2} + V(q)$,
is given by
\begin{equation}
\widehat{H} = -\frac{\hbar^{2}}{2}\sum_{i}
\frac{\partial^{2}}{\partial q_{i}^{2}} + M_{V}\nonumber
\end{equation}
where $\hbar = \frac{h}{2\pi}$, $h$ is Planck's constant and
$M_{V}$ is multiplication operator by $V$.

However, on a manifold, one aims to quantize the Hamiltonian $H = \frac{1}%
{2}g^{ij}(q)p_{i}p_{j} + V(q)$ where $q$ is some coordinate system
on $M$. Using the \textquotedblleft Feynman's (Kac) path integral
prescription \textquotedblright, one defines the operator
$e^{-T\hat{H}}$ via an integral
\begin{align}
&  (e^{-T\widehat{H}}f)(o) :=
\text{\textquotedblleft}\frac{1}{Z_{T}}\int_{H_{T}(M)}e^{-\int_{0}^{T}
H(\sigma(t), \sigma^{\prime}(t)) \ dt}f(\sigma(T))\
\mathcal{D}\sigma
,\text{\textquotedblright} \label{eq.feynmann}%
\end{align}
where $H(\sigma(t), \sigma^{\prime}(t)) = \frac{1}{2}\parallel
\sigma^{\prime}(t)\parallel^{2} + V(\sigma(t))$ is the classical
Hamiltonian. $H_{T}(M)$ is the space of
finite energy paths, $Z_{T}$ is some normalization constant and $\mathcal{D}%
\sigma$ is to be interpreted as a \text{\textquotedblleft}
Lebesque type measure\text{\textquotedblright}. The operator
$\widehat{H}$ can then be obtained by differentiating the operator
$e^{-T\hat{H}}$ with respect to $T$ at 0.

For simplicity, we will assume that $T = 1$ and set $ \hbar = 1$.
Furthermore we will set $V = 0$ since quantizing a scalar function
on the manifold $M$ corresponds to the multiplication operator
with $V$.

\begin{defn}%
Define $H(M)$, a Hilbert manifold of absolutely continuous paths
with finite energy,
\begin{equation}
H(M) = \left\{\sigma:[0,1]\rightarrow M|\ \sigma(0)=o \in M\ {\rm
and}\ E(\sigma) < \infty \right\}
\label{defn.H(M)}%
\end{equation}
where the energy $E$ is given by
\begin{equation}
E(\sigma) := \int_{0}^{1} g(\sigma^{\prime}(s),\
\sigma^{\prime}(s))\ ds.
\label{eq.energy}%
\end{equation}

\end{defn}

The tangent space $T_{\sigma}H(M)$ to $H(M)$ at $\sigma$ may be
identified with the space of absolutely continuous vector fields
along $\sigma$. On this Hilbert manifold $H(M)$, we can define a
metric $G^{1}$ as follows.

\begin{defn}%
Let $T_{\sigma}H(M)$ be the space of absolutely continuous vector
fields $X$ along $\sigma$ {\rm (i.e. $X(s) \in T_{\sigma(s)}M\
\forall s \in [0,1]$)} such that $G^1(X, X) < \infty$ where
\begin{align}
G^{1}(X,\ X)  &  := \int_{0}^{1} g\left(\frac{\nabla X(s)}{ds},\
\frac{\nabla
X(s)}{ds}\right)\ ds, \label{eq.G^1}%
\end{align}
and $\frac{\nabla}{ds}$ is defined as in Equation
(\ref{eq.coderive}).
\end{defn}

See \cite{MR0098419, MR0226681, MR0341527, MR0478069, MR0158410}
for more details. By polarization, Equation (\ref{eq.G^1}) defines
a Riemannian metric on $H(M)$.

The integral over $H(M)$, defined in Equation (\ref{eq.feynmann})
however, is highly heuristic. Firstly, the normalization constant
$Z_{1}$ may be interpreted to be 0 or $\infty$. Secondly,
$\mathcal{D}\sigma$ which is to be interpreted as
\text{\textquotedblleft}Lebesgue measure\text{\textquotedblright},
fails to exist in an infinite dimensional space. 

We would like to make sense out of the RHS of Equation
(\ref{eq.feynmann}), by writing it as a limit of a sequence of
integrals over finite dimensional spaces $H_{ {\mathcal{P}} }(M)$.

\begin{defn}
\label{defn.H_p(M)} Let \beq {\mathcal{P}} = \{0 = s_{0} < s_{1} <
s_{2} < \cdots< s_{n} = 1\} \label{e.part} \eeq be a partition of
$[0,1]$. Define $H_{ {\mathcal{P}}}(M)$ as a set of piecewise
geodesics paths in $H(M)$ which change directions only at the
partition points.
\begin{equation}
H_{ {\mathcal{P}}}(M) = \left\{ \sigma\in H(M) \cap
C^{2}(I\backslash {\mathcal{P}})\ \Big|\
\frac{\nabla\sigma^{\prime}(s)}{ds} = 0\ for\ s
\notin{\mathcal{P}} \right\}.\nonumber
\end{equation}

\end{defn}

It will be shown later that $H_{ {\mathcal{P}}}(M)$ is a finite
dimensional submanifold of $H(M)$. In fact, $H_{
{\mathcal{P}}}(M)$ is diffeomorphic to $(\mathbb{R}^{d})^{n}$. For
$\sigma\in H_{ {\mathcal{P}}}(M)$, the tangent space
$T_{\sigma}H_{ {\mathcal{P}}}(M)$ can be identified with elements
$X \in T_{\sigma}H_{ {\mathcal{P}}}(M)$ satisfying the Jacobi
equations on $I\backslash{\mathcal{P}}$. As a submanifold of
$H(M)$, $H_{ {\mathcal{P}}}(M)$ inherits the induced metric
$G^{1}|_{TH_{ {\mathcal{P}}}(M)}$ by restricting the $G^{1}$
metric on $H_{ {\mathcal{P}}}(M)$.

If $N^{p}$ is any manifold with a metric $G$, define a volume
density $Vol_{G}$ on $T_nN$ by
\begin{equation}
Vol_{G}\left(v_{1}, v_{2}, \ldots, v_{p}\right) = \sqrt{\det\
\{G(v_{i}, v_{j})\}_{i,j =
1}^{p}} \label{eq.density}%
\end{equation}
where $\{v_{1}, v_{2}, \ldots, v_{p} \} \subset T_{n}N$ is a basis
and $n \in N$.

\begin{thm}\label{t.m_G}
Given a density of the form $\rho Vol_G$, where $\rho: N
\rightarrow [0, \infty)$, there exists a unique measure $m_G$ on
$N$ such that \beq \int_{\mathcal D(y)} f\ dm_G = \int_{\mathcal
D(y)} f \rho\ Vol_{G}\left(\frac{\partial}{\partial y_{1}},
\frac{\partial}{\partial y_{2}}, \ldots, \frac{\partial}{\partial
y_{p}}\right)dy_1 \ldots dy_p \nonumber \eeq for any local
coordinates $y = (y_1, \dots, y_p) : \mathcal{D}(y) \rightarrow
\mathbb{R}^p$ and measurable function $f: N \rightarrow
[0,\infty)$. If $\rho = 1$, the associated measure will be called
Riemann volume measure.
\end{thm}

\begin{defn}
Let $Vol_{ {\mathcal{P}}}$ denote the density on $H_{
{\mathcal{P}}}(M)$ determined by $G^{1}|_{TH_{ {\mathcal{P}}}(M)
\otimes TH_{ {\mathcal{P}}}(M)}$ using Equation
(\ref{eq.density}).
\end{defn}

Given the above definition, we can now define a measure on $H_{ {\mathcal{P}}%
}(M)$.

\begin{defn}\label{d.m}
For each partition $ {\mathcal{P}} $ of $[0,1]$ as in Equation
(\ref{e.part}), let $\nu_{ {\mathcal{P}}}$ denote the unique
measure on $H_{ {\mathcal{P}}}(M)$ as in Theorem \ref{t.m_G},
defined by the density
\begin{equation}
\frac{1}{Z_{ {\mathcal{P}}}^{1}}e^{-\frac{1}{2}%
E}Vol_{ {\mathcal{P}}} \nonumber
\end{equation}
where $E : H(M) \rightarrow[0,\infty)$ is the energy functional
defined in Equation (\ref{eq.energy}) and $Z_{ {\mathcal{P}}}^{1}$
is a normalization constant given by
\begin{equation}
Z_{ {\mathcal{P}}}^{1} = (2\pi)^{\frac{dn}{2}}. \label{eq.constant}%
\end{equation}

\end{defn}

We can now write the RHS of Equation (\ref{eq.feynmann}) as a
limit of a
sequence of integrals over the finite dimensional space $H_{ {\mathcal{P}}%
}(M)$, equipped with the measure $\nu_{ {\mathcal{P}}}$. Our
result shows that this limit can be written as an integral over
the Wiener space of $M$, with Wiener measure $\nu$. (See
Definition \ref{defn.wiener} below.)

\begin{defn}\label{d.J}
Define $\Delta_{i} s = s_{i} - s_{i-1}$ and $| {\mathcal{P}}| =
\bigvee_{i=1, \ldots, n}\Delta_{i}s = \max\{\Delta_{i} s : i =
1,2, \ldots, n \}$ be the norm of the partition and $J_{i} :=
(s_{i-1}, s_{i}]$ for $i = 1,2, \ldots, n$.
\end{defn}

Let $\triangle= {\mathrm{{tr}}} \ \nabla^{2}$ denote the Laplacian
acting on $C^{\infty}(M)$ and $p_{s}(x,y)$ be the fundamental
solution to the heat equation,
\begin{equation}
\frac{\partial u}{\partial s} = \frac{1}{2} \triangle u.\nonumber
\end{equation}

In the case when $M = \mathbb{R}^d$, \beq p_s(x, y) =
\left(\frac{1}{2\pi
s}\right)^{\frac{d}{2}}e^{-\frac{1}{2s}\parallel x-y\parallel^2}.
\nonumber \eeq

\begin{defn}
\label{defn.wiener} The Wiener space $W(M)$ is the path space
\begin{equation}
W(M) = \{\sigma:\ [0,1] \rightarrow M : \sigma(0) = o\ and\ \sigma
\ is\ continuous \}.\nonumber
\end{equation}

The $Wiener\ measure\ \nu$ associated to $(M,\ g,\ o)$ is the
unique probability measure on $W(M)$ such that
\begin{equation}
\int_{W(M)} f(\sigma)\ d\nu(\sigma) = \int_{M^{n}} F(x_{1},
\ldots, x_{n})\ \prod_{i=1}^{n}\ p_{\Delta_{i}s}(x_{i-1},\ x_{i})\
dm(x_{1}) \cdots
dm(x_{n}) \label{eq.wiener}%
\end{equation}
for all functions $f$ of the form $f(\sigma) =
F\left(\sigma(s_{1}), \ldots, \sigma(s_{n})\right)$, for all
${\mathcal{P}} $, a partition of $[0,1]$ as in Equation
(\ref{e.part}) and $F:\ M^{n} \rightarrow\mathbb{R}$ is a bounded
measurable function. In Equation $(\ref{eq.wiener})$, $dm(x)$
denotes the Riemann volume measure on $M$ as in Theorem
\ref{t.m_G} and by convention $x_{0} := o$.

\end{defn}

It is known that there exists a unique probability measure $\nu$
on $W(M)$ satisfying Equation (\ref{eq.wiener}). The measure $\nu$
is concentrated on continuous but nowhere differentiable paths.

\begin{notation}
When $M = \mathbb{R}^{d}$, $g( \cdot, \cdot)$ is the usual dot
product and $o = 0$, the measure $\nu$ defined in Definition
\ref{defn.wiener} is the standard Wiener measure on
$W(\mathbb{R}^{d})$. We will denote this standard Wiener measure
by $\mu$ rather than $\nu$. We will also let $b(s) :
W(\mathbb{R}^{d}) \rightarrow\mathbb{R}^{d}$ be the coordinate map
such that
\begin{equation}
b(s)(\omega) := \omega(s)\nonumber
\end{equation}
for all $\omega\in W(\mathbb{R}^{d})$.
\end{notation}

\begin{rem}
The process $\{b(s)\}_{s \in[0,1]}$ is a standard
$\mathbb{R}^{d}$-valued Brownian motion on the probability space
$(W(\mathbb{R}^{d}), \mu)$.
\end{rem}

Suppose we now view $M$ as an imbedded submanifold of
$\mathbb{R}^d$ with the induced Riemannian structure. Let $P(m)$
be the projection on the tangent space $T_mM$ and $v \in T_mM$.
Then for a vector field $X(m)$, \beq \nabla_v X =
P(\sigma(0))\frac{d}{dt}X(\sigma(t))\Big|_{t=0}\nonumber \eeq
where $\sigma$ is a path in $M$ such that $\sigma(0) = m$ and
$\sigma'(0) = v$. Define a projection $Q$ on the orthogonal
complement of $T_mM$ by $Q = I - P$, where $I$ is the identity.
With this definition, for any vector $v \in T_mM$, one can define
parallel translation along $\sigma$ by $//_s(\sigma) v := w(s) v$
where $w$ solves the following differential equation \beq w'(t) +
dQ(\sigma'(t))w(t) = 0,\ w(0) = P(\sigma(0)). \nonumber \eeq

\begin{thm}
Let $\Sigma$ be an $M$-valued semi-martingale and $V_0(m)$ be a
measurable vector field on $M$, then there is a unique parallel
$TM$-valued semi-martingale $V$ such that $V_0 = V_0(\Sigma_0)$
and $V_s \in T_{\Sigma_s}M$ for all $s$. More explicitly, $V_s =
w_s V(\Sigma_0)$ where $w$ solves the following Stratonovich
stochastic differential equation \beq \delta w + dQ(\delta
\Sigma)w = 0,\ w_0 = P(\Sigma_0). \label{e.s//trans} \eeq
\end{thm}

For a proof of this theorem, the reader should refer to
\cite{MR2090750}. Thus we can now define a
\text{\textquotedblleft}stochastic\text{\textquotedblright}
extension of parallel translation.

\begin{defn}\label{d.s//trans}{\rm (Stochastic Parallel Translation)}
Given $v \in T_{\Sigma_0}M$ and $M$-valued semi-martingale
$\Sigma$, define stochastic parallel translation $\widetilde{//}$
by \beq \widetilde{//}_s v := w_s v \nonumber \eeq where $w$
solves Equation (\ref{e.s//trans}). This is going to be used for
the particular semi-martingale $\Sigma_s(\sigma) := \sigma(s)$ on
$(W(M), \nu)$.
\end{defn}

\begin{defn}
The curvature tensor $R$ of $\nabla$ is
\begin{equation}
R(X,\ Y)Z = \nabla_{X} \nabla_{Y} Z - \nabla_{Y} \nabla_{X} Z -
\nabla _{[X,\ Y]} Z\nonumber
\end{equation}
for all vector fields $X,\ Y$ and $Z$ on $M$. The sectional
curvature $S(V)$ where $V \subseteq T_mM$ with $dim(V) = 2$, is
defined by \beq S(V) = \frac{g( R(u, v)u, v )}{\parallel u
\parallel^2
\parallel v
\parallel^2 - g(u,v)^2} , \nonumber \eeq where $\{u, v\}$ is a basis for
$V$. It can be shown that this definition is independent of the
basis used. Let $\{e_{i}\}_{i=1}^d \subseteq T_mM$ be an
orthornormal frame at $m \in M$.
The Ricci tensor of $(M,\ g)$ is $Ric\ v  = \sum_{i=1}^{d} R(v,\ e_{i}%
)e_{i}$, the scalar curvature Scal is $Scal = \sum_{i=1}^{d}
g(Ric\ e_{i}  ,\ e_{i})$ and $\Gamma_{m}$ is given by
\begin{equation}
\Gamma_{m} = \sum_{i,j = 1}^{d} \Big(R_{} (e_{i}, R_{}(e_{i},
\cdot)e_{j} )e_{j} + R_{} (e_{i}, R_{}(e_{j}, \cdot)e_{i} )e_{j} +
R_{} (e_{i}, R_{}(e_{j}, \cdot)e_{j} )e_{i} \Big).\nonumber
\end{equation}
Define for any $\sigma \in W(M)$, $K_{\sigma}: L^{2}([0,1]
\rightarrow T_{o}M) \longrightarrow L^{2}([0,1] \rightarrow
T_{o}M)$ by
\begin{equation}
(K_{\sigma}v)(s) = \int_{0}^{1} (\min{s}{t})\
\widetilde{//}_t^{-1}(\sigma)\left(\Gamma_{\sigma(t)}
\widetilde{//}_t(\sigma) v(t)\right)\ dt \nonumber
\end{equation}
where $\widetilde{//}$ is stochastic parallel translation.
\end{defn}

It will be proved later that $K_\sigma$ is a trace class operator.
We can now state the main result.

\begin{thm}
Let $(M,\ g)$ be a compact Riemannian manifold with dimension $d$.
Let ${\mathcal{P}} = \left\{0,\ \frac{1}{n},\ \frac {2}{n},\
\ldots,\ \frac{n-1}{n},\ 1 \right\}$ be an equally spaced
partition. Suppose that $f:\ W(M) \rightarrow\mathbb{R}$ is
bounded and continuous and that $0 \leq S < \frac{3}{17d}$ , then
\begin{align}
&{\lim_{| {\mathcal{P}}| \rightarrow0}} \int_{H_{
{\mathcal{P}}}(M)}
f(\sigma)\ d\nu_{ {\mathcal{P}} }(\sigma)  \nonumber \\
&  = \int_{W(M)} f(\sigma )e^{-\frac{1}{6}\int_{0}^{1}
Scal(\sigma(s))\ ds }\sqrt{\det\left(I + \frac
{1}{12}K_{\sigma}\right)}\ d\nu(\sigma). \label{eq.intro}%
\end{align}

\end{thm}

\subsection{Known Results}

Using $H_{ {\mathcal{P}}}(M)$ to approximate the Wiener space
$W(M)$ was done in \cite{MR1698956}. However, the choice of
metrics used on $TH_{ {\mathcal{P}}}(M)$ in \cite{MR1698956} are
different from $G^1|_{TH_\cP(M)}$.

\begin{defn}
\label{defn.G^1_P} For each partition $ {\mathcal{P}} $ of $[0,1]$
as in Equation (\ref{e.part}), let $G_{ {\mathcal{P}}}^{1}$, $G_{
{\mathcal{P}}}^{0}$ be the metrics on $TH_{ {\mathcal{P}}}(M)$,
given by
\begin{align}
G_{ {\mathcal{P}}}^{1}(X,Y) & = \sum_{i=1}^{n} g\left(  \
\frac{\nabla X(s_{i-1}+)}{ds},\ \frac{\nabla Y(s_{i-1}+)}{ds}\
\right)
\Delta_{i} s,\nonumber \\
G_{ {\mathcal{P}}}^{0}(X,Y) & = \sum_{i=1}^{n} g\left(  \
X(s_{i-1}+),\  Y(s_{i-1}+)\ \right)   \Delta_{i} s.\nonumber
\end{align}
for all $X, Y \in T_{\sigma}H_{ {\mathcal{P}}}(M)$ and $\sigma\in
H_{ {\mathcal{P}}}(M)$. Note that $\frac{\nabla X(s_{i-1}+)}{ds} =
\lim_{s \rightarrow{s_{i-1}}{+}}\ \frac{\nabla X(s)}{ds} $.
\end{defn}

Observe that $G_{ {\mathcal{P}}}^{1}$ is some sort of Riemann sum
approximation to $G^{1}$.

\begin{defn}
For each partition $ {\mathcal{P}} $ of $[0,1]$ as in Equation
(\ref{e.part}), define unique measures $\nu_{G_{
{\mathcal{P}}}^{1}}$ and $\nu_{G_{ {\mathcal{P}}}^{0}}$ on $H_{
{\mathcal{P}} }(M)$, as in Theorem \ref{t.m_G}, defined by
densities
\begin{align}
&\frac{1}{Z_{ {\mathcal{P}}}^{1}}%
e^{-\frac{1}{2}E}Vol_{G_{
{\mathcal{P}}}^{1}}\ {\rm and}\ \frac{1}{Z_{ {\mathcal{P}}}^{0}}%
e^{-\frac{1}{2}E}Vol_{G_{ {\mathcal{P}}}^{0}}\nonumber
\end{align}
respectively, where $E : H(M) \rightarrow[0,\infty)$ is the energy
functional defined in Equation (\ref{eq.energy}) and $Vol_{G_{
{\mathcal{P}}}^{1}}$ and $Vol_{G_{ {\mathcal{P}}}^{0}}$ are
determined by $G_{ {\mathcal{P}}}^{1}$ and $G_{
{\mathcal{P}}}^{0}$ respectively using Equation (\ref{eq.density})
. $Z_{ {\mathcal{P}}}^{1}$ is a normalization constant given by
Equation (\ref{eq.constant}) and \beq Z^0_\mathcal{P} =
\prod_{i=1}^n \left(\sqrt{2\pi}\left(s_i -
s_{i-1}\right)\right)^d.\nonumber \eeq
\end{defn}

The following theorem was proved in \cite{MR1698956}.

\begin{thm}
Let $M$ be a compact Riemannian manifold. Suppose that $f:\ W(M)
\rightarrow\mathbb{R}$ is bounded and continuous, then
\begin{align}
{\lim_{| {\mathcal{P}}| \rightarrow0}} \int_{H_{
{\mathcal{P}}}(M)} f(\sigma)\ d\nu_{G_{
{\mathcal{P}}}^{1}}(\sigma) & = \int_{W(M)}
f(\sigma)\ d\nu(\sigma)\nonumber\ \text{and} \\
\lim_{|\mathcal{P}| \rightarrow 0} \int_{H_\mathcal{P} (M)}
f(\sigma)\ d\nu_{G^0_\mathcal{P}}(\sigma) & =  \int_{W(M)}
f(\sigma)e^{-\frac{1}{6}\int_0^1 Scal(\sigma(s))ds}\ d\nu(\sigma).
\nonumber
\end{align}

\end{thm}

Using the Feynman-Kac formula, $\widehat{H}$ in Equation
(\ref{eq.feynmann}) is given by
\begin{equation}
\widehat{H} = -\frac{1}{2}\Delta + \kappa Scal \nonumber
\end{equation}
where $\Delta$ is the Laplacian operator and $\kappa = 0,
\frac{1}{6}$ for the $G_{ {\mathcal{P}}}^{1}$ and $G_{
{\mathcal{P}}}^{0}$ metrics respectively. However, the integral in
Equation (\ref{eq.intro}) is not of the Feynman-Kac form. Hence
the interpretation of the operator $\widehat{H}$ corresponding to
this integral is unclear at this stage.

\section{Finite Dimensional Approximations}\label{s.fda}

A detailed account of this section is given in \cite{MR1698956}.

Let $\pi: O(M) \rightarrow M$ denote the bundle of orthogonal
frames on $M$. An element $u \in O(M)$ is an isometry
\begin{equation}
u: \mathbb{R}^{d} \rightarrow T_{\pi(u)}M.\nonumber
\end{equation}
Fix $u_{o} \in\pi^{-1}(o)$, which identifies $T_{o}M$ of $M$ at
$o$ with $\mathbb{R}^{d}$.

Define $\theta$, a $\mathbb{R}^{d}$-valued form on $O(M)$ by
$\theta_{u}(\xi) = u^{-1}\pi_{*} \xi$ for all $u \in O(M)$,
$\xi\in T_{u}O(M)$ and let $\vartheta$ be the $so(d)$-valued
connection form on $O(M)$ defined by $\nabla$. Explicitly, if $s
\rightarrow u(s)$ is a smooth path in $O(M)$ then
\begin{equation}
\vartheta(u^{\prime}(0)) = u(0)^{-1}\frac{\nabla
u(s)}{ds}\Big|_{s=0}\nonumber
\end{equation}
where $\frac{\nabla}{ds}$ is defined as in Equation
(\ref{eq.coderive}) with $X$ replaced by $u$. We define the
horizontal lift ${\mathcal{H}}_{u}: T_{\pi(u)}M \rightarrow
T_{u}O(M)$ by
\begin{equation}
\theta{\mathcal{H}}_{u}u = id_{\mathbb{R}^{d}},\ \vartheta_{u}
{\mathcal{H}}_{u} = 0.\nonumber
\end{equation}
Explicitly, for $v \in T_{\pi(u)}M$, $\mathcal{H}_uv =
\frac{d}{dt}|_{t=0}//_t(\sigma)u$ where $\dot{\sigma}(0) = v$.

\begin{defn}\label{defn.Omegau}
For $a,\ c \in\mathbb{R}^{d}$, let
\begin{equation}
\Omega_{u}(a, c) := u^{-1}R(ua, uc)u.\nonumber
\end{equation}

\end{defn}

Let $H(O(M))$ be the set of finite energy paths $u: [0,1]
\rightarrow O(M)$ as defined in Equation (\ref{defn.H(M)}) with
$M$ replaced by $O(M)$ and $o$ by $u_{o}$. For $\sigma\in H(M)$,
let $u$ be defined by the ordinary differential equation
\begin{equation}
u^{\prime}(s) = {\mathcal{H}}_{u(s)}\sigma^{\prime}(s),\ u(0) = u_{o}%
.\nonumber
\end{equation}
This equation implies that $\vartheta(u^{\prime}(s)) = 0$ or that
$\frac{\nabla u(s)}{ds} = 0$. Thus we have
\begin{equation}
u(s) = //_{s}(\sigma)u_{o}\nonumber
\end{equation}
where $//_{s}(\sigma)$ is parallel translation along $\sigma$.
Since $u_{o}$ is fixed, we will drop it and write $u(s) = //_{s}$.
We will call $u(s)$ a horizontal lift of $\sigma$ starting at
$u_{o}$ and use it to define $\phi$, which associates $\omega \in
H(\mathbb{R}^{d})$ with a path $\sigma\in H(M)$.

\begin{defn}
{\rm (Cartan's Development Map)} The development map, $\phi:
H(\mathbb{R}^{d}) \rightarrow H(M)$ is defined for $\omega \in
H(\mathbb{R}^{d})$ by $\phi(\omega) = \sigma\in H(M)$ where
$\sigma$ solves the functional differential equation
\begin{equation}
\sigma^{\prime}(s) = //_{s}(\sigma)\omega^{\prime}(s),\ \sigma(0)
= o.
\label{eq.cartan}%
\end{equation}

\end{defn}

The development map, $\phi$ is smooth and injective. We can define
an anti-development map, $\phi^{-1}: H(M) \rightarrow
H(\mathbb{R}^{d})$ given by $\omega = \phi^{-1}(\sigma)$ where
\begin{equation}
\omega(s) = \int_{0}^{s} //_{r}^{-1}(\sigma)\sigma^{\prime}(r)\
dr.\nonumber
\end{equation}

Again, $\phi^{-1}$ is smooth and is injective. Therefore, $\phi:
H(\mathbb{R}^{d}) \rightarrow H(M)$ is a diffeomorphism of
infinite dimensional Hilbert manifolds.

\begin{defn}
For each $h \in C^{\infty}( H(M) \rightarrow H(\mathbb{R}^{d}) )$
and $\sigma\in H(M)$, let $X^{h}(\sigma) \in T_{\sigma}H(M)$ be
given by
\begin{equation}
X_{s}^{h}(\sigma) := //_{s}(\sigma)h_{s}(\sigma) \label{eq.X^h}%
\end{equation}
for all $s \in [0,1]$, where we have written $h_{s}(\sigma)$ as
$h(\sigma)(s)$.
\end{defn}

Define $H_{ {\mathcal{P}}}(\mathbb{R}^{d}) = \{ \omega \in H \cap
C^{2}(I \backslash{\mathcal{P}})\ |\ \omega^{\prime\prime}(s) = 0\
for\ s \notin {\mathcal{P}} \}$, the set of piecewise linear paths
in $H(\mathbb{R}^{d})$, which changes directions only at the
partition points.

\begin{rem}
The development map $\phi: H(\mathbb{R}^{d}) \rightarrow H(M)$ has
the property that $\phi\left(H_{ {\mathcal{P}}}(\mathbb{R}^{d})
\right) = H_{ {\mathcal{P}}}(M)$ where $H_{ {\mathcal{P}}}(M)$ has
been defined in Definition \ref{defn.H_p(M)}. If $\sigma=
\phi(\omega)$ with $\omega \in H_{
{\mathcal{P}}}(\mathbb{R}^{d})$, then differentiating Equation
(\ref{eq.cartan}) gives
\begin{equation}
\frac{\nabla\sigma^{\prime}(s)}{ds} = \frac{\nabla}{ds}\Big(//_{s}%
(\sigma)\omega^{\prime}(s) \Big) =
//_{s}(\sigma)\omega^{\prime\prime}(s) = 0\ for\ all\ s
\notin{\mathcal{P}}.\nonumber
\end{equation}

\end{rem}

Because $\phi: H(\mathbb{R}^{d}) \rightarrow H(M)$ is a
diffeomorphism and $H_{ {\mathcal{P}}}(\mathbb{R}^{d}) \subset
H(\mathbb{R}^{d})$ is an embedded submanifold, so it follows that
$H_{ {\mathcal{P}}}(M)$ is an embedded submanifold of $H(M)$.
Therefore for each $\sigma\in H_{ {\mathcal{P}}}(M)$,
$T_{\sigma}H_{ {\mathcal{P}}}(M)$ may be viewed as a subspace of
$T_{\sigma }H(M)$. The next proposition identifies this subspace.
See \cite{MR1698956} for a proof.

\begin{prop}
\label{prop.jacobi} Let $\sigma\in H_{ {\mathcal{P}}}(M)$, then $X
\in T_{\sigma}H(M)$ is in $T_{\sigma}H_{ {\mathcal{P}}}(M)$ if and
only if
\begin{equation}
\frac{\nabla^{2}}{ds^{2}}X(s) = R(\sigma^{\prime}(s),\ X(s)
)\sigma
^{\prime}(s) \label{eq.JacobiEquation}%
\end{equation}
on $I \backslash{\mathcal{P}}$. Equivalently, letting $\omega = \phi^{-1}%
(\sigma),\ u = //(\sigma)$ and $h \in H(\mathbb{R}^{d})$, then
$X^{h} \in T_{\sigma}H(M)$ defined in Equation (\ref{eq.X^h}) is
in $T_{\sigma}H_{ {\mathcal{P}} }(M)$ if and only if
\begin{equation}
h^{\prime\prime}(s) = \Omega_{u(s)}(\omega^{\prime}(s),
h(s))\omega^{\prime}(s)
\label{eq.jacobieq}%
\end{equation}
on $I \backslash{\mathcal{P}}$.
\end{prop}

\subsection{Comparing $\nu_{\mathcal{P}}$ and
$\nu_{G_{\mathcal{P}}^{1}}$ on $M$\label{s.compare}}

\begin{defn}
For $\omega \in H_\cP(\mathbb{R}^d)$, let $\{h_{k,a}\}_{%
\genfrac{}{}{0pt}{}{k=1,2,\ldots,n}{a=1,2,\ldots,d}%
}$ be any basis in
$\phi_{\ast}^{-1}(T_{\phi(\omega)}H_{\mathcal{P}}%
(M)) \\
= H_{\mathcal{P}}\left(  T_{o}M\right)  $ and
$\sigma=\phi(\omega)$. Let
$Vol_{\mathcal{P}}$ be the density associated to $G^{1}%
|_{TH_{\mathcal{P}}(M)\otimes TH_{\mathcal{P}}(M)}$ metric and
$Vol_{G_{\mathcal{P}}^{1}}$ be the density associated to $G_{\mathcal{P}%
}^{1}$ metric. Then
\[
X_{s}^{h_{k,a}}\left(  \sigma\right)  :=//_{s}\left( \sigma\right)
h_{k,a}\left(  s\right)  \text{ for }k=1,2,\dots,n\text{ and
}a=1,2,\dots,d
\]
is a basis for $T_{\sigma}H_{\mathcal{P}}(M)$ and we define%
\begin{align}
\rho_{\mathcal{P}}&=\frac{\left\vert Vol_{\mathcal{P}}\left(
\{X^{h_{k,a}}\}_{%
\genfrac{}{}{0pt}{}{i=1,2,\ldots,n}{a=1,2,\ldots,d}%
}\right)  \right\vert }{\left\vert Vol_{G_{\mathcal{P}}^{1}}\left(
\{X^{h_{k,a}}\}_{%
\genfrac{}{}{0pt}{}{i=1,2,\ldots,n}{a=1,2,\ldots,d}%
}\right)  \right\vert }\nonumber \\
&=\frac{\sqrt{\det\left(  \left\{ G^{1}\left(
X^{h_{k,a}},X^{h_{k^{\prime},a^{\prime}}}\right) \right\} _{\left(
k,a\right)  ,\left(  k^{\prime},a^{\prime}\right) }\right)
}}{\sqrt {\det\left(  \left\{  G_{\mathcal{P}}^{1}\left(
X^{h_{k,a}},X^{h_{k^{\prime },a^{\prime}}}\right)  \right\}
_{\left( k,a\right)  ,\left(  k^{\prime
},a^{\prime}\right)  }\right)  }}. \label{eq.dense}%
\end{align}

\end{defn}

The relevance of this definition is contained in the next remark.

\begin{rem}
\label{rem.dense}First off, it is well known (and easily verified)
that the $\rho_{\mathcal{P}}\left(  \sigma\right)  $ defined in
Equation (\ref{eq.dense}) is well defined independent of the
choice of basis
$\{h_{k,a}\}_{%
\genfrac{}{}{0pt}{}{k=1,2,\ldots,n}{a=1,2,\ldots,d}%
}.$ Secondly, if $\nu_{\mathcal{P}}$ and
$\nu_{G_{\mathcal{P}}^{1}}$ are the measures associated to
$G^{1}|_{TH_{\mathcal{P}}(M)\otimes TH_{\mathcal{P}}(M)}$ and
$G_{\mathcal{P}}^{1}$ respectively, then
$d\nu_{\mathcal{P}}=\rho_{{\mathcal{P}}}\cdot d\nu_{G_{\mathcal{P}}^{1}%
}.$
\end{rem}

From \cite{MR1698956}, we know the limiting behavior of the
measure $\nu_{G_{\mathcal{P}}^{1}}.$ Hence, our proof that
$\nu_{\mathcal{P}}$ has a limit will break into two parts. Very
roughly speaking we are going to first show that $\left\{
\rho_{\mathcal{P}}:\mathcal{P}\right\}  $ is uniformly integrable
and then we will show that $\lim_{\left\vert
\mathcal{P}\right\vert \rightarrow0}\rho_{\mathcal{P}}$ exists in
$\mu$-measure. This rough outline will
have to be appropriately modified since $\left\{  \rho_{\mathcal{P}%
}:\mathcal{P}\right\}  $ are functions on different probability
spaces for each $\mathcal{P}.$ This will be remedied by pulling
$\rho_{\mathcal{P}}$ to classical Wiener space $\left(  W\left(
T_{o}M\right)  ,\mu\right)  $ using Cartan's rolling map $\phi$
and the natural piecewise linear approximation map from $W\left(
T_{o}M\right)  $ to $H_{\mathcal{P}}\left(  T_{o}M\right)$. We
will identify $T_o(M)$ with $\mathbb{R}^d$.

\par Before we move on, we would like to point out that
$\rho_\cP \circ \phi$ is only defined on
$H_\mathcal{P}(\mathbb{R}^d)$, which has $\mu$-measure zero.

\begin{defn}
\label{defn.projection} Let $\{b(s)\}_{s\in\lbrack0,1]}$ be the
standard $\mathbb{R}^{d}$ Brownian motion on $(W(\mathbb{R}^{d}),\
\mu)$ and ${\mathcal{P}}$ be any partition, i.e. $b(s) :
W(\mathbb{R}^d) \rightarrow \mathbb{R}^d$, \beq b(s)(\omega) :=
b(s, \omega) := \omega(s). \nonumber \eeq By abuse of notation,
define $b_\cP : W(\mathbb{R}^d) \rightarrow H_\cP(\mathbb{R}^d)$
by
\begin{align}
b_\cP(s)&=b(s_{i-1})+(s-s_{i-1})\frac{\Delta_{i}b}{\Delta _{i}s}\
\ if\ s\in(s_{i-1},\ s_{i}]\nonumber
\end{align}
where $\Delta_{i}b:=b(s_{i})-b(s_{i-1})$. We will write $b_n =
b_{\cP_n}$ if $\mathcal{P}_n = \{0 < \frac{1}{n} < \cdots
<\frac{n}{n}=1 \}$ is an equally spaced partition
\end{defn}

Thus by composing with $b_\cP$, we can now view $\rho_\cP \circ
\phi \circ b_\cP$ as a random variable on $(W(\mathbb{R}^d),
\mu)$.

\section{Uniform Integrability of $\{\rho_{ n} \circ \phi \circ b_n\}_{n=1}^\infty%
$}\label{s.ui}

\par We will first show that $\rho_\cP
\circ \phi \circ b_\cP$ is uniformly integrable. But first, we
need to write down a formula for $\rho_\cP$.

\subsection{A First Formula for $\rho_n$}\label{s.uid}

We will now only consider equally spaced partitions ${\mathcal{P}_n}%
=\{0=s_{0}<s_{1}<s_{2}<\dots<s_{n}=1\}$, such that $\Delta_{i}s=\frac{1}%
{n},\ i=1,\dots n$ and write $\rho_n = \rho_{\cP_n}$. Let
$\sigma\in H_{\mathcal{P}_n}(M)$ and consider
$\omega=\phi^{-1}(\sigma)$. On each $J_{i}=(s_{i-1}, s_i],\
i=1,2,\ldots,n$, $\omega_{i}^{\prime }:=\omega^{\prime}(s_{i-1}+)$
is a constant. Thus $\Delta_{i}\omega:=\omega_{i}^{\prime
}\Delta_{i}s$.

From Proposition \ref{prop.jacobi}, we know that for $s\in J_{i}$,
for each $h$ such that $X^h \in TH_{\cP_n}(M)$, $h(\omega, s)$
satisfies the ordinary differential equation
\begin{equation}
\frac{d^{2}h(\omega,
s)}{ds^{2}}=\Omega_{u(s)}(\omega_{i}^{\prime},h\left(\omega,
s\right)
)\omega_{i}^{\prime}. \label{eq.Jacobi}%
\end{equation}

Let $\left\{  e_{a}\right\}  _{a=1}^{d}$ be the standard basis for
$\mathbb{R}^{d}$ and for $i=1,2,\dots,n,$ let
\[
e_{i,a}=\left(  0,\dots,0,\overset{\text{i}^{\text{th }}\text{--spot}}{e_{a}%
},0,\dots,0\right)  \in\left(  \mathbb{R}^{d}\right)
^{n}=\mathbb{R}^{nd}.
\]
Then $\{e_{i,a}\}_{%
\genfrac{}{}{0pt}{}{{i=1,\ldots n}}{a=1,\ldots,d}%
}$ is an indexing of the standard basis for $\mathbb{R}^{nd}$ such
that all the components of $e_{i,a}$ are $0$ except at the
$a+(i-1)d$ position, which is $1$.

\begin{notation}
\label{not.a}Let $h_{i,a}(\omega, s)$ denote the continuous function in $\mathbb{R}%
^{d}$ which solves Equation (\ref{eq.Jacobi}) on $\left[
0,1\right]
\setminus\mathcal{P}_n$ and satisfies%
\begin{equation}
h_{i,a}(\omega, 0)=0,\text{ and }h_{i,a}^{\prime}(\omega,
s_{j-1}+)=\delta_{ij}e_{a}\text{ for }j=1,\ldots,n.\nonumber
\end{equation}
It is easily seen that $\{h_{i,a}\}_{%
\genfrac{}{}{0pt}{}{i=1,\ldots,n}{a=1,\ldots,d}%
}$ forms a basis for $\phi_{\ast}^{-1}(T_{\phi(\omega)}H_{\mathcal{P}_n%
}(M))=TH_{{\mathcal{P}_n}}(\mathbb{R}^{d})\cong
H_{\mathcal{P}_n}\left( \mathbb{R}^{d}\right)  .$ Further let
$\mathcal{Q}^{n}$ denote the $nd\times nd$ matrix which is given
in $d\times d$ blocks, $\mathcal{Q}^{n}:=\left\{ \left(
\mathcal{Q}_{mk}^{n}\right) \right\}  _{m,k=1}^{n},$ with
\[
\left(  \mathcal{Q}_{mk}^{n}e_{a},e_{c}\right)(\omega)
:=\int_{0}^{1}\left\langle h_{ma}^{\prime}\left(\omega,  s\right)
,h_{kc}^{\prime}\left(\omega,  s\right) \right\rangle ds\text{ for
}a,c=1,2,\dots,d.
\]

\end{notation}


\begin{notation}
Unless stated otherwise, upper case letters without a superscript
$n$ will denote $d\times d$ matrices. Upper case letters and
scripted upper case letters with a superscript $n$ will denote
$n\times n$ block matrices with entries being $d\times d$ blocks.
We will reserve ${\mathcal{I}^{n}}$ and $I$ for a $nd \times nd$
identity matrix and a $d \times d$ identity matrix respectively.
Matrices with a superscript $T$ will denote the matrix transpose.
To avoid confusion, we will use $\mathrm{Tr}$ and $\mathrm{tr}$ to
denote taking the trace of a $nd \times nd$ matrix and a $d \times
d$ matrix respectively.\newline $b$ will denote a standard $d$ --
dimensional Brownian path in $\mathbb{R}^{d}$. For a piecewise
continuous
function on $[0,1]$, we will use the notation $f(s+)=\lim_{r\downarrow s%
}f(r)$. We also let $\left\langle \cdot,\cdot\right\rangle $
denote $g(\cdot,\cdot)|_{o}$ at base point $o$.
\end{notation}


\begin{rem}
All norms used for matrices will be the operator norm. Norms used
for vectors will be the euclidean norm.
\end{rem}

\begin{lem}
\label{lem.firstformula}The relationship between $\rho_n \circ \phi  $ and $\mathcal{Q}^{n}$ is%
\begin{equation}
\rho_{n} \circ \phi=\sqrt{\det\ (n\mathcal{Q}^{n})}. \label{eq.rhoq}%
\end{equation}

\end{lem}

\begin{proof}
Observe that
\begin{align}
G^{1}\left(  X^{h_{k,a}},\ X^{h_{m,c}}\right)(\omega)   &
=\int_{0}^{1}g\left( \frac{\nabla X^{h_{k,a}}(\omega, s)}{ds},\
\frac{\nabla X^{h_{m,c}}(\omega, s)}{ds}\right)
\ ds\nonumber\\
&  =\int_{0}^{1}\left\langle h_{k,a}^{\prime}(\omega, s),\
h_{m,c}^{\prime
}(\omega, s)\right\rangle \ ds\nonumber \\
&=\langle\mathcal{Q}_{mk}^{n}(\omega)\ e_{a},\ e_{c}%
\rangle=\langle\mathcal{Q}^{n}(\omega)\ e_{k,a},\ e_{m,c}\rangle
\nonumber
\end{align}
Hence
\begin{align*}
\left\vert Vol_{\mathcal{P}}\left(  \{X^{h_{k,a}}:\ k=1,2,\ldots
,n,a=1,2,\ldots,d\}\right)  \right\vert  &  :=\sqrt{\det\left\{
G^{1}\left(
X^{h_{k,a}},X^{h_{m,c}}\right)  \right\}  }\\
&  =\sqrt{\det\ \mathcal{Q}^{n}}.
\end{align*}
On the other hand,%
\begin{align*}
G_{{\mathcal{P}}}^{1}\left(
X^{h_{k,a}},X^{h_{m,c}}\right)(\omega) & =\sum _{i=1}^{n}g\left(
\frac{\nabla X^{h_{k,a}}(\omega, s_{i-1}+)}{ds},\ \frac{\nabla
X^{h_{m,c}}(\omega, s_{i-1}+)}{ds}\right)  \Delta_{i}s\\
&  =\sum_{i=1}^{n}\langle h_{k,a}^{\prime}(\omega, s_{i-1}+),h_{m,c}^{\prime}%
(\omega, s_{i-1}+)\rangle\ \Delta_{i}s\\
&  =\sum_{i=1}^{n}\delta_{k,i}\delta_{m,i}\langle
e_{a},e_{c}\rangle \
\Delta_{i}s=\delta_{k,m}\delta_{a,c}\Delta_{k}s,
\end{align*}
i.e. $\left\{  G_{{\mathcal{P}}}^{1}\left(
X^{h_{k,a}},X^{h_{m,c}}\right)
\right\}  =\frac{1}{n}\mathcal{I}^n$ and hence%
\begin{align}
\left\vert Vol_{\mathcal{P}}\left(  \{X^{h_{k,a}}:\ k=1,2,\ldots
,n,a=1,2,\ldots,d\}\right)  \right\vert  &  :=\sqrt{\det\left\{G_{{\mathcal{P}}}%
^{1}\left(X^{h_{k,a}},X^{h_{m,c}}\right)\right\}}\nonumber\\
&  =\sqrt{\prod_{i=1}^{n}\left(  \Delta_{i}s\right)
^{d}}=n^{-dn/2}.\nonumber
\end{align}
Hence it follows that
\begin{equation}
\rho_{n} \circ \phi=\frac{\sqrt{\det\ \mathcal{Q}^{n}}}{\sqrt{\prod_{i=1}%
^{n}\left(  \Delta_{i}s\right)  ^{d}}}\ =\sqrt{\det(n\mathcal{Q}^{n}%
)}.\nonumber
\end{equation}

\end{proof}

\par Our next goal is to estimate
the size of $\det(n\mathcal{Q}^{n}).$ For this we are going to
apply Proposition \ref{prop.linalgebra1} in the Appendix as
follows. For $\alpha\geq1,$ to be chosen later, we have, from
Equation (\ref{eq.2}) with
$M=n\mathcal{Q}^{n},$ that%
\begin{align}
\det(n\mathcal{Q}^{n})  &  \leq\alpha^{nd}\exp\left(  \alpha^{-1}%
\operatorname{Tr}\left(  n\mathcal{Q}^{n}-
{\mathcal{I}^{n}}\right) \right) =\alpha^{nd}\exp\left(
\alpha^{-1}\sum_{m=1}^{n}\operatorname{tr}\left(
n\mathcal{Q}_{mm}^{n}-I\right)  \right) \nonumber\\
&  \leq\alpha^{nd}\exp\left(  \alpha^{-1}d\sum_{m=1}^{n}\left\Vert
n\mathcal{Q}_{mm}^{n}-I\right\Vert \right)  , \label{eq.est}%
\end{align}
where we have used Equation (\ref{eq.b}) of Proposition
\ref{prop.linalgebra2} in the last inequality. So according to
Equation
(\ref{eq.est}) we need to estimate $\left\Vert n\mathcal{Q}_{mm}%
^{n}-I\right\Vert $ for each $m.$ The first step in the process is
to record a formula for $\mathcal{Q}_{mm}^{n}.$

\subsection{A Formula for $\mathcal{Q}^{n}$\label{s.formq}}

\begin{notation}
\label{not.brack}Given any partition, $\mathcal{P}:=\left\{  0=s_{0}<s_{1}%
<\dots<s_{n}=1\right\}  $ of $\left[  0,1\right]  ,$ $j\in\left\{
1,2,\dots,n\right\}  ,$ and $s\in\left[  0,1\right]  ,$ let
\[
\left[  s\right]  _{j}:=\left[  \left(  s-s_{j-1}\right)
\vee0\right] \wedge\Delta_{j}s=\left\{
\begin{array}
[c]{ccc}%
0 & \text{if} & s\leq s_{j-1}\\
s-s_{j-1} & \text{if} & s\in J_{j}\\
\Delta_{j}s & \text{if} & s\geq s_{j}%
\end{array}
\right.  .
\]

\end{notation}

\begin{defn}
\label{defn.C&S} For $i=1,2,\dots,n$ and $0\leq s\leq\Delta_{i}s,$ let%
\begin{equation}
A_{i}(\omega,s):=\Omega_{u(s_{i-1}+s)}(\omega_{i}^{\prime},\cdot)\omega_{i}^{\prime}
\label{eq.ai}%
\end{equation}
and $C_{i}\left(  \omega,s\right)  $ and $S_{i}\left(
\omega,s\right) $ be the $\operatorname*{End}\left(
\mathbb{R}^{d}\right)  $ -- valued functions
determined by%
\[
\frac{d^{2}C_{i}\left(  \omega,s\right)
}{ds^{2}}=A_{i}(\omega,s)C_{i}\left( \omega,s\right)  \text{ with
}C_{i}\left(  \omega,0\right)  =I\text{ and }C_{i}^{\prime }\left(
\omega,0\right)  =0
\]
and%
\[
\frac{d^{2}S_{i}\left(  \omega,s\right)
}{ds^{2}}=A_{i}(\omega,s)S_{i}\left( \omega,s\right)  \text{ with
}S_{i}\left(  \omega,0\right)  =0\text{ and }S_{i}^{\prime }\left(
\omega,0\right)  =I.
\]
Let $K > 0$ be a constant such that
\begin{align}
\parallel A_{i}(\omega, s) \Delta_{i}s^{2} \parallel &  = \parallel\Omega
_{u(s_{i-1}+s)}(\omega_{i}^{\prime}, \cdot)\omega_{i}^{\prime}\parallel\Delta_{i}%
s^{2}\nonumber\\
&  \leq K \parallel\Delta_{i}\omega \parallel^{2} .\nonumber
\end{align}
Refer to Section \ref{s.bds} of Appendix for the existence of such
a $K$.
\end{defn}

Using this notation it follows that $h(\omega, \cdot)$ is a
continuous function which solves
Equation (\ref{eq.Jacobi}) on $\left[  0,1\right]  \setminus\mathcal{P}_n,$ then%
\begin{equation}
h(\omega, s)=C_{i}(\omega,s-s_{i-1})\ h(\omega, s_{i-1})+S_{i}(\omega,s-s_{i-1})\ h^{\prime}%
(\omega, s_{i-1}+)\text{ when }s\in J_{i}.\nonumber
\end{equation}

\begin{notation}
\label{not.prod}For $m,l\in\left\{  1,2,\dots,n\right\}  $ and
$d\times d$ matrices, $\left\{  M_{k}\right\}  _{k=1}^{n},$ let
\[
\prod_{k=l}^{m}M_{k}:=\left\{
\begin{array}
[c]{ccc}%
I & \text{if} & m<l\\
M_{m}M_{m-1}\dots M_{l+1}M_{l} & \text{if} & m\geq l
\end{array}
.\right.
\]

\end{notation}

With all of this notation, we may write $h_{i,a}\left(\omega,
s\right) $ as in the following lemma.

\begin{lem}
\label{lem.hia}Continuing to use the notation introduced above, we have%
\begin{align}
h_{m,a}(\omega, s)  &  =\left[  \prod_{k=m+1}^{n}C_{k}\left(
\omega,\left[ s\right] _{k}\right)  \right]  S_{m}\left(
\omega,\left[ s\right] _{m}\right)
e_{a}\label{eq.hia1}\\
&  =\left\{
\begin{array}
[c]{lll}%
0 & \text{if} & s\leq s_{m-1}\\
S_{m}(\omega,s-s_{m-1})e_{a} & \text{if} & s\in J_{m}\\
C_{j}(\omega,s-s_{j-1})V_{mj}\left(  \omega\right)  e_{a} &
\text{if} & s\in J_{j}\text{ and }j\geq m+1
\end{array}
\right.  \label{eq.hia}%
\end{align}
where
\begin{equation}
V_{mj}(\omega):=\left[  \prod_{k=m+1}^{j-1}C_{k}(\omega,\Delta_{k}s)\right]  S_{m}%
(\omega,\Delta_{m}s). \label{eq.vij}%
\end{equation}
Differentiating Equation (\ref{eq.hia}),%
\begin{equation}
h_{m,a}^{\prime}\left(\omega,  s\right)  =1_{J_{m}}\left(
s\right) S_{m}^{\prime }\left(  \omega,s-s_{m-1}\right)
e_{a}+\sum_{j=m+1}^{n}1_{J_{j}}\left(  s\right)
C_{j}^{\prime}(\omega,s-s_{j-1})V_{mj}e_{a}. \label{eq.hip}%
\end{equation}

\end{lem}

From Equation (\ref{eq.hip}), we learn that
\begin{align*}
\left\langle \mathcal{Q}_{mm}^{n}(\omega)e_{a},e_{c}\right\rangle  &  =\int_{J_{m}%
}\left\langle S_{m}^{\prime}\left(  \omega,s-s_{m-1}\right)
e_{a},S_{m}^{\prime
}\left(  \omega,s-s_{m-1}\right)  e_{c}\right\rangle ds\\
&  +\sum_{j=m+1}^{n}\int_{J_{j}}\left\langle C_{j}^{\prime}(\omega,s-s_{j-1}%
)V_{mj}(\omega)e_{a},C_{j}^{\prime}(\omega,s-s_{j-1})V_{mj}(\omega)e_{c}\right\rangle ds\\
&  =\left\langle \int_{0}^{\Delta_{m}s}S_{m}^{\prime}\left(
\omega,s\right)
^{T}S_{m}^{\prime}\left(  \omega,s\right)  ds~e_{a},e_{c}\right\rangle \\
&  +\sum_{j=m+1}^{n}\left\langle V_{mj}^{T}(\omega)\left[  \int_{0}^{\Delta_{j}s}%
C_{j}^{\prime}(\omega,s)^{T}C_{j}^{\prime}(\omega,s)ds\right]  V_{mj}(\omega)e_{a}%
,e_{c}\right\rangle
\end{align*}
and hence we have shown,%
\begin{align}
&\mathcal{Q}_{mm}^{n}(\omega)  \nonumber \\
&=\int_{0}^{\Delta_{m}s}S_{m}^{\prime}\left(
\omega,s\right)  ^{T}S_{m}^{\prime}\left(  \omega,s\right)  ds+\sum_{j=m+1}^{n}%
V_{mj}^{T}(\omega)\left[
\int_{0}^{\Delta_{j}s}C_{j}^{\prime}(\omega,s)^{T}C_{j}^{\prime
}(\omega,s)ds\right]  V_{mj}(\omega)\nonumber\\
&  =\int_{0}^{1/n}S_{m}^{\prime}\left(  \omega,s\right)
^{T}S_{m}^{\prime}\left( \omega,s\right)
ds+\sum_{j=m+1}^{n}V_{mj}^{T}(\omega)\left[
\int_{0}^{1/n}C_{j}^{\prime
}(\omega,s)^{T}C_{j}^{\prime}(\omega,s)ds\right]  V_{mj}(\omega). \label{eq.qmm}%
\end{align}
Noting that
\[
\left\Vert V_{mj}(\omega)\right\Vert =\left\Vert
V_{mj}^{T}(\omega)\right\Vert \leq\left\Vert
S_{m}(\omega,\Delta_{m}s)\right\Vert
^{2}\prod_{k=m+1}^{j-1}\left\Vert
C_{k}(\omega,\Delta_{k}s)\right\Vert ^{2},
\]
it follows from Equation (\ref{eq.qmm}) that%
\begin{equation}
\left\Vert n\mathcal{Q}_{mm}^{n}-I\right\Vert \leq A_{m}+B_{m} \label{eq.nqn}%
\end{equation}
where
\begin{equation}
A_{m}(\omega):=n\int_{0}^{1/n}\left\Vert S_{m}^{\prime}\left(  \omega,s\right)  ^{T}%
S_{m}^{\prime}\left(  \omega,s\right)  -I\right\Vert ds \label{eq.A}%
\end{equation}
and
\begin{align}
B_{m}(\omega)  &  =\sum_{j=m+1}^{n}\left\Vert
S_{m}(\omega,\Delta_{m}s)\right\Vert
^{2}\prod_{k=m+1}^{j-1}\left\Vert
C_{k}(\omega,\Delta_{k}s)\right\Vert ^{2}\left\Vert
n\int_{0}^{1/n}C_{j}^{\prime}(\omega,s)^{T}C_{j}^{\prime
}(\omega,s)ds\right\Vert \nonumber
\\[0.01in]
&  \leq\sum_{j=m+1}^{n}\left\Vert S_{m}(\omega,\Delta_{m}s)\right\Vert ^{2}%
\prod_{k=m+1}^{j-1}\left\Vert C_{k}(\omega,\Delta_{k}s)\right\Vert
^{2}n\int
_{0}^{1/n}\left\Vert C_{j}^{\prime}(\omega,s)\right\Vert ^{2}ds. \label{eq.B2}%
\end{align}
Thus we are now left to estimate the quantities comprising $A_{m}$
and $B_{m}.$

\subsection{Estimates for Solutions to Jacobi's Equation\label{s.jac-est}}

\begin{rem}
In what follows we will make use of the following elementary
estimates without further comment.

\begin{enumerate}
\item $\cosh
x=\sum_{n=0}^{\infty}\frac{x^{2n}}{(2n)!}\leq\sum_{n=0}^{\infty
}\frac{x^{2n}}{2^{n}n!}=e^{x^{2}/2}.$

\item $\cosh x=\cosh\left\vert x\right\vert =\frac{e^{\left\vert
x\right\vert }+e^{-\left\vert x\right\vert }}{2}\leq e^{\left\vert
x\right\vert }$ so that
\[
\cosh x\leq\mathrm{min}\left(  e^{x^{2}/2},e^{\left\vert
x\right\vert }\right)  ,
\]
and

\item for $x\geq0,$%
\[
\sinh x=\int_{0}^{x}\cosh t~dt\leq\int_{0}^{x}\cosh x~dt=x\cosh x.
\]
This estimate is also easily understood using the power series
expansions for $\sinh$ and $\cosh.$
\end{enumerate}
\end{rem}


\begin{lem}
[Global Estimate]\label{lem.jacobibound} Let $A(s)$ be a $d\times
d$ matrix for all $s\geq0,$ $\kappa:=\sup_{s\geq0}\left\Vert
A(s)\right\Vert <\infty,$ and let $Z(s)$ be either a
$\mathbb{R}^{d}$ or $d\times d$ matrix valued solution to the
second order differential equation
\begin{equation}
Z^{\prime\prime}(s)=A(s)Z(s).\nonumber
\end{equation}
Then
\begin{equation}
\left\Vert Z(s)-Z(0)\right\Vert \leq\left\Vert Z(0)\right\Vert
\left( \cosh\sqrt{\kappa}s-1\right)  +\left\Vert
Z^{\prime}(0)\right\Vert \frac
{\sinh\sqrt{\kappa}s}{\sqrt{\kappa}}. \label{eq.odebound1}%
\end{equation}

\end{lem}

\begin{proof}
By Taylor's theorem with integral remainder,
\begin{align}
Z(s)  &  =Z(0)+sZ^{\prime}(0)+\int_{0}^{s}Z^{\prime\prime}%
(u)(s-u)\ du\nonumber\\
&  =Z(0)+sZ^{\prime}(0)+\int_{0}^{s}A(u)Z(u)(s-u)\ du \label{eq.z-int}%
\end{align}
and therefore
\begin{align}
\left\Vert Z(s)-Z(0)\right\Vert  &  \leq s\left\Vert
Z^{\prime}(0)\right\Vert +\kappa\int_{0}^{s}\left\Vert Z\left(
u\right)  \right\Vert
(s-u)\ du\nonumber\\
&  \leq s\left\Vert Z^{\prime}(0)\right\Vert
+\kappa\int_{0}^{s}\left\Vert Z(u)-Z(0)\right\Vert (s-u)\
du+\frac{1}{2}s^{2}\kappa\left\Vert
Z(0)\right\Vert \nonumber \\
&:=f(s). \label{eq.odebound2}%
\end{align}
Note that $f(0)=0,$%
\begin{equation}
f^{\prime}(s)=\left\Vert Z^{\prime}(0)\right\Vert +\kappa\int_{0}%
^{s}\left\Vert Z(u)-Z(0)\right\Vert (s-u)\ du+s\kappa\left\Vert
Z(0)\right\Vert ,\nonumber
\end{equation}
$f^{\prime}(0)=\left\Vert Z^{\prime}(0)\right\Vert $, and
\begin{equation}
f^{\prime\prime}(s)=\kappa\parallel
Z(s)-Z(0)\parallel+\kappa\left\Vert Z(0)\right\Vert \leq\kappa
f(s)+\kappa\left\Vert Z(0)\right\Vert .\nonumber
\end{equation}
That is,
\begin{equation}
f^{\prime\prime}(s)=\kappa f(s)+\eta(s),\hspace{1cm}f(0)=0,\hspace
{1cm}\mathrm{and}\hspace{1cm}f^{\prime}(0)=\left\Vert
Z^{\prime}(0)\right\Vert
, \label{eq.odebound3}%
\end{equation}
where $\eta(s):=f^{\prime\prime}(s)-\kappa
f(s)\leq\kappa\left\Vert Z(0)\right\Vert $. Equation
(\ref{eq.odebound3}) may be solved by variation of parameters to
find
\begin{align}
f(s)  &  =\left\Vert Z^{\prime}(0)\right\Vert \frac{\sinh\sqrt{\kappa}s}%
{\sqrt{\kappa}}+\int_{0}^{s}\frac{\sinh\sqrt{\kappa}(s-r)}{\sqrt{\kappa}}%
\eta(r)\ dr\nonumber\\
&  \leq\left\Vert Z^{\prime}(0)\right\Vert \frac{\sinh\sqrt{\kappa}s}%
{\sqrt{\kappa}}+\left\Vert Z(0)\right\Vert
\int_{0}^{s}\sqrt{\kappa}\sinh
\sqrt{(}s-r)\ dr\nonumber\\
&  =\left\Vert Z^{\prime}(0)\right\Vert \frac{\sinh\sqrt{\kappa}s}%
{\sqrt{\kappa}}+\left\Vert Z(0)\right\Vert
(\cosh\sqrt{\kappa}s-1).\nonumber
\end{align}
\qquad Combining this with Equation (\ref{eq.odebound2}) proves
Equation (\ref{eq.odebound1}).
\end{proof}

\begin{thm}
\label{t.cs-est}Suppose that $A\left(  s\right)  $ above
satisfies, $0\leq-A\left(  s\right)  \leq\kappa I$ for all $s$ or
equivalently that
$-\kappa I\leq A\left(  s\right)  \leq0$ for all $s.$ Let%
\[
\psi\left(  s\right)  :=\mathrm{min}\left(  1+\frac{\cosh\left(
s\right) }{16}s^{4},\cosh\left(  s\right)  \right)
\]
whose graph is shown below and $C(s)$ and $S(s)$ be the matrix functions defined by%
\begin{align}
C''(s) &= A(s)C(s),\ \text{with}\ C(0)=I,\ C'(0) = 0, \nonumber \\
S''(s) &= A(s)S(s),\ \text{with}\ S(0)=0,\ S'(0) = I. \nonumber
\end{align}

\begin{figure}
[ptbh]
\begin{center}
\fbox{\includegraphics[ natheight=2.177000in, natwidth=3.265500in,
height=2.177in, width=3.2655in
]%
{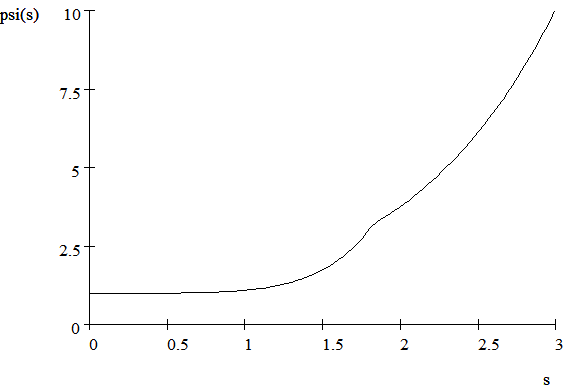}}
\caption{Graph of $\psi.$}%
\label{f.phi}%
\end{center}
\end{figure}

\newpage
Then

\begin{enumerate}
\item $\left\Vert C\left(  s\right)  \right\Vert \leq\psi\left(
\sqrt{\kappa }s\right)  ,$

\item $\left\Vert S\left(  s\right)  \right\Vert \leq s\psi\left(
\sqrt{\kappa}s\right)  ,$

\item $\left\Vert C^{\prime}\left(  s\right)  \right\Vert
\leq\kappa s\psi\left(  \sqrt{\kappa}s\right)  ,$

\item $\left\Vert S^{\prime}\left(  s\right)  \right\Vert \leq 1 + \frac{1}%
{2}\kappa s^{2}\psi\left(  \sqrt{\kappa}s\right)  ,$ and

\item $\left\Vert S^{\prime}\left(  s\right)  ^{T}S^{\prime}\left(
s\right)
-I\right\Vert \leq\psi\left(  \sqrt{\kappa}s\right)  \kappa s^{2}+\frac{1}%
{3}\psi^{2}\left(  \sqrt{\kappa}s\right)  \kappa^{2}s^4.$
\end{enumerate}

Moreover, if we only assume that $\left\Vert A\left(  s\right)
\right\Vert \leq\kappa,$ all of the above estimates still hold
provided that $\psi$ is replaced by $\cosh.$
\end{thm}

\begin{proof}
If $Z$ solves, $Z^{\prime\prime}\left(  s\right)  =A\left(
s\right)  Z\left(
s\right)  ,$ then iterating Equation (\ref{eq.z-int}) shows%
\begin{align*}
Z(s)  &  =Z(0)+sZ^{\prime}(0)+\int_{0}^{s}A(u)\left[  Z(0)+uZ^{\prime}%
(0)+\int_{0}^{u}A(r)Z(r)(u-r)\ dr\right]  (s-u)\ du\\
&  =\left(  I+\int_{0}^{s}(s-u)A(u)du\right)  Z(0)+\left(  sI+\int_{0}%
^{s}u(s-u)A(u)\ du\right)  Z^{\prime}(0)\\
&  +\int_{0\leq r\leq u\leq s}(u-r)(s-u)A(u)A(r)Z(r)drdu.
\end{align*}
In particular this shows%
\begin{align}
C\left(  s\right) & =\left(  I-\int_{0}^{s}(s-u)\left[
-A(u)\right]
du\right)  +\int_{0\leq r\leq u\leq s}(u-r)(s-u)A(u)A(r)C(r)drdu \label{eq.ic}%
\end{align}
and%
\begin{align}
S\left(  s\right) & =\left(  sI-\int_{0}^{s}u(s-u)\left[
-A(u)\right] \ du\right)  +\int_{0\leq r\leq u\leq
s}(u-r)(s-u)A(u)A(r)S(r)drdu.
\label{eq.is}%
\end{align}
From Equation (\ref{eq.ic}),
\begin{align*}
\left\Vert C\left(  s\right)  \right\Vert  &  \leq\left\Vert I+\int_{0}%
^{s}(s-u)A(u)du\right\Vert +\int_{0\leq r\leq u\leq
s}(u-r)(s-u)\left\Vert
A(u)A(r)C(r)\right\Vert drdu\\
&  \leq\left\Vert I+\int_{0}^{s}(s-u)A(u)du\right\Vert
+\frac{\cosh\left( \sqrt{\kappa}s\right)  }{24}\kappa^{2}s^{4}.
\end{align*}
Moreover,
\[
\left(  1-\kappa\frac{s^{2}}{2}\right)  I=I-\kappa\left[  \int_{0}%
^{s}(s-u)du\right]  I\leq I+\int_{0}^{s}(s-u)A(u)du\leq I
\]
from which it follows that
\[
\left\Vert I+\int_{0}^{s}(s-u)A(u)du\right\Vert \leq\max\left(
1,\kappa
\frac{s^{2}}{2}-1\right)  =1_{s^{2}\leq4/\kappa}+\left(  \kappa\frac{s^{2}}%
{2}-1\right)  1_{s^{2}\geq4/\kappa}%
\]
and hence we have%
\begin{align*}
\left\Vert C\left(  s\right)  \right\Vert  &
\leq1_{s^{2}\leq4/\kappa }+\left(  \kappa\frac{s^{2}}{2}-1\right)
1_{s^{2}\geq4/\kappa}+\frac {\cosh\left(  \sqrt{\kappa}s\right)
}{24}\kappa^{2}s^{4}  \leq1+\frac{\cosh\left(
\sqrt{\kappa}s\right)  }{16}\kappa^{2}s^{4}.
\end{align*}
This is because $ f(s) = s^4\frac{\cosh s}{48} -
\left(\frac{s^2}{2}-1\right) = 1 + s^2 \left[ s^2\frac{\cosh
s}{48} - \frac{1}{2} \right] $ is an increasing function and for
$s \geq 2$, $f(s) \geq f(2) > 0.25$. Recall that we also know that
$\left\Vert C\left( s\right) \right\Vert
\leq\cosh\left(  \sqrt{\kappa}s\right)  $ and therefore we have%
\begin{align*}
\left\Vert C\left(  s\right)  \right\Vert  &
\leq\mathrm{min}\left(
1_{s^{2}\leq4/\kappa}+\left(  \kappa\frac{s^{2}}{2}-1\right)  1_{s^{2}%
\geq4/\kappa}+\frac{\cosh\left(  \sqrt{\kappa}s\right)  }{24}\kappa^{2}%
s^{4},\cosh\left(  \sqrt{\kappa}s\right)  \right) \\
&  \leq\mathrm{min}\left(  1+\frac{\cosh\left(  \sqrt{\kappa}s\right)  }%
{16}\kappa^{2}s^{4},\cosh\left(  \sqrt{\kappa}s\right)  \right)  .
\end{align*}
It will be convenient to define%
\[
\psi\left(  s\right)  :=\mathrm{min}\left(  1 +\frac{\cosh\left(  s\right)  }{16}%
s^{4},\cosh\left(  s\right)  \right)
\]
and hence
\[
\left\Vert C\left(  s\right)  \right\Vert \leq \psi\left(
\sqrt{\kappa}s\right).
\]
Similarly, from Equation (\ref{eq.is})%
\begin{align*}
\left\Vert S\left(  s\right)  \right\Vert  &  \leq\left\Vert sI+\int_{0}%
^{s}u(s-u)A(u)\ du\right\Vert +\int_{0\leq r\leq u\leq
s}(u-r)(s-u)\left\Vert
A(u)A(r)S(r)\right\Vert drdu\\
&  \leq s\left\Vert I+\frac{1}{s}\int_{0}^{s}u(s-u)A(u)\
du\right\Vert
+\frac{\kappa^{2}}{24}s^{4}\frac{\sinh\left(  \sqrt{\kappa}s\right)  }%
{\sqrt{\kappa}}.
\end{align*}
In this case,%
\[
\left(  1-\frac{\kappa}{6}s^{2}\right)  I=\left(
1-\kappa\frac{1}{s}\int _{0}^{s}u(s-u)\ du\right)  I\leq
I+\frac{1}{s}\int_{0}^{s}u(s-u)A(u)\ du\leq I.
\]
Combining this with the previous equation shows%
\begin{align*}
\left\Vert S\left(  s\right)  \right\Vert  &  \leq s\left[  1_{s^{2}%
\leq12/\kappa}+\left(  \kappa\frac{s^{2}}{6}-1\right)
1_{s^{2}\geq12/\kappa }\right]
+\frac{\kappa^{2}}{24}s^{4}\frac{\sinh\left( \sqrt{\kappa}s\right)
}{\sqrt{\kappa}}\\
&  \leq s\left[  1_{s^{2}\leq12/\kappa}+\left(  \kappa\frac{s^{2}}%
{6}-1\right)  1_{s^{2}\geq12/\kappa}+\frac{\cosh\left(
\sqrt{\kappa}s\right)
}{24}\kappa^{2}s^{4}\right] \\
&  \leq s\left[  1+\frac{\cosh\left(  \sqrt{\kappa}s\right)
}{16}\kappa ^{2}s^{4}\right]  .
\end{align*}
Since we also have
\[
\left\Vert S\left(  s\right)  \right\Vert \leq\frac{\sinh\left(
\sqrt{\kappa }s\right)  }{\sqrt{\kappa}}\leq s\cosh\left(
\sqrt{\kappa}s\right),
\]
we may conclude that
\[
\left\Vert S\left(  s\right)  \right\Vert \leq
s~\mathrm{min}\left( 1+\frac{\cosh\left(  \sqrt{\kappa}s\right)
}{16}\kappa^{2}s^{4},\cosh\left( \sqrt{\kappa}s\right)  \right)
=s\psi\left(  \sqrt{\kappa}s\right)  .
\]

Furthermore,%
\[
\left\Vert C^{\prime}\left(  s\right)  \right\Vert =\left\Vert \int_{0}%
^{s}A\left(  r\right)  C\left(  r\right)  dr\right\Vert \leq\int_{0}%
^{s}\left\Vert A\left(  r\right)  \right\Vert \left\Vert C\left(
r\right) \right\Vert dr\leq\kappa s\psi\left(
\sqrt{\kappa}s\right)
\]
and%
\begin{align*}
\left\Vert S^{\prime}\left(  s\right)  \right\Vert  & = \left\Vert
I +\int
_{0}^{s}A\left(  r\right)  S\left(  r\right)  dr\right\Vert \leq 1 + \int_{0}%
^{s}\left\Vert A\left(  r\right)  \right\Vert \left\Vert S\left(
r\right)
\right\Vert dr\\
&  \leq 1 + \kappa\int_{0}^{s}r\psi\left(  \sqrt{\kappa}r\right)  dr\leq 1 + \frac{1}%
{2}\kappa s^{2}\psi\left(  \sqrt{\kappa}s\right)  .
\end{align*}
Finally%
\[
\frac{d}{ds}\left[  S^{\prime}\left(  s\right)
^{T}S^{\prime}\left( s\right)  \right]  =S\left(  s\right)
^{T}A\left(  s\right)  S^{\prime }\left(  s\right)
+S^{\prime}\left(  s\right)  ^{T}A\left(  s\right) S\left(
s\right)
\]
and therefore,%
\begin{align*}
\left\Vert \frac{d}{ds}\left[  S^{\prime}\left(  s\right)
^{T}S^{\prime }\left(  s\right)  \right]  \right\Vert  &
\leq2\left\Vert S\left(  s\right) \right\Vert \left\Vert A\left(
s\right)  \right\Vert \left\Vert S^{\prime
}\left(  s\right)  \right\Vert \\
&  \leq2\left\Vert S\left(  s\right)  \right\Vert \left\Vert
A\left(
s\right)  \right\Vert \left\Vert S^{\prime}\left(  s\right)  \right\Vert \\
&  \leq2\kappa s\psi\left(  \sqrt{\kappa}s\right)  \left(
1+\frac{1}{2}\kappa s^{2}\psi\left(  \sqrt{\kappa}s\right) \right)
.
\end{align*}
Integrating this equation then implies,%
\begin{align*}
\left\Vert S^{\prime}\left(  s\right)  ^{T}S^{\prime}\left(
s\right) -I\right\Vert  &  \leq\int_{0}^{s}2\kappa r\psi\left(
\sqrt{\kappa}r\right) \left(  1+\frac{1}{2}\kappa r^{2}\psi\left(
\sqrt{\kappa}r\right)  \right)
dr\\
&  \leq\int_{0}^{s}2\kappa\psi\left(  \sqrt{\kappa}s\right) \left(
r+\frac{1}{2}\kappa r^{3}\psi\left(  \sqrt{\kappa}s\right)  \right)  dr\\
&  =\kappa\psi\left(  \sqrt{\kappa}s\right)  \left(
s^{2}+\frac{1}{3}\kappa
s^{4}\psi\left(  \sqrt{\kappa}s\right)  \right) \\
&  =\psi\left(  \sqrt{\kappa}s\right)  \kappa
s^{2}+\frac{1}{3}\psi^{2}\left( \sqrt{\kappa}s\right)
\kappa^{2}s^{4}.
\end{align*}

\end{proof}

\begin{prop}
\label{p.expand}If $\left\Vert A\left(  s\right)  \right\Vert
\leq\kappa,$
then the following estimates hold:%
\begin{equation}
\left\Vert S^{\prime}\left(  s\right)  -\left(
I+\int_{0}^{s}rA\left( r\right)  dr\right)  \right\Vert \leq
s^{4}\kappa^{2}\cosh\left(  \sqrt
{\kappa}s\right)  \label{eq.sp}%
\end{equation}%
\begin{equation}
\left\Vert \frac{S\left(  s\right)
}{s}-\left(I+\frac{1}{s}\int_{0}^{s}\left( s-r\right)  rA\left(
r\right) dr\right)\right\Vert \leq s^{4}\kappa^{2}\cosh\left(
\sqrt{\kappa}s\right)  \label{eq.s_s}%
\end{equation}%
\begin{equation}
\left\Vert C\left(  s\right)  -\left(
I+\int_{0}^{s}(s-u)A(u)du\right) \right\Vert \leq
s^{4}\kappa^{2}\cosh\left(  \sqrt{\kappa}s\right)  .
\label{eq.cs}%
\end{equation}

\end{prop}

\begin{proof}
Now%
\begin{align}
S^{\prime}\left(  s\right)   &  =I+\int_{0}^{s}A\left(  r\right)
S\left( r\right)  dr=I+\int_{0}^{s}A\left(  r\right)  \left[
\int_{0}^{r}S^{\prime
}\left(  u\right)  du\right]  dr\nonumber\\
&  =I+\int_{0}^{s}A\left(  r\right)  \left[  \int_{0}^{r}\left[
I+\int _{0}^{u}A\left(  v\right)  S\left(  v\right)  dv\right]
du\right]
dr\nonumber\\
&  =I+\int_{0}^{s}rA\left(  r\right)  dr+\int_{0\leq v\leq u\leq
r\leq
s}A\left(  r\right)  A\left(  v\right)  S\left(  v\right)  dvdudr\nonumber\\
&  =I+\int_{0}^{s}rA\left(  r\right)  dr+\int_{0\leq v\leq r\leq
s}\left( r-v\right)  A\left(  r\right)  A\left(  v\right)  S\left(
v\right)  dvdr. \nonumber 
\end{align}
Thus 
\begin{align*}
\left\Vert S^{\prime}\left(  s\right)  -\left(
I+\int_{0}^{s}rA\left( r\right)  dr\right)  \right\Vert  &
\leq\int_{0\leq v\leq r\leq s}\left(
r-v\right)  \kappa^{2}v\cosh\left(  \sqrt{\kappa}v\right)  dvdr\\
&  \leq s^{4}\kappa^{2}\cosh\left(  \sqrt{\kappa}s\right)  ,
\end{align*}
Integrating this estimate implies%
\[
\left\Vert S\left(  s\right)  -Is-\int_{0}^{s}\left(  s-r\right)
A\left( r\right)  dr\right\Vert \leq s\cdot
s^{4}\kappa^{2}\cosh\left(  \sqrt{\kappa }s\right)
\]
which is equivalent to Equation (\ref{eq.s_s}). Similarly from
Equation (\ref{eq.ic}) we have
\begin{align*}
\left\Vert C\left(  s\right)  -\left(
I+\int_{0}^{s}(s-u)A(u)du\right)
\right\Vert  &  =\left\Vert \int_{0\leq r\leq u\leq s}%
(u-r)(s-u)A(u)A(r)C(r)drdu\right\Vert \\
&  \leq s^{4}\kappa^{2}\cosh\left(  \sqrt{\kappa}s\right)  .
\end{align*}

\end{proof}

\begin{figure}
[ptbh]
\begin{center}
\fbox{\includegraphics[ natheight=2.159400in, natwidth=3.239600in,
height=1.9594in, width=3.0396in
]%
{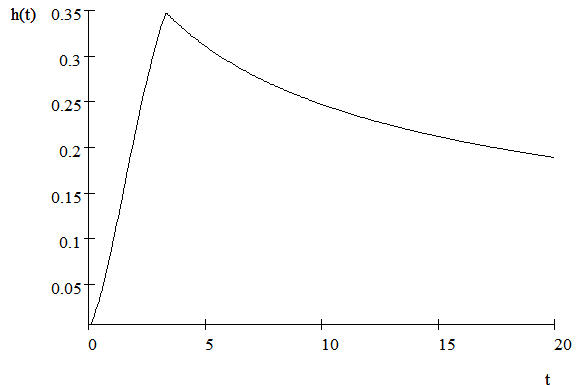}}
\caption{Here is the graph of $h\left(  t\right)  .$}%
\label{fig.h-graph}%
\end{center}
\end{figure}

\begin{defn}
\label{def.h}Let%
\[
h\left(  t\right)  =\frac{1}{t}\ln\psi\left(  \sqrt{t}\right)  =t^{-1}%
\ln\left(  \mathrm{min}{\left(  1+\frac{\cosh\left(  \sqrt{t}\right)  }%
{16}t^{2},\cosh\left(  \sqrt{t}\right)  \right)  }\right)
\]
whose graph is given in Figure \ref{fig.h-graph}, i.e.
\[
e^{s^{2}h\left(  s^{2}\right)  }=\psi\left(  s\right)  .
\]

\end{defn}%

Let $\varphi$ be a function, which we will specify shortly, such
that $\varphi\geq\psi.$ Further let $g\left(  t\right)
:=\frac{1}{t}\ln \varphi\left(  \sqrt{t}\right)  $ so that
$\varphi\left(  s\right)
=e^{s^{2}g\left(  s^{2}\right)  }$ and define%
\[
u\left(  t\right)  :=\frac{\psi^{2}\left(  \sqrt{t}\right)
}{\varphi ^{2}\left(  \sqrt{t}\right)  }=\frac{e^{2th\left(
t\right)  }}{e^{2tg\left( t\right)  }}=e^{-2t\left(  g\left(
t\right)  -h\left(  t\right)  \right)  }.
\]

We will specify $\varphi\left(  t\right)  $ by requiring $g\left(
t\right) $ to be a smooth function such that $g\left(  t\right)
=h\left( t\right)  $ for $t$ near zero which then rises rapidly to
a height of $.6$ as $t$ increases. For later purposes, let us
observe that with this definition
$tu\left(  t\right)  $ is bounded by $.63$ as the graph below indicates.%
\begin{figure}
[ptbhptbh]
\begin{center}
\fbox{\includegraphics[ natheight=2.316000in, natwidth=3.474000in,
height=2.016in, width=3.174in
]%
{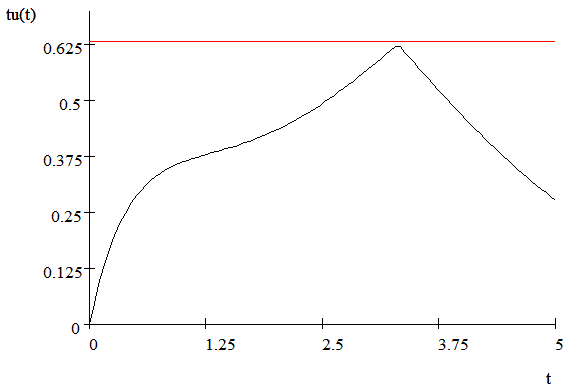}}
\caption{Graph of $tu\left(  t\right)  $ and the line $y=.63.$ We
see that
$tu\left(  t\right)  \leq.63.$}%
\label{f.tut}%
\end{center}
\end{figure}

\subsection{Proof of Uniform Integrability\label{s.unif-i-d}}

\begin{prop}
\label{prop.ints1}Suppose $N$ is a $d$ -- dimensional standard
Gaussian normal random variable and $G\left(  x,\left\Vert
N\right\Vert \right)  $ is a $C^{1}$ -- function in $x\in\left(
-\varepsilon,\varepsilon\right)  $ such that $G\left( 0,\left\Vert
N\right\Vert \right)  \equiv0,$
\[
\alpha:=\sup\left\{  G\left(  x,t\right)  :x\in\left( -\varepsilon
,\varepsilon\right)  \text{ and }t\geq0\right\}  <\frac{1}{2}%
\]
and there exist constants $C<\infty$ and $\beta<\infty$ such that
\[
G_{x}\left(  x,t\right) = \frac{\partial G}{\partial x}(x,t)  \leq C\left(  1+t\right)  ^{\beta}\text{ for }%
x\in\left(  -\varepsilon,\varepsilon\right)  \text{ and }t\geq0.
\]
Then
\[
\lim_{n\rightarrow\infty}\left(  \mathbb{E}\left[  e^{\left\Vert
N\right\Vert ^{2}G\left(  \frac{1}{n},\left\Vert N\right\Vert
\right)  }\right]  \right) ^{n}=e^{\mathbb{E}\left[  \left\Vert
N\right\Vert ^{2}G_{x}\left( 0,\left\Vert N\right\Vert \right)
\right]  }.
\]

\end{prop}

\begin{proof}
Let%
\[
f\left(  x\right)  :=\mathbb{E}\left[  e^{\left\Vert N\right\Vert
^{2}G\left( x,\left\Vert N\right\Vert \right)  }\right]  .
\]
Since
\[
e^{\left\Vert N\right\Vert ^{2}G\left(  x,\left\Vert N\right\Vert
\right)
}\leq e^{\alpha\left\Vert N\right\Vert ^{2}}%
\]
and%
\[
G_{x}\left(  x,\left\Vert N\right\Vert \right)  \left\Vert
N\right\Vert ^{2}e^{\left\Vert N\right\Vert ^{2}G\left(
x,\left\Vert N\right\Vert \right) }\leq C\left(  1+\left\Vert
N\right\Vert \right)  ^{\beta}\left\Vert
N\right\Vert ^{2}e^{\alpha\left\Vert N\right\Vert ^{2}}%
\]
with the right hand members of these inequalities being integrable
functions, it follows that $f\left(  x\right)  $ is a $C^{1}$ -
function for $x$ near $0$
with%
\[
f^{\prime}\left(  x\right)  =\mathbb{E}\left[  G_{x}\left(
x,\left\Vert N\right\Vert \right)  \left\Vert N\right\Vert
^{2}e^{\left\Vert N\right\Vert ^{2}G\left(  x,\left\Vert
N\right\Vert \right)  }\right]  .
\]
In particular we have $f\left(  0\right)  =\mathbb{E}\left[
e^{\left\Vert N\right\Vert ^{2}G\left(  0,\left\Vert N\right\Vert
\right)  }\right]  =1$ and
\[
f^{\prime}\left(  0\right)  =\mathbb{E}\left[  G_{x}\left(
0,\left\Vert N\right\Vert \right)  \left\Vert N\right\Vert
^{2}\right]
\]
Therefore,%
\begin{align*}
\lim_{n\rightarrow\infty}\left(  \mathbb{E}\left[  e^{\left\Vert
N\right\Vert ^{2}G\left(  \frac{1}{n},\left\Vert N\right\Vert
\right)  }\right]  \right) ^{n}  &
=\lim_{n\rightarrow\infty}\left(  f\left(  n^{-1}\right)  \right)
^{n}=\lim_{n\rightarrow\infty}e^{n\ln f\left(  n^{-1}\right)
}=\lim
_{x\rightarrow0}e^{\frac{1}{x}\ln f\left(  x\right)  }\\
&  =\lim_{x\rightarrow0}e^{f^{\prime}\left(  x\right)  /f\left(
x\right) }=e^{f^{\prime}\left(  0\right)  }=e^{\mathbb{E}\left[
G_{x}\left( 0,\left\Vert N\right\Vert \right)  \left\Vert
N\right\Vert ^{2}\right]  }.
\end{align*}

\end{proof}

\begin{thm}
\label{thm.uniformint} \label{thm.intergab}Suppose that $(M,g)$ is
a Riemannian manifold with non-negative sectional curvatures which
are bounded
above by $K=\frac{1}{2d}.$ Then for all $p$ sufficiently close to $1,$%
\[
\sup_{n}\mathbb{E}\left[  \det{}^{p/2}(n\mathcal{Q}^{n}\circ
b_n)\right] <\infty.
\]

\end{thm}

\begin{proof}
For $\omega \in H_{\cP}(\mathbb{R}^d)$, from Theorem \ref{t.cs-est},%
\begin{align*}
A_{m}(\omega)  &  =n\int_{0}^{1/n}\left\Vert S_{m}^{\prime}\left(
\omega,s\right)
^{T}S_{m}^{\prime}\left(  \omega,s\right)  -I\right\Vert ds\\
&  \leq\psi\left(  \sqrt{K}\left\Vert \Delta_{m}\omega\right\Vert
\right) K\left\Vert \Delta_{m}\omega\right\Vert
^{2}+\frac{1}{3}\psi^{2}\left(  \sqrt {K}\left\Vert
\Delta_{m}\omega\right\Vert \right)  K^{2}\left\Vert \Delta
_{m}\omega\right\Vert ^{4},
\end{align*}
and if we write $\tau_j(\omega) = \sqrt{K}\left\Vert
\Delta_{j}\omega\right\Vert$,
\begin{align*}
B_{m}(\omega)  &  =\sum_{j=m+1}^{n}\left\Vert
S_{m}(\omega,\Delta_{m}s)\right\Vert
^{2}\prod_{k=m+1}^{j-1}\left\Vert C_{k}(\omega,\Delta_{k}s)\right\Vert ^{2}%
n\int_{0}^{1/n}\left\Vert C_{j}^{\prime}(\omega,s)\right\Vert ^{2}ds\\
&  \leq\sum_{j=m+1}^{n}\left(  \Delta_{m}s\right)
^{2}\psi^{2}(\tau_m(\omega) )\left[  \prod_{k=m+1}^{j-1}\psi
^{2}\left( \tau_k(\omega) \right)  \right] K^{2}\left(
\Delta_{j}s\right) ^{2}\left\Vert \omega_{j}^{\prime}\right\Vert
^{4}\psi^{2}(\tau_j(\omega) )\\
&  =K^{2}\sum_{j=m+1}^{n}\left[  \prod_{k=m}^{j}\psi^{2}\left(
\tau_k(\omega) \right)  \right] \left\Vert \Delta
_{j}\omega\right\Vert ^{4} \nonumber \\
&=K^{2}\sum_{j=m+1}^{n}\left[ \prod_{k=m}^{j}\psi^{2}\left( \sqrt
{K}\left\Vert \Delta_{k}\omega\right\Vert \right)  \right]
\left\Vert \Delta _{j}\omega\right\Vert ^{4}.
\end{align*}
Hence if we let%
\[
\alpha(\omega):=\prod_{k=1}^{n}\varphi^{2}\left(
\sqrt{K}\left\Vert \Delta
_{k}\omega\right\Vert \right)  =\prod_{k=1}^{n}e^{2K\left\Vert \Delta_{k}%
\omega\right\Vert ^{2}g\left(  K\left\Vert
\Delta_{k}\omega\right\Vert ^{2}\right)  },
\]
then%
\begin{align*}
A_{m}(\omega)  &  \leq\alpha (\omega)K\left\Vert
\Delta_{m}\omega\right\Vert ^{2}+\frac{1}{3}\alpha(\omega) u\left(
K\left\Vert \Delta_{m}\omega\right\Vert ^{2}\right)
K^{2}\left\Vert
\Delta_{m}\omega\right\Vert ^{4}\\
&  \leq\alpha(\omega)\left[  K\left\Vert \Delta_{m}\omega\right\Vert ^{2}+\frac{1}%
{3}u\left(  K\left\Vert \Delta_{m}\omega\right\Vert ^{2}\right)
K^{2}\left\Vert \Delta_{m}\omega\right\Vert ^{4}\right]
\end{align*}
and
\begin{align*}
B_{m}(\omega)  &  \leq\alpha(\omega) K^{2}\sum_{j=m+1}^{n}\left[
\prod_{k=m}^{j}u\left( \sqrt{K}\left\Vert
\Delta_{k}\omega\right\Vert \right)  \right]  \left\Vert
\Delta_{j}\omega\right\Vert ^{4} \nonumber \\
&\leq\alpha(\omega) K^{2}\sum_{j=1}^{n}u\left( \sqrt{K}\left\Vert
\Delta _{j}\omega\right\Vert \right)  \left\Vert
\Delta_{j}\omega\right\Vert ^{4}.
\end{align*}
For $\omega \in (W(\mathbb{R}^d), \mu)$, \beq \Delta_ib_n(\omega)
=\Delta_ib(\omega) = b(s_i)(\omega) - b(s_{i-1})(\omega) =
\omega(s_i) - \omega(s_{i-1})
= \Delta_i\omega. \nonumber \eeq Therefore, on $W(\mathbb{R}^d)$,%
\begin{align*}
& \sum_{m=1}^{n}\left(  A_{m}+B_{m}\right) \circ b_n  \nonumber \\
& \leq\alpha\circ b_n \cdot\left[
\begin{array}
[c]{c}%
\sum_{m=1}^{n}\left[  K\left\Vert \Delta_{m}b\right\Vert ^{2}+\frac{1}%
{3}u\left(  K\left\Vert \Delta_{m}b\right\Vert ^{2}\right)
K^{2}\left\Vert
\Delta_{m}b\right\Vert ^{4}\right] \\
+K^{2}n\sum_{m=1}^{n}u\left(  K\left\Vert \Delta_{m}b\right\Vert
^{2}\right)
\left\Vert \Delta_{m}b\right\Vert ^{4}%
\end{array}
\right] \\
&  =\alpha\circ b_n \cdot\sum_{m=1}^{n}\left[  K\left\Vert
\Delta_{m}b\right\Vert ^{2}+\left(  \frac{1}{3}+n\right)
K^{2}u\left( \sqrt{K}\left\Vert \Delta _{m}b\right\Vert \right)
\left\Vert \Delta_{m}b\right\Vert ^{4}\right]  .
\end{align*}

Now let $x=n^{-1}$ and $N_{m}:=\sqrt{n}\Delta_{m}b,$ so that
$\left\{ N_{m}\right\}  _{m=1}^{n}$ is a collection of
$\mathbb{R}^{d}$ -- valued
independent standard normal random variables. With this notation we have%
\begin{align*}
&\sum_{m=1}^{n}\left(  A_{m}+B_{m}\right)\circ b_n  \nonumber \\
 &
\leq\alpha\circ b_n \cdot \sum_{m=1}^{n}\left[ Kx\left\Vert
N_{m}\right\Vert ^{2}+\left(  \frac{1}{3}+n\right)  K^{2}u\left(
Kx\left\Vert N_{m}\right\Vert ^{2}\right)  x^{2}\left\Vert
N_{m}\right\Vert
^{4}\right] \\
&  =\alpha\circ b_n \cdot \sum_{m=1}^{n}\left[  Kx\left\Vert
N_{m}\right\Vert ^{2}+K\left( \frac{1}{3}x+1\right)  u\left(
Kx\left\Vert N_{m}\right\Vert ^{2}\right)
Kx\left\Vert N_{m}\right\Vert ^{4}\right] \\
&  =\alpha\circ b_n \cdot \sum_{m=1}^{n}\left[  Kx+K\left(
\frac{1}{3}x+1\right) u\left( Kx\left\Vert N_{m}\right\Vert
^{2}\right)  Kx\left\Vert N_{m}\right\Vert
^{2}\right]  \left\Vert N_{m}\right\Vert ^{2}%
\end{align*}
and therefore, using Equation (\ref{eq.est}),%
\begin{align*}
&\det(n\mathcal{Q}^{n}\circ b_n)  = \rho_n \circ \phi \circ b_n
\nonumber \\
& \leq\left(\alpha^{nd}\exp\left( \alpha^{-1}d\sum
_{m=1}^{n}\left\Vert \left(n\mathcal{Q}_{mm}^{n}-I\right)\right\Vert \right) \right)\circ b_n \\
& \leq\left(\alpha^{nd}\exp\left(  \alpha^{-1}d\sum_{m=1}^{n}\left(  A_{m}%
+B_{m}\right)  \right) \right)\circ b_n \\
&  \leq\left(\alpha\circ b_n\right)^{nd}\exp\left(  d\sum_{m=1}^{n}\left[  Kx+K\left(  \frac{1}%
{3}x+1\right)  u\left(  Kx\left\Vert N_{m}\right\Vert ^{2}\right)
Kx\left\Vert N_{m}\right\Vert ^{2}\right]  \left\Vert
N_{m}\right\Vert ^{2}\right)
\end{align*}
where
\begin{align*}
(\alpha \circ b_n)^{nd}  &  =\prod_{m=1}^{n}\varphi^{2nd}\left(
\sqrt{K}\left\Vert \Delta_{m}b\right\Vert \right)
=\prod_{m=1}^{n}\varphi^{2nd}\left(  \sqrt
{Kx}\left\Vert N_{m}\right\Vert \right) \\
&  =\prod_{m=1}^{n}e^{2ndK\left\Vert \Delta_{m}b\right\Vert
^{2}g\left(
K\left\Vert \Delta_{m}b\right\Vert ^{2}\right)  }=\prod_{m=1}^{n}%
e^{2dK\left\Vert N_{m}\right\Vert ^{2}g\left(  Kx\left\Vert
N_{m}\right\Vert ^{2}\right)  }.
\end{align*}
Let $\theta\left( x, \left\Vert N \right \Vert^2 \right) =
x+\left( \frac{1}{3}x+1\right)  u\left( Kx\left\Vert N\right\Vert
^{2}\right)  Kx\left\Vert N\right\Vert ^{2} $.
Then,%
\begin{align*}
\det(n\mathcal{Q}^{n}\circ b_n)  &  \leq\prod_{m=1}^{n}
e^{2dK\left\Vert N_{m}\right\Vert ^{2}g\left(  Kx\left\Vert
N_{m}\right\Vert ^{2} \right)} \exp\left(  dK \theta\left(x, \left
\Vert N_m \right \Vert^2 \right)
\left\Vert N_{m}\right\Vert ^{2}  \right) \\
&  =\prod_{m=1}^{n}\exp\left(  dK\left[  2g\left(  Kx\left\Vert N_{m}%
\right\Vert ^{2}\right)  +\theta\left(x, \left \Vert N_m \right
\Vert^2 \right)\right]  \left\Vert N_{m}\right\Vert ^{2}\right)  .
\end{align*}
Hence
\begin{align*}
\mathbb{E}\left[  \det{}^{p/2}(n\mathcal{Q}^{n}\circ b_n)\right]
& =\left[ \mathbb{E}\exp\left(  dpK\left[  g\left(  Kx\left\Vert
N\right\Vert ^{2}\right)  +\frac{1}{2}\theta\left(x, \left \Vert
N_m \right \Vert^2 \right)\right]
\left\Vert N\right\Vert ^{2}\right)  \right]  ^{n}\\
&  =\left[  \mathbb{E}\exp\left(  dpK\cdot G\left(  x,\left\Vert
N\right\Vert
^{2}\right)  \left\Vert N\right\Vert ^{2}\right)  \right]  ^{1/x}%
\end{align*}
where
\begin{align*}
G\left(  x,\left\Vert N\right\Vert ^{2}\right) &=g\left(
Kx\left\Vert
N\right\Vert ^{2}\right)  +\frac{x}{2}+\left(  \frac{1}{6}x+\frac{1}%
{2}\right)  u\left(  Kx\left\Vert N\right\Vert ^{2}\right)
Kx\left\Vert
N\right\Vert ^{2}.
\end{align*}
By our choice of $g$ and hence $u$, we know
\[
g\left(  Kx\left\Vert N\right\Vert ^{2}\right) +\frac{1}{2}u\left(
Kx\left\Vert N\right\Vert ^{2}\right)  Kx\left\Vert N\right\Vert ^{2}%
\leq.6+\frac{1}{2}.63=.915<1.
\]
Therefore, for small $x,$ $G\left(  x,\left\Vert N\right\Vert
^{2}\right) \leq.92<1$ for small $x.$ Hence, if $p$ is
sufficiently close to $1,$ we will have
\[
dpK\cdot G\left(  x,\left\Vert N\right\Vert ^{2}\right)
=\frac{1}{2}pG\left( x,\left\Vert N\right\Vert ^{2}\right)
\leq\frac{1}{2}p\cdot0.92<\frac{1}{2}.
\]
Therefore we may apply Proposition \ref{prop.ints1} to conclude
that
\begin{align*}
&\limsup_{n\rightarrow\infty}\mathbb{E}\left[  \det{}^{p/2}(n\mathcal{Q}%
^{n}\circ b_n)\right]     \leq\lim_{x\rightarrow0}\left[
\mathbb{E}\exp\left( dpK\cdot G\left(  x,\left\Vert N\right\Vert
^{2}\right)  \left\Vert
N\right\Vert ^{2}\right)  \right]  ^{1/x}\\
&  =\exp\left(  \mathbb{E}\left[  dpK\cdot G_{x}\left(
0,\left\Vert
N\right\Vert ^{2}\right)  \left\Vert N\right\Vert ^{2}\right]  \right) \\
&  =\exp\left(  \mathbb{E}\left[  dpK\left(  \frac{2K}%
{16}\left\Vert N\right\Vert ^{2}+
\frac{1}{2}+\frac{1}{2}K\left\Vert N\right\Vert ^{2}\right)
\left\Vert N\right\Vert ^{2}\right] \right)  <\infty.
\end{align*}

\end{proof}

\section{Second Formula for $\rho_{ n}$} \label{s.2nd-rho}

\begin{defn}
For any $\epsilon> 0$ and any partition $\cP$ of $[0,1]$, let
\begin{align}
H_{ {\mathcal{P}}}^{\epsilon}(\mathbb{R}^{d})  &  = \Big\{ \omega
\in H_{ {\mathcal{P}}}(\mathbb{R}^{d})\ \Big|\
\int_{s_{i-1}}^{s_{i}}
\parallel \omega^{\prime
}(s) \parallel ds < \epsilon\ for\ i = 1, \dots, n \Big\}\nonumber\\
&  = \{ \omega \in H_{ {\mathcal{P}}}(\mathbb{R}^{d})\ |
\parallel\Delta_{i}\omega
\parallel< \epsilon\ \forall i \},\nonumber
\end{align}
where $\Delta_i\omega = \omega(s_i) - \omega(s_{i-1})$. The second
equality holds since $\omega_{i}^{\prime}(s)$ is a constant in
$J_{i}$ for each $i$.
\end{defn}

\begin{rem}
We will now consider all the $\omega \in H_{ {\mathcal{P}}}^{\epsilon}%
(\mathbb{R}^{d})$, with $\epsilon$ sufficiently small, specified
in the next lemma.
\end{rem}

\begin{lem}
\label{lem.Sinv} There exists an $\epsilon$ with
$\bigvee_{i=1,\ldots ,n}
\parallel\Delta_{i}\omega \parallel< \epsilon$ such that for $i = 1,2, \ldots n$,
$S_{i}(\omega, s)$ is invertible for $0< s \leq \Delta_is$.
\end{lem}

\begin{proof}
From Equation (\ref{eq.odebound1}) with $\kappa= \frac{K}{(\Delta_{i}s)^{2}%
}\parallel\Delta_{i}\omega \parallel^{2}$, we see that
\begin{equation}
\parallel S_{i}(\omega, s) \parallel\leq s \left( \frac{\sinh(\sqrt K \parallel
\Delta_{i}\omega \parallel)}{\sqrt K \parallel\Delta_{i}\omega
\parallel} \right) \leq
s\cosh(\sqrt K \parallel\Delta_{i}\omega \parallel), \label{eq.Sbd}%
\end{equation}
where we have used the inequality $\frac{\sinh x}{x} \leq\cosh x$.
By Taylor's Theorem with integral remainder,
\begin{equation}
S_{i}(\omega, s) = sI + \int_{0}^{s} (s-u) S_{i}^{(2)}(\omega,u)\
du = sI + \int_{0}^{s} (s-u) A_{i}(\omega,u) S_{i}(\omega,u)\
du.\nonumber
\end{equation}
Now using Theorem \ref{t.cs-est},
\begin{align}
\left \Vert \int_{0}^{s} (s-u) A_{i}(\omega,u) S_{i}(\omega,u)\ du
\right \Vert & \leq\int_{0}^{s}
(s-u)\frac{K\parallel\Delta_{i}\omega
\parallel^{2}}{(\Delta_{i}s)^{2}} u\cosh(\sqrt K \parallel\Delta_{i}\omega \parallel)\ du\nonumber\\
&  = \Big(\int_{0}^{s} (s-u)u\ du
\Big)\frac{K\parallel\Delta_{i}\omega
\parallel^{2}}{(\Delta_{i}s)^{2}}
 \cosh(\sqrt K \parallel\Delta_{i}\omega
\parallel)\nonumber\\
&
=\frac{s^{3}}{6}\frac{K}{(\Delta_{i}s)^{2}}\parallel\Delta_{i}\omega
\parallel^{2} \cosh(\sqrt K \parallel\Delta_{i}\omega \parallel)\nonumber\\
&  \leq s K \parallel\Delta_{i}\omega \parallel^{2} \cosh(\sqrt K
\parallel \Delta_{i}\omega \parallel).\nonumber
\end{align}
Hence, if we choose an $\epsilon$ such that for $0 < x
\leq\epsilon$,
\begin{equation}
K x^{2}\cosh(\sqrt K x) < 1,\nonumber
\end{equation}
then
\begin{equation}
S_{i}(\omega, s) =  s \left( I + \frac{1}{s}\int _{0}^{s} (s-u)
A_{i}(\omega,u) S_{i}(\omega,u)\ du \right)\nonumber
\end{equation}
is invertible for $0< s \leq \Delta_is$.
\end{proof}

In order to compute $\lim_{n\rightarrow0}\rho_{n}$, we will first
derive another formula for $\rho_{n}$. Define a set
of tangent vectors $\{f_{i,a}(\omega,s)\}_{%
\genfrac{}{}{0pt}{}{i=1,2,\ldots,n}{a=1,2,\ldots,d}%
}$ on $T_{\omega}H_{{\mathcal{P}}}^{\epsilon}(\mathbb{R}^{d})$
such that $f_{i,a}(\omega,s)$ is the solution to Equation
(\ref{eq.Jacobi}) with the given initial conditions
\begin{equation}
f_{i,a}(\omega,0)=0,\nonumber
\end{equation}
for $j=1,\ldots,n$,
\begin{equation}
f_{i,a}^{\prime}(\omega,s_{j-1})=\left\{
\begin{array}
[c]{ll}%
e_{a}, & \hbox{$j = i$}\\
-F_{i}(\omega)\ e_{a}, & \hbox{$j = i+1$}\\
0, & \hbox{otherwise.}
\end{array}
\right. \nonumber
\end{equation}
where
\begin{equation}
F_{i}(\omega):=\left(  S_{i+1}(\omega,\Delta_{i+1}s)\right)
^{-1}C_{i+1}(\omega,\Delta
_{i+1}s)\ S_{i}(\omega,\Delta_{i}s), \label{eq.F_i}%
\end{equation}
where $S_{i}$ and $C_{i}$ are as in Definition \ref{defn.C&S}. By
Lemma \ref{lem.Sinv}, we can choose an $\epsilon$ such that
$F_{i}(\omega)$ is defined on
$H_{{\mathcal{P}}}^{\epsilon}(\mathbb{R}^{d})$. Therefore,
\[
f_{i,a}(\omega,s)=\left\{
\begin{array}
[c]{ll}%
S_{i}(\omega,s-s_{i-1})\ e_{a}, & \hbox{$s \in J_i$}\\
C_{i+1}(\omega,s-s_{i})S_{i}(\omega,\Delta_{i}s)\ e_{a}-S_{i+1}(\omega,s-s_{i})F_{i}%
(\omega)\ e_{a}, & \hbox{$s \in J_{i+1}$}\\
0, & \hbox{otherwise,}
\end{array}
\right.
\]
and hence
\[
f_{i,a}^{\prime}(\omega,s)=\left\{
\begin{array}
[c]{ll}%
F_{ii}(\omega,s)\ e_{a}, & \hbox{$s \in J_i$}\\
F_{i+1,i}(\omega,s)\ e_{a}, & \hbox{$s \in J_{i+1}$}\\
0, & \hbox{otherwise,}
\end{array}
\right.
\]
where for $i=1,\ldots,n$,
\begin{align}
F_{ii}(\omega,s)  &  =S_{i}^{\prime}(\omega,s-s_{i-1}),\nonumber\\
F_{i+1,i}(\omega,s)  &  =C_{i+1}^{\prime}(\omega,s-s_{i})S_{i}(\omega,\Delta_{i}%
s)\ -S_{i+1}^{\prime}(\omega,s-s_{i})F_{i}(\omega).\nonumber
\end{align}

\begin{rem}
For all $\sigma\in H_{{\mathcal{P}}}^{\epsilon}(M),$ the
vectors, $\{X^{f_{i,a}}\left(  \sigma\right)  \}_{%
\genfrac{}{}{0pt}{}{i=1,\ldots,n}{a=1,\ldots,d}%
}$ form a basis for $T_{\sigma}H_{{\mathcal{P}}}^{\epsilon}(M)$,
where $\sigma=\phi(\omega)$.
\end{rem}


At this point, we will now assume $\mathcal{P}_n=\left\{
s_{i}=\frac{i}{n}\right\}
_{i=0}^{n}$ and $\Delta=1/n$ throughout this section. We may now write%
\begin{equation}
f_{i,a}^{\prime}(\omega,s)=\left[  1_{J_{i}}\left(  s\right)
S_{i}^{\prime }(\omega,s-s_{i-1})+1_{J_{i+1}}\left(  s\right)
V_{i+1}^{\prime}\left(
\omega,s-s_{i}\right)  \right]  e_{a}\label{eq.fia}%
\end{equation}
where $V_{1}\equiv0\equiv V_{n+1}$ and for $2\leq i\leq n,$%
\[
V_{i}\left(  \omega,s\right)
:=C_{i}(\omega,s)S_{i-1}(\omega,\Delta)-S_{i}(\omega,s)F_{i-1}(\omega)
\]
and%
\begin{equation}
F_{i}(\omega):=S_{i+1}(\omega,\Delta)^{-1}C_{i+1}(\omega,\Delta)S_{i}(\omega,\Delta).\label{e.fib}%
\end{equation}
Observe that
\begin{equation}
V_{i}\left(  \omega,\Delta\right)  =C_{i}\left(
\omega,\Delta\right) S_{i-1}\left( \omega,\Delta\right)
-S_{i}\left( \omega,\Delta\right) S_{i}\left(
\omega,\Delta\right) ^{-1}C_{i}\left( \omega,\Delta\right)
S_{i-1}\left( \omega,\Delta\right)
=0\label{e.vd}%
\end{equation}
and that
\begin{equation}
V_{i}\left(  \omega,0\right)  =S_{i-1}\left(  \omega,\Delta\right)  .\label{e.v0}%
\end{equation}

\begin{lem}
\label{l.ff}Continuing the notation above, we have%
\begin{align*}
& \int_{0}^{1}\langle f_{i,a}^{\prime}(\omega,s),
f_{j,c}^{\prime}(\omega,s)\rangle\ ds \nonumber \\
& =\left\{
\begin{array}
[c]{ccc}%
\int_{0}^{\Delta}\langle S_{i}^{\prime}(\omega,s)e_{a},
V_{i}^{\prime}\left(
\omega,s\right)  e_{c}\rangle\ ds & \text{if}\  j=i-1\\
\int_{0}^{\Delta}\left[ \langle S_{i}^{\prime}(\omega,s)e_{a},
S_{i}^{\prime}\left( \omega,s\right)  e_{c}\rangle+ \langle
V_{i+1}^{\prime}\left( \omega,s\right)  e_{a},
V_{i+1}^{\prime}\left(  \omega,s\right)  e_{c} \rangle \right]\  ds & \text{if}\  j=i\\
\int_{0}^{\Delta}\langle V_{i+1}^{\prime}(\omega,s)e_{a},
S_{i+1}^{\prime}\left(
\omega,s\right)  e_{c}\rangle\ ds & \text{if}\  j=i+1\\
0 & \text{otherwise} &
\end{array}
\right. \\
&  =\left\{
\begin{array}
[c]{ccc}%
\int_{0}^{\Delta}\langle e_{a},
S_{i}^{\prime}(\omega,s)^{T}V_{i}^{\prime}\left(
\omega,s\right)  e_{c}\rangle\ ds & \text{if} \ j=i-1\\
\int_{0}^{\Delta}\left[ \langle e_{a},
S_{i}^{\prime}(\omega,s)^{T}S_{i}^{\prime }\left(  \omega,s\right)
e_{c}\rangle+ \langle e_{a}, V_{i+1}^{\prime}\left(
\omega,s\right)
^{T}V_{i+1}^{\prime}\left(  \omega,s\right)  e_{c} \rangle \right]\  ds & \text{if} \ j=i\\
\int_{0}^{\Delta}\langle e_{a},
V_{i+1}^{\prime}(\omega,s)^{T}S_{i+1}^{\prime}\left(
\omega,s\right)  e_{c}\rangle\ ds & \text{if} \ j=i+1\\
0 & \text{otherwise.} &
\end{array}
\right.
\end{align*}

\end{lem}

Let us define the block matrix function of $s\in\left[  0,\Delta\right]  $ by,%
\[
\mathcal{F}^n_{ij}\left(  \omega,s\right)
=\delta_{ij}S_{i}^{^{\prime}}\left( \omega,s\right)
+\delta_{i,j+1}V_{i}^{\prime}\left(  \omega,s\right)
\]
or equivalently as%
\begin{equation}
\mathcal{F}^n\left(  \omega,s\right)  :=\left[
\begin{array}
[c]{ccccc}%
S_{1}^{^{\prime}}\left(  \omega,s\right)  & 0 & \dots & 0 & 0\\
V_{2}^{\prime}\left(  \omega,s\right)  & S_{2}^{^{\prime}}\left(
\omega,s\right)  & 0 &
& 0\\
0 & V_{3}^{\prime}\left(  \omega,s\right)  &
S_{3}^{^{\prime}}\left( \omega,s\right)  &
\ddots & \vdots\\
\vdots & \ddots & \ddots & \ddots & 0\\
0 & \dots & 0 & V_{n}^{\prime}\left(  \omega,s\right)  &
S_{n}^{^{\prime}}\left( \omega,s\right)
\end{array}
\right]  . \label{e.fb}%
\end{equation}
where
\[
V_{i}^{\prime}\left(  \omega,s\right)
:=C_{i}^{\prime}(\omega,s)S_{i-1}(\omega,\Delta
)-S_{i}^{\prime}(\omega,s)F_{i-1}(\omega).
\]

\begin{rem}
When $\omega=0$ we have $S_{i}\left(  0,s\right)  =sI,$
$C_{i}\left( 0,s\right) =I,$ $F_{i}(0):=\Delta^{-1}I\Delta=I,$
$V_{i}\left( 0,s\right)  :=\Delta
I-sI,$ and $\mathcal{F}^n\left(  0,s\right)  =\mathcal{T}^n$ where $\mathcal{T}^n%
_{ij}=\left(  \delta_{ij}-\delta_{i,j+1}\right)  I,$ i.e.
\begin{equation}
\mathcal{F}^n\left(  0,s\right)  =\mathcal{T}^n:=\left[
\begin{array}
[c]{ccccc}%
I & 0 & \dots & 0 & 0\\
-I & I & 0 & \dots & 0\\
0 & -I & I & \ddots & \vdots\\
\vdots & \ddots & \ddots & \ddots & 0\\
0 & \dots & 0 & -I & I
\end{array}
\right]  . \label{e.ft}%
\end{equation}
It is also worth observing that%
\begin{equation}
\mathcal{F}^n\left(  \omega,0\right)  :=\left[
\begin{array}
[c]{ccccc}%
I & 0 & \dots & 0 & 0\\
-F_{1}(\omega) & I & 0 & \dots & 0\\
0 & -F_{2}(\omega) & I & \ddots & \vdots\\
\vdots & \ddots & \ddots & \ddots & 0\\
0 & \dots & 0 & -F_{n-1}(\omega) & I
\end{array}
\right]  , \label{e.f0}%
\end{equation}
or equivalently that
\begin{equation}
\mathcal{F}^n_{ki}\left(  \omega,0\right)
=\delta_{ki}I-\delta_{i,k-1}F_{i}\left(
\omega\right)  . \label{e.f1}%
\end{equation}

\end{rem}

\begin{thm}\label{thm.rho}
Let $\mathcal{F}^n\left(  s\right)
_{ij}=\delta_{ij}S_{i}^{^{\prime}}\left( s\right)
+\delta_{i,j+1}V_{i}^{\prime}\left(  s\right)  $ {\rm (See
Equation
(\ref{e.fb}).)} Then%
\begin{equation}
\det\left(  \left\{  \int_{0}^{1}\left\langle
f_{i,a}^{\prime}(\omega,s), f_{j,c}^{\prime
}(\omega,s)\right\rangle\ ds\right\}  \right) =\det\left(
\int_{0}^{\Delta}\left( \mathcal{F}^n\left(
\omega,s\right)\right) ^{T}\mathcal{F}^n\left( \omega,s\right)
ds\right),\label{e.det}
\end{equation}
\begin{align}
&\det\left(  \left\{ \sum_{k=1}^{n}\left\langle
f_{i,a}^{\prime}(\omega,s_{k-1}+),
f_{j,c}^{\prime}(\omega,s_{k-1}+)\right\rangle
\Delta\right\}  \right)    =\Delta^{nd}, \label{e.rdet}%
\end{align}%
and hence
\begin{equation} \left(\rho_n \circ \phi \right)^2(\omega) =\det\left(
\frac{1}{\Delta}\int_{0}^{\Delta }(\mathcal{F}^n\left(
\omega,s\right))^{T}\mathcal{F}^n\left( \omega,s\right)  ds\right)
. \label{eq.det}%
\end{equation}
\end{thm}

\begin{proof}
Since%
\[
(\mathcal{F}^n_{ij}\left(  \omega,s\right))  ^{T}  =\delta_{ij}%
S_{i}^{^{\prime}}\left(  \omega,s\right)
^{T}+\delta_{j,i+1}V_{j}^{\prime}\left( \omega,s\right)  ^{T},
\]%
\begin{align*}
&\left[ ( \mathcal{F}^n\left(  \omega,s\right) )
^{T}\mathcal{F}^n\left( \omega,s\right) \right]  _{ij}\nonumber \\
&= \sum_{k}\left[ \delta_{ik}S_{i}^{^{\prime}}\left(
\omega,s\right) ^{T}+\delta_{k,i+1}V_{k}^{\prime}\left(
\omega,s\right) ^{T}\right]
\left[  \delta_{kj}S_{k}^{^{\prime}}\left(  \omega,s\right)  +\delta_{k,j+1}%
V_{k}^{\prime}\left(  \omega,s\right)  \right]  \\
&=   \delta_{ij}S_{i}^{^{\prime}}\left(  \omega,s\right)  ^{T}S_{i}^{^{\prime}%
}\left(  \omega,s\right)  +\delta_{i,j+1}S_{i}^{^{\prime}}\left(
\omega,s\right)
^{T}V_{k}^{\prime}\left(  \omega,s\right)  +\delta_{i+1,j}V_{i+1}^{\prime}\left(  \omega,s\right)  ^{T}S_{i+1}^{^{\prime}%
}\left(  \omega,s\right)  +\delta_{ij}V_{i+1}^{\prime}\left(
\omega,s\right)
^{T}V_{i+1}^{\prime}\left(  \omega,s\right)  \\
 &=  \delta_{ij}\left[  S_{i}^{^{\prime}}\left(  \omega,s\right)  ^{T}%
S_{i}^{^{\prime}}\left(  \omega,s\right)  +V_{i+1}^{\prime}\left(
\omega,s\right) ^{T}V_{i+1}^{\prime}\left(  \omega,s\right)
\right]    +\delta_{i,j+1}S_{i}^{^{\prime}}\left(  \omega,s\right)
^{T}V_{i}^{\prime }\left(  \omega,s\right)
+\delta_{i+1,j}V_{i+1}^{\prime}\left(  \omega,s\right)
^{T}S_{i+1}^{^{\prime}}\left(  \omega,s\right)  .
\end{align*}
So comparing with the results from Lemma \ref{l.ff} it follows
that
\begin{align}
& G^1\left(X^{f_{i,a}}, X^{f_{j,c}} \right)(\omega)  = \int_0^1 g
\left(\frac{\nabla X^{f_{i,a}}(\omega,s)}{ds}, \frac{\nabla
X^{f_{j,c}}(\omega,s)}{ds}
\right)\ ds \nonumber \\
& = \int_{0}^{1}\left\langle f_{i,a}^{\prime}(\omega,s),
f_{j,c}^{\prime}(\omega,s)\right\rangle\ ds
=\int_{0}%
^{1}\left\langle e_{a},\left[ ( \mathcal{F}^n\left(
\omega,s\right)) ^{T}\mathcal{F}^n\left( \omega,s\right)  \right]
_{ij}e_{c} \right\rangle\ ds, \nonumber
\end{align}
from which Equation (\ref{e.det}) follows.

In order to prove Equation (\ref{e.rdet}) we begin by observing that%
\begin{align*}
f_{i,a}^{\prime}(\omega,s_{k-1}+) &  =\left[  1_{J_{i}}\left(
s_{k-1}+\right)
S_{i}^{\prime}(\omega,s_{k-1}-s_{i-1})+1_{J_{i+1}}\left(
s_{k-1}+\right)
V_{i+1}^{\prime}\left(  \omega,s_{k-1}-s_{i}\right)  \right]  e_{a}\\
&  =\left[
\delta_{ik}S_{i}^{\prime}(\omega,0)+\delta_{i,k-1}V_{i+1}^{\prime
}\left(  \omega,0\right)  \right]  e_{a}\\
&  =\left[  \delta_{ik}I-\delta_{i,k-1}F_{i}\left(  \omega\right)
\right]
e_{a}=\mathcal{F}^n_{ik}\left(  \omega,0\right)  e_{a}%
\end{align*}
where the last equality follows from Equation (\ref{e.f1}). Hence
it follows that
\begin{align*}
& G^1_{\mathcal{P}} \left( X^{f_{i,a}}, X^{f_{j,c}}
\right)(\omega)
= \sum_{k=1}^n g \left( \frac{\nabla X^{f_{i,a}}(\omega,s_{k-1}+)}{ds}, \frac{\nabla X^{f_{j,c}}(\omega,s_{k-1}+)}{ds}\right)  \nonumber \\
& = \sum_{k=1}^{n}\langle f_{i,a}^{\prime}(\omega,s_{k-1}+), f_{j,c}^{\prime}%
(\omega,s_{k-1}+)\rangle \Delta =\Delta\sum_{k=1}^{n}\left\langle
\mathcal{F}^n_{ki}\left( \omega,0\right)
e_{a},\mathcal{F}^n_{kj}\left( \omega,0\right)
e_{c}\right\rangle\nonumber \\
&=\Delta\sum_{k=1}^{n}\left\langle e_{a}, (
\mathcal{F}^n_{ik}\left( \omega,0\right)) ^{T}
\mathcal{F}^n_{kj}\left( \omega,0\right)
e_{c}\right\rangle=\Delta\left\langle e_{a},\left[ (
\mathcal{F}^n\left( \omega,0\right)) ^{T}\mathcal{F}^n\left(
\omega,0\right) \right] _{ij}e_{c}\right\rangle
\end{align*}
and therefore %
\[
\det\left(  \left\{
\sum_{k=1}^{n}f_{i,a}^{\prime}(\omega,s_{k-1}+)\cdot
f_{j,c}^{\prime}(\omega,s_{k-1}+)\Delta\right\}  \right)
=\det\left( \Delta (\mathcal{F}^n\left(  \omega,0\right) )
^{T}\mathcal{F}^n\left( \omega,0\right) \right)  =\Delta^{nd}.
\]
Equation (\ref{eq.det}) now follows from Equation
(\ref{eq.dense}).
\end{proof}

\subsection{Some Identities}

\begin{defn}
For real square matrix functions, $A\left(  s\right)  $ and
$B\left(
s\right)  ,$ of $s\in\left[  0,\Delta\right]  ,$ let%
\[
\left\langle A\right\rangle
:=\frac{1}{\Delta}\int_{0}^{\Delta}A\left( s\right)  ds
\]
and
\begin{align*}
\mathrm{Cov}\left(  A,B\right)   &
=\frac{1}{\Delta}\int_{0}^{\Delta}A\left( s\right)  ^{T}B\left(
s\right)  ds-\left(  \frac{1}{\Delta}\int_{0}^{\Delta }A\left(
s\right)  ds\right)  ^{T}\left(  \frac{1}{\Delta}\int_{0}^{\Delta
}B\left(  s\right)  ds\right)  \\
&  =\left\langle A^{T}B\right\rangle -\left\langle A\right\rangle
^{T}\left\langle B\right\rangle .
\end{align*}
Notice that $\langle A \rangle$ and $\mathrm{Cov}\left( A,B\right)
$ is again a square matrix.
\end{defn}

The following proposition summarizes some basic and easily proved
properties of $\mathrm{Cov}\left(  A,B\right)  .$

\begin{prop}
\label{p.cov}The \textbf{covariance functional, }$\mathrm{Cov},$
has the following properties:

\begin{enumerate}
\item $\mathrm{Cov}\left(  A,B\right)  $ is bilinear in $A$ and
$B.$

\item $\mathrm{Cov}\left(  A,B\right)  $ may be computed as%
\begin{align*}
\mathrm{Cov}\left(  A,B\right)   &
=\frac{1}{\Delta}\int_{0}^{\Delta}\left[ A\left(  s\right)
-\left\langle A\right\rangle \right]  ^{T}\left[  B\left(
s\right)  -\left\langle B\right\rangle \right]  ds\\
&  =\left\langle \left[  A\left(  \cdot\right)  -\left\langle
A\right\rangle \right]  ^{T}\left[  B\left(  \cdot\right)
-\left\langle B\right\rangle \right]  \right\rangle .
\end{align*}

\item $\mathrm{Cov}\left(  A,B\right)  =0$ if either $A\left(
s\right)  $ or $B\left(  s\right)  $ is a constant function.

\item $\mathrm{Cov}\left(  A,A\right)  $ is always a symmetric
non-negative matrix.
\end{enumerate}
\end{prop}

\textbf{Note : }To simplify notation, for the rest of this section
we will typically be omitting the argument, $\omega,$ from the
expressions to follow.

\begin{defn}
Define $\mathcal{G}^n\left(  s\right)  :=\mathcal{F}^n\left(
s\right) -\mathcal{T}^n=\mathcal{F}^n\left(  s\right)
-\mathcal{F}^n\left( 0\right)  ,$
i.e.%
\begin{align}
\mathcal{G}^n_{ij}\left(  s\right)   &  =\delta_{ij}\left[  S_{i}^{^{\prime}%
}\left(  s\right)  -I\right]  +\delta_{i,j+1}\left[
V_{i}^{\prime}\left(
s\right)  +I\right]  \nonumber\\
&  =\delta_{ij}\left[  S_{i}^{^{\prime}}\left(  s\right) -I\right]
+\delta_{i,j+1}\left[
C_{i}^{\prime}(s)S_{i-1}(\Delta)-S_{i}^{\prime
}(s)F_{i-1}+I\right]  .\label{e.gij}%
\end{align}
Also let
\begin{equation}
\mathcal{Y}^n:=\left\langle \mathcal{G}^n\right\rangle
=\frac{1}{\Delta}\int
_{0}^{\Delta}\mathcal{G}^n\left(  s\right)  ds.\label{e.y}%
\end{equation}

\end{defn}

\begin{lem}
\label{l.FF}Let
$\mathcal{V}^n_{ij}=\delta_{ij}\frac{1}{\Delta}S_{i}\left(
\Delta\right)  ,$ i.e.
\begin{equation}
\mathcal{V}^n:=\frac{1}{\Delta}\left[
\begin{array}
[c]{cccc}%
S_{1}\left(  \Delta\right)   & 0 & \dots & 0\\
0 & S_{2}\left(  \Delta\right)   & \ddots & \vdots\\
\vdots & \ddots & \ddots & 0\\
0 & \dots & 0 & S_{n}\left(  \Delta\right)
\end{array}
\right]  \label{e.vv}%
\end{equation}
and $\mathcal{D}^n_{ij}:=\delta_{ij}\left[
\frac{1}{\Delta}S_{i}\left( \Delta\right)  -I\right]  ,$ i.e.
$\mathcal{D}^n=\mathcal{V}^n-\mathcal{I}^n.$ Then
$\mathcal{Y}^n=\mathcal{T}^n\mathcal{D}^n$ and%
\begin{align}
\left\langle (\mathcal{F}^n)^T\mathcal{F}^n\right\rangle &=\left(
\mathcal{T}^n+\mathcal{Y}^n\right)^T  \left(  \mathcal{T}^n%
+\mathcal{Y}^n\right)  +\mathrm{Cov}\left(
\mathcal{G}^n,\mathcal{G}^n\right)\nonumber \\
&=(\mathcal{V}^n)^T(\mathcal{T}^n)^T\mathcal{T}^n\mathcal{V}^n%
+\mathrm{Cov}\left(  \mathcal{G}^n,\mathcal{G}^n\right)  .
\end{align}

\end{lem}

\begin{proof}
By the fundamental theorem of calculus along with Equations
(\ref{e.vd}) and
(\ref{e.v0}) we have%
\begin{align*}
\mathcal{Y}^n_{ij} &  =\frac{1}{\Delta}\int_{0}^{\Delta}\left(
\delta _{ij}\left[  S_{i}^{^{\prime}}\left(  s\right)  -I\right]
+\delta
_{i,j+1}\left[  V_{i}^{\prime}\left(  s\right)  +I\right]  \right)  ds\\
&  =\frac{1}{\Delta}\left(  \delta_{ij}\left[  S_{i}\left(
\Delta\right) -I\Delta\right]  +\delta_{i,j+1}\left[  V_{i}\left(
\Delta\right)
-V_{i}\left(  0\right)  +\Delta I\right]  \right)  \\
&  =\delta_{ij}\left[  \frac{S_{i}\left(  \Delta\right)
}{\Delta}-I\right] -\delta_{i,j+1}\left[  \frac{S_{i-1}\left(
\Delta\right)  }{\Delta}-I\right] .
\end{align*}
On the other hand%
\begin{align*}
\left(  \mathcal{T}^n\mathcal{D}^n\right)  _{ij} &  =\sum_{k}\left[  \delta_{ik}%
-\delta_{i,k+1}\right]  I\delta_{kj}\left[
\frac{1}{\Delta}S_{k}\left(
\Delta\right)  -I\right]  \\
&  =\delta_{ij}\left[  \frac{1}{\Delta}S_{i}\left(  \Delta\right)
-I\right] -\delta_{i-1,j}\left[  \frac{1}{\Delta}S_{i-1}\left(
\Delta\right) -I\right]  =\mathcal{Y}^n_{ij}.
\end{align*}
The second assertion is a consequence of the following simple manipulations,%
\begin{align*}
\left\langle (\mathcal{F}^n)^T\mathcal{F}^n\right\rangle  &
=\left\langle \left(  \mathcal{T}^n+\mathcal{G}^n\right)^T  %
\left(  \mathcal{T}^n+\mathcal{G}^n\right)  \right\rangle \nonumber \\
&=(\mathcal{T}^n)%
^T\mathcal{T}^n+(\mathcal{T}^n)^T\mathcal{Y}^n%
+(\mathcal{Y}^n)^T\mathcal{T}^n+\left\langle (\mathcal{G}^n)^T
\mathcal{G}^n\right\rangle \\
&  =\left(  \mathcal{T}^n+\mathcal{Y}^n\right)^T \left(
\mathcal{T}^n+\mathcal{Y}^n\right)  +\left\langle (\mathcal{G}^n)^T%
\mathcal{G}^n\right\rangle -(\mathcal{Y}^n)^T\mathcal{Y}^n\\
&  =\left(  \mathcal{T}^n+\mathcal{Y}^n\right) ^T\left(
\mathcal{T}^n+\mathcal{Y}^n\right)  +\mathrm{Cov}\left(  \mathcal{G}^n%
,\mathcal{G}^n\right)  .
\end{align*}
This completes the proof since
\[
\mathcal{T}^n+\mathcal{Y}^n=\mathcal{T}^n\left(
\mathcal{I}^n+\mathcal{D}^n\right) =\mathcal{T}^n\mathcal{V}^n.
\]

\end{proof}

\begin{cor}
\label{c.m}Letting $\mathcal{M}^n=\mathrm{Cov}\left(  \mathcal{G}^n,\mathcal{G}^n%
\right)  ,$%
\begin{equation}
\mathcal{S}^n:=(\mathcal{T}^n)^{-1}=\left[
\begin{array}
[c]{cccc}%
I & 0 & \dots & 0\\
I & I & \ddots & \vdots\\
\vdots & \ddots & \ddots & 0\\
I & \dots & I & I
\end{array}
\right]  ,\label{e.tinv}%
\end{equation}
and $\mathcal{V}^n$ be as in Equation (\ref{e.vv}), dropping the
superscript $n$, we have
\[
\left\langle \mathcal{F}^T\mathcal{F}\right\rangle
=\mathcal{V}^T\mathcal{T}^T\left(
\mathcal{I}+\mathcal{S}^T\left(  \mathcal{V}^T%
\right)  ^{-1}\mathcal{M}\mathcal{V}^{-1}\mathcal{S}\right)
\mathcal{T}\mathcal{V}%
\]
and
\begin{equation}
\det\left(  \left\langle \mathcal{F}^T\mathcal{F}%
\right\rangle \right)  =\left[  \det\left(  \mathcal{V}\right)
\right]
^{2}\cdot\det\left(  \mathcal{I+S}^T\left(  \mathcal{V}%
^T\right)  ^{-1}\mathcal{MV}^{-1}\mathcal{S}\right)
.\label{e.det1}%
\end{equation}

\end{cor}

\subsection{The Key Determinant Formula\label{s.exp}}

Our next goal is to expand out $\mathcal{V}^n$ and $\left(  (\mathcal{V}^n)%
^T\right)  ^{-1}\mathcal{M}^n(\mathcal{V}^n)^{-1}$ in powers of
$\omega.$ It turns out that we need the expansion of
$\mathcal{V}^n$ to order
$\parallel \omega \parallel^{3}$ and $\left(  (\mathcal{V}^n)%
^T\right)  ^{-1}\mathcal{M}^n(\mathcal{V}^n)^{-1}$ to order
$\parallel \omega\parallel^{5}.$

\begin{notation}\label{n.A_i}
Recall the definition of $A_i(\omega,s)$ in Equation
(\ref{eq.ai}). We will write $A_i(0) = A_i(\omega,0)$ and also
$A_i(t) = A_i(\omega,t)$ to simplify the notation. And when we
write $y = O(x)$, we mean there exists some constant $C>0$
independent of $i$, $n$ and $\omega$ such that $\parallel y
\parallel \leq C\parallel x
\parallel$.
\end{notation}

\begin{prop}
\label{p.est}There exists a $C<\infty$ such that the following
estimates hold
for $0\leq s\leq\Delta;$%
\begin{align}
\left\Vert S_{i}^{\prime}\left(  s\right)  -\left(
I+\frac{1}{2}A_{i}\left(
0\right)  s^{2}\right)  \right\Vert  &  \leq C\parallel \Delta_{i}%
\omega\parallel ^{3},\label{eq.a1}\\
\left\Vert \frac{S_{i}\left(  s\right)  }{s}-\left(
I+\frac{1}{6}A_{i}\left(
0\right)  s^{2}\right)  \right\Vert  &  \leq C\left\Vert \Delta_{i}%
\omega\right\Vert ^{3},\label{eq.a2}\\
\left\Vert C_{i}\left(  s\right)  -\left( I+\frac{1}{2}A_{i}\left(
0\right) s^{2}\right)  \right\Vert  & \leq C\left\Vert
\Delta_{i}\omega\right\Vert
^{3}\ \text{and}\label{eq.a3}\\
\left\Vert C_{i}^{\prime}\left(  s\right)  -A_{i}\left(  0\right)
s\right\Vert  &  \leq Cs^{-1}\left\Vert \Delta_{i}\omega\right\Vert ^{3}.%
\label{eq.a4}%
\end{align}
In the sequel we will abbreviate these type of estimates by writing%
\begin{align}
S_{i}^{\prime}\left(  s\right)   &  =I+\frac{1}{2}A_{i}\left(
0\right)
s^{2}+O\left(  \left\Vert \Delta_{i}\omega\right\Vert ^{3}\right),  \label{eq.b1}\\
\frac{S_{i}\left(  s\right)  }{s} &  =I+\frac{1}{6}A_{i}\left(
0\right)
s^{2}+O\left(  \left\Vert \Delta_{i}\omega\right\Vert ^{3}\right),  \label{eq.b2}\\
C_{i}\left(  s\right)   &  =I+\frac{1}{2}A_{i}\left(  0\right)
s^{2}+O\left(
\left\Vert \Delta_{i}\omega\right\Vert ^{3}\right)\ \text{and}\label{eq.b3}\\
C_{i}^{\prime}\left(  s\right)   &  =s^{-1}\left[  A_{i}\left(
0\right) s^{2}+O\left(  \left\Vert \Delta_{i}\omega\right\Vert
^{3}\right)  \right]
.\label{eq.b4}%
\end{align}

\end{prop}

\begin{proof}
Let $(DA_i)(\omega'(s), \cdot, \cdot) :=
\frac{d}{ds}\Omega_{u(s_{i-1} + s)}$. Since $M$ is compact, there
exists a constant $C > 0$ such that $\sup_{u \in O(M)} \parallel
D\Omega_u
\parallel < C$. Then
\begin{eqnarray}
\parallel (DA_i)(\omega'(s), \omega'(s),
\cdot)\omega'(s) \parallel & \leq &
\frac{1}{\Delta^3}\left(\parallel
(DA_i)(\Delta_i\omega, \Delta_i\omega, \cdot)\Delta_i\omega \parallel \right) \nonumber \\
& = & O\left(\frac{1}{\Delta^3}\parallel \Delta_i\omega
\parallel^3\right). \nonumber
\end{eqnarray}
Thus%
\[
\int_{0}^{s}rA_{i}\left(  r\right)  dr=\int_{0}^{s}r\left[
A_{i}\left( 0\right)  +\int_{0}^{r}A_{i}^{\prime}\left( t\right)
dt\right] dr=\frac{1}{2}A_{i}\left(  0\right) s^{2}+O\left(
\left\Vert \Delta _{i}\omega\right\Vert ^{3}\right)
\]
and similarly that%
\begin{align*}
\frac{1}{s}\int_{0}^{s}\left(  s-r\right)  rA_{i}\left(  r\right)
dr & =\frac{1}{s}\int_{0}^{s}\left(  s-r\right)  rA_{i}\left(
0\right)
dr+O\left(  \left\Vert \Delta_{i}\omega\right\Vert ^{3}\right)  \\
&  =\frac{1}{6}A_{i}\left(  0\right)  s^{2}+O\left(  \left\Vert
\Delta _{i}\omega\right\Vert ^{3}\right)
\end{align*}
and finally%
\[
\int_{0}^{s}(s-u)A_{i}(u)du=\int_{0}^{s}(s-u)A_{i}(0)du+O\left(
\left\Vert \Delta_{i}\omega\right\Vert ^{3}\right)
=\frac{1}{2}A_{i}\left(  0\right) s^{2}+O\left(  \left\Vert
\Delta_{i}\omega\right\Vert ^{3}\right)  .
\]
Combining these results with the three estimates in Proposition
\ref{p.expand} with $\kappa:=K\left\Vert
\omega_{i}^{\prime}\right\Vert^2 $ proves
(\ref{eq.a1}) -- (\ref{eq.a3}). For Equation (\ref{eq.a4}), we have%
\begin{align*} C_{i}^{\prime}\left(  s\right)   &
=\int_{0}^{s}A_{i}\left( r\right) C_{i}\left(  r\right)
dr=\int_{0}^{s}\left[  A_{i}\left( 0\right)  +\int
_{0}^{r}A_{i}^{\prime}\left(  t\right) dt\right] \left(
I+O\left(  \left\Vert \Delta_{i}\omega\right\Vert^2 \right)  \right)  dr\\
&  =A_{i}\left(  0\right)  s+s^{-1}O\left(  \left\Vert
\Delta_{i}\omega\right\Vert ^{3}\right)
\end{align*}
as desired.
\end{proof}

\begin{cor}
\label{cor.v_est}With $\mathcal{V}^n=\mathcal{I}^n+\mathcal{D}^n$
as in Equation (\ref{e.vv}) we have
\begin{equation}
\mathcal{V}^n_{ij}=\delta_{ij}\left(  I+\frac{1}{6}A_{i}\left(
0\right) \Delta^{2}+ \eta_i(\Delta)  \right)
\label{e.uij}%
\end{equation}
where \begin{equation} \eta_i(s) := \frac{S_{i}\left(  s\right)
}{s}-\left( I+\frac{1}{6}A_{i}\left( 0\right)  s^{2}\right)
\label{e.eta} \end{equation} and $\parallel \eta_i(\Delta)
\parallel = O\left( \left\Vert \Delta_{i}\omega\right\Vert
^{3}\right)$.
\end{cor}

\begin{lem}
\label{lem.fib}The function, $F_{i}  $ in Equation (\ref{e.fib})
satisfies,
\begin{equation}
F_{i}=I+\frac{1}{6}A_{i}(0)\Delta^{2}+\frac{1}{3}A_{i+1}\left(
0\right) \Delta^{2}+O\left(  \left\Vert
\Delta_{i}\omega\right\Vert ^{3}\vee\left\Vert
\Delta_{i+1}\omega\right\Vert ^{3}\right)  .\label{eq.fib}%
\end{equation}

\end{lem}

\begin{proof}
In order to simplify notation, let $a_{i}:=A_{i}\left(  0\right)
\Delta^{2}$ and $\beta_{i} = \left\Vert
\Delta_{i}\omega\right\Vert ^{3}$.
Then%
\begin{align*}
 F_{i}
&  =S_{i+1}\left(  \Delta\right) ^{-1}C_{i+1}\left( \Delta\right)
\ S_{i}\left(  \Delta\right) =\left( \frac{S_{i+1}\left(
\Delta\right) }{\Delta}\right) ^{-1}C_{i+1}\left(  \Delta\right) \
\frac{S_{i}\left(
\Delta\right)  }{\Delta}\\
&  =\left(  I+\frac{1}{6}a_{i+1}+O\left(\beta_{i+1}\right) \right)
^{-1}\left( I+\frac{1}{2}a_{i+1}+O\left(\beta_{i+1}\right)  \right)  \left(  I+\frac{1}{6}%
a_{i}+O\left(\beta_{i}\right)  \right)  \\
&  =\left(  I-\frac{1}{6}a_{i+1}+O\left(\beta_{i+1}\right) \right)
\left(
I+\frac{1}{2}a_{i+1}+O\left(\beta_{i+1}\right)  \right)  \left(  I+\frac{1}{6}%
a_{i}+O\left(\beta_i\right)  \right)  \\
&  =I+\left(  \frac{1}{2}-\frac{1}{6}\right)  a_{i+1}+\frac{1}{6}%
a_{i}+O\left(  \left\Vert \Delta_{i}\omega\right\Vert ^{3}\right)
+O\left( \left\Vert \Delta_{i+1}\omega\right\Vert ^{3}\right)
\end{align*}
which is equivalent to Equation (\ref{eq.fib}).
\end{proof}

\begin{thm}
\label{t.h_exp}Let%
\[
\mathcal{H}^n_{ij}\left(  s\right)  :=\delta_{ij}A_{i}\left(
0\right) \frac{s^{2}}{2}+\delta_{i,j+1}\left[  A_{i}\left(
0\right) \left( s\Delta -
\frac{s^{2}}{2}-\frac{\Delta^{2}}{3}\right) -A_{i-1}\left(
0\right) \frac{\Delta^{2}}{6}\right]
\]
and $\Upsilon^n_{ij}\left(  s\right) :=\mathcal{G}^n_{ij}\left(
s\right) -\mathcal{H}^n_{ij}\left( s\right)  .$ Then $\left\Vert
\Upsilon^n_{ij}\left(
s\right)  \right\Vert =O\left(  \left\Vert \Delta_{i-1}\omega\right\Vert ^{3}%
\vee\left\Vert \Delta_{i}\omega\right\Vert ^{3}\right)  $ and
since $\Upsilon^n_{ij}=0$ unless $i\in\left\{  j,j+1\right\}  ,$
it follows that
\begin{equation}
\mathcal{G}^n\left(  s\right)  =\mathcal{H}^n\left(  s\right)
+\Upsilon^n\left( s\right)  \text{ and }\left\Vert
\Upsilon^n\right\Vert =O\left(  \bigvee_{i=1, \ldots, n}
\left\Vert \Delta_{i}\omega\right\Vert ^{3}\right)  .\label{eq.ghe}%
\end{equation}

\end{thm}

\begin{proof}
Let $\gamma_i =  \left\Vert \Delta_{i-1}\omega\right\Vert
^{3}\vee\left\Vert \Delta_{i}\omega\right\Vert ^{3}$. By
Proposition \ref{p.est} and Lemma \ref{lem.fib},
\begin{align}
C_{i}^{\prime}(s)S_{i-1}(\Delta) &  =\left[  A_{i}\left(  0\right)
s\Delta+O\left(  \left\Vert \Delta_{i}\omega\right\Vert
^{3}\right) \right]
\Delta^{-1}S_{i-1}(\Delta)\nonumber\\
&  =\left[  A_{i}\left(  0\right)  s\Delta+O\left(  \left\Vert
\Delta _{i}\omega\right\Vert ^{3}\right)  \right]  \left[
I+\frac{1}{6}A_{i-1}\left( 0\right)  \Delta^{2}+O\left( \left\Vert
\Delta_{i-1}\omega\right\Vert ^{3}\right)
\right]  \nonumber\\
&  =A_{i}\left(  0\right)  s\Delta+O\left(\gamma_i\right)  \label{eq.csi}%
\end{align}
and%
\begin{align}
& S_{i}^{\prime}(s)F_{i-1} \nonumber \\
&  =\left( I+\frac{1}{2}A_{i}\left( 0\right) s^{2}+O\left(
\left\Vert \Delta_{i}\omega\right\Vert ^{3}\right)  \right)
\left( I+\frac{1}{6}A_{i-1}(0)\Delta^{2}+\frac{1}{3}A_{i}\left(
0\right)
\Delta ^{2}+O\left(\gamma_i\right)  \right)  \nonumber\\
&  =I+A_{i-1}\left(  0\right)  \frac{\Delta^{2}}{6}+A_{i}\left(
0\right) \left(  \frac{s^{2}}{2}+\frac{\Delta^{2}}{3}\right)
+O\left(\gamma_i\right) .\label{eq.sfi}%
\end{align}
Combining the last two equations shows,%
\begin{align}
V_{i}^{\prime}\left(  s\right)    &
=C_{i}^{\prime}(s)S_{i-1}(\Delta
)+I-S_{i}^{\prime}(s)F_{i-1}\nonumber\\
& =A_{i}\left(  0\right)  \left( s\Delta
-\frac{s^{2}}{2}-\frac{\Delta^{2}}{3}\right) -A_{i-1}\left(
0\right) \frac{\Delta^{2}}{6}+O\left(  \left\Vert
\Delta_{i-1}\omega\right\Vert ^{3}\vee\left\Vert
\Delta_{i}\omega\right\Vert
^{3}\right)  .\label{eq.vi}%
\end{align}
This equation along with Equations (\ref{e.gij}) and (\ref{eq.b1}) shows%
\begin{align*}
\mathcal{G}^n_{ij}\left(  s\right)    & =\delta_{ij}\left[  S_{i}^{^{\prime}%
}\left(  s\right)  -I\right]  +\delta_{i,j+1}\left[
V_{i}^{\prime}\left(
s\right)  +I\right]  \\
& =\delta_{ij}\left[  A_{i}\left(  0\right)
\frac{s^{2}}{2}+O\left(
\left\Vert \Delta_{i}\omega\right\Vert ^{3}\right)  \right]  \\
& +\delta_{i,j+1}\left[  A_{i}\left(  0\right)  \left( s\Delta - \frac{s^{2}}{2}%
-\frac{\Delta^{2}}{3}\right)  -A_{i-1}\left(  0\right)  \frac{\Delta^{2}}%
{6}+O\left(  \left\Vert \Delta_{i-1}\omega\right\Vert
^{3}\vee\left\Vert \Delta
_{i}\omega\right\Vert ^{3}\right)  \right]  \\
& =\mathcal{H}^n_{ij}\left(  s\right)  +O\left(  \left\Vert \Delta
_{i-1}\omega\right\Vert ^{3}\vee\left\Vert
\Delta_{i}\omega\right\Vert ^{3}\right)  .
\end{align*}
\end{proof}

\begin{thm}
\label{t.m_est}The matrix $\mathcal{M}^n$ of Corollary \ref{c.m} satisfies, $\mathcal{M}^n=\mathcal{C}^n%
+\mathcal{\tilde{E}}^n$ where $\mathcal{\tilde{E}}^n$ is a
tri-block-diagonal matrix such that $\left\Vert
\mathcal{\tilde{E}}^n\right\Vert =O\left( \bigvee_{i=1, \ldots,
n}\parallel \Delta_i\omega \parallel ^{5}\right)  $ and
$\mathcal{C}^n$ is the non-negative
tri-block-diagonal matrix given by%
\begin{equation}
\mathcal{C}^n_{ij}:=  \delta_{ij}\frac{\Delta^{4}}{45}\left[  A_{i}%
^{2}\left(  0\right)  +A_{i+1}^{2}\left(  0\right)  \right]
+\frac{7}{360}\left[\delta _{i,j+1}A_{i}^{2}\left( 0\right)
+\delta_{i,j-1}A_{j}^{2}\left(  0\right)
\right]\Delta^{4},\label{e.cij}%
\end{equation}
where $A_{n+1}^2(0) := 0$.
Equivalently,%
\[
\mathcal{C}^n:=\left(
\begin{array}
[c]{ccccc}%
\frac{1}{45}(a_{1}^{2}+a_{2}^{2}) & \frac{7}{360}a_{2}^{2} & 0 & \cdots & 0\\
\frac{7}{360}a_{2}^{2} & \frac{1}{45}(a_{2}^{2}+a_{3}^{2}) & \frac{7}{360}a_{3}^{2} & \ddots & \vdots\\
0 & \frac{7}{360}a_{3}^{2} & \ddots & \ddots & 0\\
\vdots & \ddots & \ddots & \frac{1}{45}(a_{n-1}^{2}+a_{n}^{2}) & \frac{7}{360}a_{n}^{2}\\
0 & \cdots & 0 & \frac{7}{360}a_{n}^{2} & \frac{1}{45}a_{n}^{2}%
\end{array}
\right)
\]
where as before, $a_{i}:=A_{i}\left(  0\right)  \Delta^{2}.$
\end{thm}

\begin{proof}
Since $\mathcal{G}^n = \mathcal{H}^n + \Upsilon^n$ with
$\left\Vert \mathcal{H}^n\right\Vert =O\left( \bigvee_{i=1,
\ldots, n}\left\Vert \Delta_{i}\omega\right\Vert ^{2}\right)  $ and\\
$\left\Vert \Upsilon^n\right\Vert =O\left( \bigvee_{i=1, \ldots,
n}\left\Vert \Delta_{i}\omega\right\Vert ^{3}\right)  ,$ it
follows
that%
\begin{align*}
\mathcal{M}^n  & =\mathrm{Cov}\left(
\mathcal{G}^n,\mathcal{G}^n\right) =\mathrm{Cov}\left(
\mathcal{H}^n+\Upsilon^n,\mathcal{H}^n+\Upsilon^n\right)
\\
& =\mathrm{Cov}\left(  \mathcal{H}^n,\mathcal{H}^n\right)
+\mathrm{Cov}\left( \mathcal{H}^n,\Upsilon^n\right)
+\mathrm{Cov}\left(  \Upsilon^n
,\mathcal{H}^n\right)  +\mathrm{Cov}\left(  \Upsilon^n,\Upsilon^n\right)  \\
& =\mathrm{Cov}\left(  \mathcal{H}^n,\mathcal{H}^n\right)  +\mathcal{\tilde{E}}^n%
=\mathcal{C}^n+\mathcal{\tilde{E}}^n%
\end{align*}
where $\mathcal{C}^n:=\mathrm{Cov}\left(
\mathcal{H}^n,\mathcal{H}^n\right)  $ and
\[
\mathcal{\tilde{E}}^n=\mathrm{Cov}\left(  \mathcal{H}^n,\Upsilon^n\right)  +\mathrm{Cov}%
\left(  \Upsilon^n,\mathcal{H}^n\right)  +\mathrm{Cov}\left(
\Upsilon^n ,\Upsilon^n\right)  =O\left(  \bigvee_{i=1, \ldots,
n}\left\Vert \Delta_{i}\omega\right\Vert ^{5}\right)  .
\]
Since $\mathcal{H}^n=\mathcal{K}^n+\mathcal{L}^n$ where%
\[
\mathcal{K}^n_{ij}\left(  s\right)  := \delta_{ij} A_{i}\left(
0\right) \frac{s^{2}}{2} +
\delta_{i,j+1}A_i(0)\left[s\Delta-\frac{s^2}{2}\right]%
\]
and $\mathcal{L}^n$ is the constant (in $s)$ matrix,%
\[
\mathcal{L}^n_{ij}:=-\delta_{i,j+1}\left[  A_{i}\left(  0\right)  \frac{\Delta^{2}}%
{3}+A_{i-1}\left(  0\right)  \frac{\Delta^{2}}{6}\right]  ,
\]
it follows, using Proposition \ref{p.cov}, that%
\[
\mathcal{C}^n=\mathrm{Cov}\left(  \mathcal{H}^n,\mathcal{H}^n\right)  =\mathrm{Cov}%
\left(  \mathcal{K}^n,\mathcal{K}^n\right)  .
\]
Since%
\begin{align*}
\left\langle \mathcal{K}^n_{ij}\right\rangle & :=  \left[ \frac{1}{6} \delta_{ij}%
+\frac{1}{3}\delta_{i,j+1}\right]  A_{i}\left(  0\right) \Delta^2,
\end{align*}%
we have
\begin{align*}
\mathcal{\hat{K}}^n_{ij}\left(  s\right)   &
:=\mathcal{K}^n_{ij}\left(
s\right)  -\left\langle \mathcal{K}^n_{ij}\right\rangle \\
&  =  \delta_{ij}  A_{i}\left(  0\right)  \left(
\frac{s^{2}}{2}-\frac{\Delta^{2}}{6}\right)+\delta_{i,j+1}A_{i}\left(
0\right)\left(s\Delta - \frac{s^2}{2} - \frac{1}{3}\Delta^2\right)
.
\end{align*}
Let
\begin{align*}
\alpha & = \frac{s^{2}}{2}-\frac{\Delta^{2}}{6},\  \beta  =
s\Delta - \frac{s^2}{2} - \frac{1}{3}\Delta^2.
\end{align*}
By direct integration, we get
\begin{align*}
&\langle\alpha^2\rangle = \frac{1}{\Delta}\int_{0}^{\Delta}
\alpha ^{2}\ ds
= \frac{1}{\Delta}\int_{0}^{\Delta}\left(  \frac{s^{2}}{2}-\frac{\Delta^{2}}%
{6}\right)  ^{2}ds  =\frac{\Delta^{4}}{45}, & \\
&\langle \beta^2 \rangle = \frac{1}{\Delta}\int_{0}^{\Delta} \beta
^{2}\ ds = \frac{1}{\Delta}\int_{0}^{\Delta}\left( s\Delta -
\frac{s^2}{2} -
\frac{1}{3}\Delta^2\right)  ^{2}ds  = \frac{\Delta^{4}}{45} \ \text{and}& \\
&\langle \alpha \beta \rangle = \frac{1}{\Delta}\int_{0}^{\Delta} \alpha \beta\ ds = \frac{1}{\Delta}\int_{0}^{\Delta}\left(  \frac{s^{2}}{2}-\frac{\Delta^{2}}%
{6}\right)\left( s\Delta - \frac{s^2}{2} -
\frac{1}{3}\Delta^2\right)  ds  = \frac{7\Delta^4}{360}. &
\end{align*}
 Then we may conclude that
\begin{align*}
\mathcal{C}^n_{ij} &  =\mathrm{Cov}\left(
\mathcal{K}^n,\mathcal{K}^n\right)
_{ij}=\sum_{k}\left\langle \left(  (\mathcal{\hat{K}}^n)^T%
\right)  _{ik}\mathcal{\hat{K}}^n_{kj}\right\rangle
=\sum_{k}\left\langle
\left(  \mathcal{\hat{K}}^n_{ki}\right)^T\mathcal{\hat{K}%
}^n_{kj}\right\rangle \\
& =  \sum_{k} \Big\langle \left[ \alpha
\delta_{ki}+\beta\delta_{k,i+1}\right] A_{k}\left( 0\right) \left[
\alpha\delta_{kj}+\beta\delta_{k,j+1}\right] A_{k}\left(
0\right) \Big\rangle\\
&  =\frac{\Delta^{4}}{45}\sum_{k}\left[ \delta_{ki}\delta_{kj} +
\delta_{k,i+1}\delta_{k,j+1}\right] A_{k}^2\left( 0\right) +
\Delta^4\frac{7}{360}\sum_{k}\left[ \delta_{ki}\delta_{k,j+1}+
\delta_{k,i+1}\delta_{kj}\right] A_{k}^2\left(
0\right)\\
&  =\frac{\Delta^{4}}{45}\delta_{ij}\left(     A_{i}^{2}\left(
0\right)  +  A_{i+1}^{2}\left( 0\right) \right) +
\Delta^4\frac{7}{360}\left[ \delta_{i,j+1}+ \delta_{i+1,j}\right]
A_{i}^2\left( 0\right)
\end{align*}
which is equivalent to Equation (\ref{e.cij}).
\end{proof}

Putting together the previous estimates leads to the following key
determinant formula.

\begin{thm}
\label{t.ket}As above, let $\mathcal{S}^n=(\mathcal{T}^n)^{-1},$%
\begin{align*}
\mathcal{V}^n_{ij}  &
=\delta_{ij}+\mathcal{D}^n_{ij}=\delta_{ij}\frac{1}{\Delta
}S_{i}\left(  \Delta\right)  =\delta_{ij}\left(
I+\frac{1}{6}A_{i}\left( 0\right)  \Delta^{2}+\eta_i(\Delta) \right) ,\\
\mathcal{C}^n_{ij}  & =\frac{\Delta^{4}}{45} \delta_{ij}\left[  A_{i}%
^{2}\left(  0\right)  +A_{i+1}^{2}\left(  0\right)  \right] +
\frac{7}{360}\left[\delta _{i,j+1}A_{i}^{2}\left(  0\right)
+\delta_{i,j-1}A_{j}^{2}\left( 0\right)\right]\Delta^4
\end{align*}
where $\parallel \eta_i(\Delta) \parallel \leq C\parallel
\Delta_i\omega
\parallel^3$, $C$ is independent of $i$, $n$, $\omega$
and define%
\begin{equation}
\mathcal{U}^n:=(\mathcal{S}^n)^T\mathcal{C}^n\mathcal{S}^n.\label{eq.u}%
\end{equation}
Then there exists a tri-block-diagonal matrix $\mathcal{E%
}^n=\mathcal{E}^{n}\left(  \omega\right)  $ such that \beq
\left\Vert \mathcal{E}^n\right\Vert =O\left(  \bigvee_{i=1,\ldots
,n}\left\Vert \Delta
_{i}\omega\right\Vert ^{5}\right)  \nonumber \eeq and%
\begin{equation}
\det\left(  \left\langle (\mathcal{F}^n)^T\mathcal{F}%
\right\rangle \right)  =\left[  \det\left(  \mathcal{V}^n\right)
\right] ^{2}\cdot\det\left(  \mathcal{I}^n+\mathcal{U}^n\right)
\cdot \det\ \left(
\mathcal{I}^n+ \mathcal{X}^n\right)  ,\label{e.key}%
\end{equation}
where $\mathcal{X}^n := \left(\mathcal{I}^n+\mathcal{U}^n\right)^{-1}(\mathcal{S}^n)%
^T\mathcal{E}^n\mathcal{S}^n$.

\end{thm}

\begin{proof}
To ease the notation, we will drop the superscript $n$ in this
proof. From Equation (\ref{e.det1}) of Corollary \ref{c.m},
Corollary \ref{cor.v_est}, and Theorem \ref{t.m_est},
\begin{equation}
\det\left(  \left\langle \mathcal{F}^{T}\mathcal{F}%
\right\rangle \right)  =\left[  \det\left(  \mathcal{V}\right)
\right] ^{2}\cdot\det\left(  \mathcal{I}+\mathcal{S}^{T}\left(
\mathcal{V}^{T}\right)  ^{-1}\left(  \mathcal{C}%
+\mathcal{\tilde{E}}\right)  \mathcal{V}^{-1}\mathcal{S}\right)  .\label{e.det2}%
\end{equation}
Now write $\mathcal{V}^{-1}=\mathcal{I}+\Psi,$ where
\[
\Psi:=\sum_{n=1}^{\infty}\left(  -1\right)  ^{n}\mathcal{D}^{n}%
\]
is a block-diagonal matrix satisfying, $\left\Vert \Psi\right\Vert
=O\left( \bigvee_{i=1, \ldots, n}\left\Vert
\Delta_{i}\omega\right\Vert ^{2}\right)  .$ Hence for $\epsilon>0$
sufficiently small, if $\bigvee_{i=1, \ldots, n}\left\Vert \Delta
_{i}\omega\right\Vert ^{2}\leq\epsilon,$ then $\left\Vert \mathcal{V}%
^{-1}\right\Vert \leq2.$ Furthermore,%
\begin{align*}
\left(  \mathcal{V}^{T}\right)  ^{-1}\left(  \mathcal{C}%
+\mathcal{\tilde{E}}\right)  \mathcal{V}^{-1}  & =\left(  \mathcal{V}%
^{T}\right)  ^{-1}\mathcal{CV}^{-1}+\left(  \mathcal{V}%
^{T}\right)  ^{-1}\mathcal{\tilde{E}V}^{-1}\\
& =\left(  \mathcal{I}+\Psi^{T}\right)  \mathcal{C}\left(
\mathcal{I}+\Psi\right)  +\left( \mathcal{V}^{T}\right)
^{-1}\mathcal{\tilde{E}V}^{-1}=\mathcal{C}+\mathcal{E}%
\end{align*}
where $\mathcal{E}$ is the tri-block-diagonal matrix defined by%
\[
\mathcal{E}=\mathcal{C}\Psi+\Psi^{T}\mathcal{C}%
+\Psi^{T}\mathcal{C}\Psi+\left(  \mathcal{V}%
^{T}\right)  ^{-1}\mathcal{\tilde{E}V}^{-1}%
\]
and $\mathcal{E}$ satisfies the norm estimate, $\left\Vert
\mathcal{E}\right\Vert =O\left(  \bigvee_{i=1, \ldots, n}\parallel
\Delta_i\omega \parallel ^{5}\right)  .$ Putting these results
back into Equation
(\ref{e.det2}) shows%
\begin{align*}
\det\left(  \left\langle \mathcal{F}^{T}\mathcal{F}%
\right\rangle \right)    & =\left[  \det\left(  \mathcal{V}\right)
\right] ^{2}\cdot\det\left(  \mathcal{I}+\mathcal{S}^{T}\left[
\mathcal{C}+\mathcal{E}\right]  \mathcal{S}\right)  \\
& =\left[  \det\left(  \mathcal{V}\right)  \right]
^{2}\cdot\det\left(
\mathcal{I}+\mathcal{U}+\mathcal{S}^{T}\mathcal{E%
S}\right)  \\
& =\left[  \det\left(  \mathcal{V}\right)  \right]
^{2}\cdot\det\left( \left(  \mathcal{I}+\mathcal{U}\right)  \left(
\mathcal{I}+\left(
\mathcal{I}+\mathcal{U}\right)  ^{-1}\mathcal{S}^{T%
}\mathcal{ES}\right)  \right)
\end{align*}
from which the desired result follows.
\end{proof}

\section{Convergence of $\{\rho_{ n} \circ \phi \circ b_n\}_{n=1}^\infty$ in
$\mu$-measure}\label{s.rho}

Recall that $H_\cP^\epsilon(\mathbb{R}^d)$ was defined as \beq
H_\cP^\epsilon(\mathbb{R}^d) = \left\{ \omega \in
H_\cP(\mathbb{R}^d) :\ \parallel \Delta_i\omega \parallel <
\epsilon\ \forall i \right\} \nonumber \eeq where $\Delta_i\omega
= \omega(s_i) - \omega(s_{i-1})$. Note that $\mathcal{V}^n$ and
$\mathcal{X}^n$ are only defined on $H_\cP^\epsilon(\mathbb{R}^d)$
for some $\epsilon$ satisfying Lemma \ref{lem.Sinv}.

\begin{notation}
By abuse of notation, we will now write $A_i(\omega, 0) =
A_i(\omega)$ from now on. This should not be confused with the
notation $A_i(t)$ as described in Notation \ref{n.A_i}, where $t
\in [0, \frac{1}{n})$ and we suppressed the argument $\omega$.
\end{notation}

Unless stated otherwise, we will only consider equally spaced
partitions $\cP_n = \{0 < \frac{1}{n} < \cdots < \frac{n}{n}=1\}$.
By Theorem \ref{thm.rho} and \ref{t.ket}, $\rho_n $ has been
written as a product of 3 terms, namely $\left[\det
(\mathcal{V}^n)\right]^2$, $\det( {\mathcal{I}^{n}} +
{\mathcal{U}}^n )$ and $\det\left( {\mathcal{I}^{n}} +
\mathcal{X}^n \right)$ on $H_{\cP_n}^\epsilon(\mathbb{R}^d)$. We
will now show that the determinant given in Theorem \ref{t.ket}
has a limit as $| {\mathcal{P}}_n| \rightarrow0$. The limit for
each term will be computed in this order.

\par The following theorem is the Wong-Zakai type approximation
theorem for solutions to Stratonovich stochastic differential
equations. This theorem is a special case of Theorem 5.7.3 and
Example 5.7.4 in \cite{MR1472487}. Theorems of this type have a
long history starting with Wong and Zakai \cite{MR0183023,
MR0386015}. The following version maybe found in \cite{Driver&Hu}.

\begin{thm}
\label{thm.wongkai} Let $f : \mathbb{R}^{d} \times\mathbb{R}^{n}
\rightarrow End(\mathbb{R}^{d}, \mathbb{R}^{n})$ and $f_{0}:
\mathbb{R}^{d} \times \mathbb{R}^{n} \rightarrow\mathbb{R}^{n}$ be
twice differentiable with bounded continuous derivatives. Let
$\xi_{0} \in\mathbb{R}^{n}$ and ${\mathcal{P}}$ be any partition
of $[0,1]$. Further, let $b$ and $b_{ {\mathcal{P}}}$ be as in
Definition \ref{defn.projection} and $\xi_{ {\mathcal{P}}}(s)$
denote the solution to the ordinary differential equation
\begin{equation}
\xi^{\prime}_{ {\mathcal{P}}}(s) = f(\xi_{ {\mathcal{P}}}(s))
b^{\prime}_{ {\mathcal{P}}}(s) + f_{0}(\xi_{ {\mathcal{P}} }(s)),\
\xi_{ {\mathcal{P}}}(0) = \xi_{0}\nonumber
\end{equation}
and $\xi$ denote the solution to the Stratonovich stochastic
differential equation,
\begin{equation}
d\xi(s) = f(\xi(s))\delta b(s) + f_{0}(\xi(s)) ds,\ \xi(0) =
\xi_{0}.\nonumber
\end{equation}
Then for any $\alpha\in(0, \frac{1}{2}),\ p \in[1, \infty)$, there
is a constant $C(p, \alpha) < \infty$ depending only on $f, f_{0}$
and $M$ so that
\begin{equation}
{{\mathbb{E}}} [\sup_{s \leq1}
\parallel\xi_{ {\mathcal{P}}}(s) - \xi(s) \parallel^{p}] \leq C(p, \alpha)|
{\mathcal{P}}|^{\alpha p}.\nonumber
\end{equation}

\end{thm}

\begin{defn}
\label{defn.stratonovich}

\begin{enumerate}

\item Let $u$ be the solution to the Stratonovich stochastic
differential equation
\begin{equation}
\delta u = {\mathcal{H}}_{u} u \delta b,\ u(0) =
u_{0}.\label{e.par.trans}
\end{equation}
Hence $u$ maybe viewed as $\mu$-a.s. defined function from
$W(\mathbb{R}^{d}) \rightarrow W(O(M))$.

\item Let $\tilde{\phi} := \pi\circ u: W(\mathbb{R}^{d})
\rightarrow W(M)$. This map will be called the stochastic
development map.


\end{enumerate}
\end{defn}

The following facts will be stated without any proof. See
\cite{MR1194990}.

\begin{rem}\label{rem.s//trans}

\begin{enumerate}
\item\label{rem.s//trans.1} Using Theorem \ref{thm.wongkai}, one
may show that $\tilde{\phi}$ is a "stochastic extension" of
$\phi$, i.e., $\tilde{\phi} = \lim_{| {\mathcal{P}}|
\rightarrow0}\phi \circ b_\cP$. Moreover, the law of
$\tilde{\phi}$ (i.e., $\mu\tilde{\phi}^{-1}$) is the Wiener
measure $\nu$ on $W(M)$.

\item \label{rem.s//trans.2}One can prove that $u_s =
\widetilde{//}_s(\tilde{\phi})$ where $\widetilde{//}$ is
stochastic parallel transport defined in Definition
\ref{d.s//trans}.

\item \label{rem.s//trans.3}The law of $u$ under $\mu$ on
$W(\mathbb{R}^{d})$ and the law of $\widetilde{//}$ under $\nu$ on
$W(M)$ are equal. \end{enumerate}
\end{rem}

\begin{defn}\label{defn.walpha}
Let $0 < \alpha < \frac{1}{2}$ and define for $\cP_n = \{0 =
\frac{1}{n} < \cdots < \frac{n}{n} = 1 \}$, \beq
W_\alpha(\mathbb{R}^d) = \left\{ \omega \in W(\mathbb{R}^d) \Big|\
\exists\ n(\omega)\ s.t.\ \forall n \geq n(\omega),\ \bigvee_{i=1,
\ldots, n}
\parallel \Delta_i\omega \parallel \leq n^{-\alpha} < \epsilon \right\}
\nonumber \eeq for some $\epsilon$ satisfying Lemma \ref{lem.Sinv}
and $n_0 : W_\alpha(\mathbb{R}^d) \rightarrow \mathbb{N}$ by \beq
n_0(\omega) := \inf \left\{ n(\omega) \in \mathbb{N}\ | \forall n
\geq n(w),\ \bigvee_{i=1, \ldots, n}
\parallel \Delta_i\omega \parallel \leq n^{-\alpha} < \epsilon
\right\}.
\nonumber \eeq
\end{defn}

Note that \beq \left( W_\alpha(\mathbb{R}^d) \right)^c = \left\{
\omega \in W(\mathbb{R}^d) \Big| \bigvee_{i=1, \ldots, n}
\parallel \Delta_i\omega \parallel > n^{-\alpha}\ i.o. \right\}.
\nonumber \eeq

\begin{lem} \label{lem.walpha}
Let $0 < \alpha < \frac{1}{2}$, then \beq \mu\left(
W_\alpha(\mathbb{R}^d)\right) = 1. \nonumber \eeq
\end{lem}

\begin{proof}
Now
\begin{align}
& \mu\left( \left\{ \omega \in W(\mathbb{R}^d) \Big| \bigvee_{i=1,
 \ldots, n} \parallel \Delta_i\omega \parallel > n^{-\alpha}
\right\}\right)  \leq \frac{1}{n^{-\alpha p}} \bE \left[
\bigvee_{i=1, \ldots, n}
\parallel \Delta_ib \parallel^p \right] \nonumber \\
& \leq \frac{1}{n^{-\alpha p}} \bE \left[ \sum_{i=1}^n
\parallel \Delta_ib \parallel^p \right]
= \frac{C_p}{n^{-\alpha p}}\left(\frac{1}{n}
\right)^{\frac{p}{2}}n, \nonumber
\end{align}
where $C_p > 0$ is some constant. Therefore,
\begin{align}
&\sum_{n=1}^\infty \mu\left(\left\{ \omega \in W(\mathbb{R}^d)
\Big| \bigvee_{i=1, \ldots, n} \parallel \Delta_i\omega \parallel
>
n^{-\alpha} \right\} \right) \nonumber \\
& \leq C_p\sum_{n=1}^\infty
\left(\frac{1}{n}\right)^{\frac{p}{2}-1-\alpha p}  = C_p
\sum_{n=1}^\infty \left( \frac{1}{n} \right)^{(\frac{1}{2}-\alpha
)p - 1} < \infty, \nonumber
\end{align}
if $p$ is sufficiently large. Hence by Borel Cantelli Lemma, \beq
\mu\left( \left\{ \omega \in W(\mathbb{R}^d) \Big| \bigvee_{i=1,
\ldots, n} \parallel \Delta_i\omega \parallel > n^{-\alpha}\ i.o.
\right\}\right) = 0, \nonumber \eeq and hence the proof.
\end{proof}

Using Lemma \ref{lem.walpha}, we can now extend the definition of
$\det (\mathcal{V}^n) $ and $\det (\mathcal{I}^n + \mathcal{U}^n)$
to be a $\mu$-a.s. function. Let $\omega \in
W_\alpha(\mathbb{R}^d)$. By Definition \ref{defn.walpha}, for $n
\geq n_0(\omega)$, \beq \bigvee_{i=1, \ldots, n} \parallel
\Delta_i \omega \parallel < n^{-\alpha} < \epsilon. \nonumber \eeq
Therefore $b_n(\omega)
 \in H_{\cP_n}^\epsilon(\mathbb{R}^d)$ for all $n
\geq n_0(\omega)$.

\begin{defn}
Define \beq \widetilde{\det (\mathcal{V}^n \circ b_n )}= 1_{\{ n
\geq n_0\} }\cdot \det (\mathcal{V}^n \circ b_n ) \nonumber \eeq
and \beq \widetilde { \det (\mathcal{I}^n + \mathcal{X}^n)\circ
b_n } = 1_{\{ n \geq n_0\}} \cdot \det (\mathcal{I}^n +
\mathcal{X}^n)\circ b_n  . \nonumber \eeq
\end{defn}

\begin{notation}
\label{n.bm}Throughout the next few sections, let $b\left(
s\right) :W\left( \mathbb{R}^{d}\right) \rightarrow\mathbb{R}^{d}$
be the projection map, $b\left( s\right)  \left(  \omega\right)
=\omega\left( s\right)  $ for all $0\leq s\leq1$ and $\omega\in
W\left( \mathbb{R}^{d}\right)  .$ Note that when $W\left(
\mathbb{R}^{d}\right) $ is equipped with Wiener measure, $\mu,$
$\left\{  b\left( s\right) :0\leq s\leq1\right\} $ is a Brownian
motion. We
further let $\phi_{n}=\phi \circ b_{n}  $ and $u_{n}%
=//(\phi_{n})$.
\end{notation}

\subsection{Convergence of $\widetilde{\det (\mathcal{V}^n \circ b_n )}$ } \label{s.V}

\begin{lem}
\label{lem.trD}
\begin{equation}
\sum_{i=1}^{n} \left( {\mathrm{{tr}}} \ A_{i} (b_n)\Delta_{i}
s^{2} - \left(-\sum_{i = 1}^{n} Scal(\phi_n(s_{i-1}))\frac{1}{n}
\right) \right) \ \longrightarrow0 \nonumber
\end{equation}
$\mu-a.s.$ as $n \rightarrow\infty$.
\end{lem}

\begin{proof}
Note that \beq \Delta_ib_n = b_n(s_i) - b_n(s_{i-1}) = b(s_i) -
b(s_{i-1}) = \Delta_ib \nonumber \eeq and \beq A_i(b_n)\Delta_is^2
= \Omega_{u_n(s_{i-1})}\left((b_n)_i^\prime,\cdot
\right)(b_n)_i^\prime\Delta_is^2 =
\Omega_{u_n(s_{i-1})}(\Delta_ib_n, \cdot )\Delta_ib_n. \nonumber
\eeq Let $Ric_{u_n(s)} := \sum_{i=1}^d \Omega_{u_n(s)}\left(\cdot,
e_i\right)e_i$. Using the symmetry of $Ric$,
\begin{align}
\sum_{i=1}^{n} {\mathrm{{tr}}} \ A_{i} (b_n)\Delta_{i} s^{2} & =
\sum _{i=1}^{n} {\mathrm{{tr}}} \
\Omega_{u_n(s_{i-1})}(\Delta_ib_n, \cdot )\Delta_ib_n  =
-\sum_{i=1}^{n} \Big\langle Ric_{u_n(s_{i-1})}\Delta_{i}b_n,
\Delta_{i}b_n \Big\rangle \nonumber \\
&  = -\sum_{i=1}^{n} \Big\langle Ric_{u_n(s_{i-1})}\Delta_{i}b,
\Delta_{i}b \Big\rangle.\nonumber
\end{align}
By Ito's formula,
\begin{align}
&\Big\langle Ric_{u_n(s_{i-1})}\Delta_{i}b, \Delta_{i}b
\Big\rangle
\nonumber \\
& = \Big\langle Ric_{u_n(s_{i-1})}(b(s_{i}) - b(s_{i-1})),
b(s_{i}) - b(s_{i-1})
\Big\rangle\nonumber\\
&  = 2\int_{s_{i-1}}^{s_{i}} \Big\langle Ric_{u_n(s_{i-1})}((b(s)
- b(s_{i-1})),
db(s) \Big\rangle + \int_{s_{i-1}}^{s_{i}} {\mathrm{{tr}}} \ Ric_{u_n(s_{i-1}%
)}\ ds\nonumber\\
&  = 2\int_{s_{i-1}}^{s_{i}} \Big\langle Ric_{u_n(s_{i-1})}(b(s) -
b(s_{i-1})), db(s) \Big\rangle +
Scal(\phi_n(s_{i-1}))\Delta_{i}s.\nonumber
\end{align}
Thus
\begin{align}
\sum_{i=1}^{n} \left( {\mathrm{{tr}}} \ A_{i} (b_n)\Delta_{i}
s^{2} +
Scal(\phi_n(s_{i-1}))\frac{1}{n}\right)   &  = -2\int_{s_{i-1}}%
^{s_{i}} \Big\langle Ric_{u_n(s_{i-1})}(b(s) - b(s_{i-1})), db(s)
\Big\rangle.\nonumber
\end{align}
Define
\begin{align}
\xi_{n}  &  = 2\sum_{i=1}^{n} \int_{s_{i-1}}^{s_{i}}
\Big\langle Ric_{u_n(s_{i-1})}(b(s) - b(s_{i-1})), db(s) \Big\rangle\nonumber\\
&  = 2\sum_{i=1}^{n} \left(  \int_{0}^{1} \Big\langle 1_{J_{i}}(s)Ric_{u_n(s_{i-1}%
)}(b(s) - b(s_{i-1})), db(s) \Big\rangle \right)\nonumber\\
&  = 2\int_{0}^{1} \left\langle M_{ n}(s),\ db(s)
\right\rangle,\nonumber
\end{align}
where $M_{ n}(s) = \sum_{i=1}^{n}
1_{J_{i}}(s)Ric_{u_n(s_{i-1})}(b(s) - b(s_{i-1}))$. To complete
the proof, it suffices to show that $\xi_n$ converges to 0
$\mu$-a.s.. We will make use of Burkholder's Inequality,
\begin{equation}
{{\mathbb{E}}} \Big[ \sup_{0 \leq t \leq T} |M_{t}|^{p} \Big] \leq
C {{\mathbb{E}}} \Big[ \langle M \rangle_{T}^{\frac{p}{2}}
\Big],\nonumber
\end{equation}
where C is some constant. Thus applying this with $T = 1$ and $p =
4$, we have
\begin{align}
{{\mathbb{E}}}\ |\xi_{n}|^{4}  &  \leq{{\mathbb{E}}} \left[ \Big(
C \int _{0}^{1} \parallel M_{ n}(s)\parallel^{2}\ ds \Big)^{2}
\right]\nonumber\\
&  \leq {{\mathbb{E}}} \left[ C^{2} \int_{0}^{1} \parallel M_{ n}(s)%
\parallel^{4}\ ds \right]\ (By\ Jensen^{\prime}s\ Inequality)\nonumber\\
&  = {{\mathbb{E}}} \left[ C^{2} \int_{0}^{1} \sum_{i=1}^{n}
1_{J_{i}}(s)\parallel
Ric_{u_n(s_{i-1})}(b(s) - b(s_{i-1})) \parallel^{4}\ ds \right]\nonumber\\
&  = O\left( \int_{0}^{1} \sum_{i=1}^{n} 1_{J_{i}}(s)
{{\mathbb{E}}}
\parallel(b(s) - b(s_{i-1})) \parallel^{4}\ ds \right)\nonumber\\
&  = O\left( \int_{0}^{1} \sum_{i=1}^{n} 1_{J_{i}}(s) (s -
s_{i-1})^{2}\ ds \right)= O \left(\frac{1}{n^{2}}\right).\nonumber
\end{align}
Hence $\sum_{n=1}^{\infty}\mathbb{E} | \xi_{n} |^{4} < \infty$ and
thus
\begin{equation}
\xi_{n} \longrightarrow0\ \mu-a.s..\nonumber
\end{equation}

\end{proof}

\begin{prop}
\beq \sum_{i = 1}^{n} Scal(\phi_n(s_{i-1}))\frac{1}{n}
\longrightarrow \int_0^1 Scal(\tilde{\phi}(s))\ ds \nonumber \eeq
$\mu$-a.s. as $n \rightarrow \infty$.
\end{prop}

\begin{proof}
Note that we can write
\begin{eqnarray}
\sum_{i = 1}^{n} Scal(\phi_n(s_{i-1}))\frac{1}{n}  & = & \int_0^1
\sum_{i=1}^n 1_{J_i}(s)Scal(\phi_n(s_{i-1}))\ ds. \nonumber
\end{eqnarray}
Since $\phi_n = \phi \circ b_n \rightarrow \tilde{\phi}$ in the
sup norm $\mu$-a.s. as $n \rightarrow \infty$ and $Scal$ is a
continuous function, thus
\begin{eqnarray}
\sum_{i=1}^n 1_{J_i}(s)Scal(\phi_n(s_{i-1})) & = & \sum_{i=1}^n
1_{J_i}(s)\left[Scal(\phi_n(s_{i-1})) -
Scal(\tilde{\phi}(s_{i-1}))  + Scal(\tilde{\phi}(s_{i-1})) \right] \nonumber \\
&   \longrightarrow & Scal(\tilde{\phi}(s)). \nonumber
\end{eqnarray}
Hence we can apply the dominated convergence theorem to obtain
\beq \lim_{n \rightarrow \infty}\int_0^1 \sum_{i=1}^n
1_{J_i}(s)Scal(\phi_n(s_{i-1}))\ ds = \int_0^1
Scal(\tilde{\phi}(s))\ ds \nonumber \eeq $\mu-a.s.$.
\end{proof}

\begin{lem}
\label{lem.prodscal}
\begin{equation}
\widetilde { \det (\mathcal{V}^n \circ b_n) } = 1_{\{n \geq n_o\}}
 \prod_{i=1}^{n}\ \det \left[ I + \frac{1}{6}A_{i}
(b_n)\Delta_{i} s^{2} + \eta_{i}(b_n) \right] \longrightarrow
e^{-\frac{1}{6}\int_{0}^{1} Scal(\tilde{\phi}(s))\ ds}\nonumber
\end{equation}
$\mu-a.s.$ as $n \rightarrow\infty$.
\end{lem}

\begin{proof}
It suffices to consider on $W_\alpha(\mathbb{R}^d)$ with
$\frac{1}{3} < \alpha < \frac{1}{2}$, since $\mu \left(
W_\alpha(\mathbb{R}^d) \right) = 1$. For $n \geq n_0(\cdot)$,
write \beq \zeta_i^n = \frac{1}{6}A_{i} (b_n)\Delta_{i} s^{2} +
\eta_{i}(b_n) \nonumber \eeq where $\eta_i(b_n)$ was defined in
Equation (\ref{e.eta}) and $\parallel \eta_i(b_n) \parallel =
O\left( \bigvee_{i=1, \ldots, n} \parallel\Delta_{i}b \parallel^3
\right)$. Now using the perturbation formula in Equation
(\ref{eq.detformula}) of the Appendix with $r = 2$,
\begin{align}
\det \left[  I + \zeta_i^n \right] &  = \exp\left[ {\mathrm{{tr}}}
\ \zeta_i^n
+ R_2(\zeta_i^n) \right] \nonumber \\
&= \exp\left[\frac{1}{6}\tr\ A_i(b_n)\Delta_is^2 + \psi_i^n
\right] \nonumber
\end{align}
where \begin{align}\psi_{i}^n &:= {\mathrm{{tr}}} \  \eta_{i}(b_n)
+ R_2(\zeta_i^n) = {\mathrm{{tr}}} \  \eta_{i}(b_n) +
\sum_{k=2}^{\infty}(-1)^{k+1} {\mathrm{{tr}}} \ \left(
\frac{1}{6}A_{i} (b_n)\Delta_{i} s^{2} + \eta_{i}(b_n) \right)
^{k}.\nonumber
\end{align}
Using Equation (\ref{e.bound}),
\begin{equation}
\left| R_2(\zeta_i^n) \right| \leq \frac{d\parallel \zeta_i^n
\parallel^{2}}{1 - \parallel \zeta_i^n \parallel} = O\left( \bigvee_{i=1,
\ldots, n} \parallel\Delta_{i}b \parallel^4 \right) \nonumber
\end{equation}
and hence
\begin{align}
\left|\psi_i^n \right| & \leq |\tr\ \eta_i(b_n)| + \left|
R_2(\zeta_i^n) \right| = O\left(\bigvee _{i=1, \ldots, n}
\parallel\Delta_{i}b \parallel^{3}\right).\nonumber
\end{align}
Since we chose $\alpha > \frac{1}{3}$, on $
W_\alpha(\mathbb{R}^d)$,
\begin{align}
\left \Vert \sum_{i=1}^{n} \psi_{i}^n \right \Vert  &  =
O\left(n\bigvee_{i=1,
\ldots, n} \parallel\Delta_{i}b \parallel^{3} \right)\nonumber\\
&  = O\left(n \cdot n^{-3\alpha} \right) = O\left(n^{1-3\alpha}
\right) \nonumber \\
&  \longrightarrow0\ \nonumber
\end{align}
as $n \rightarrow\infty$. Together with Lemma \ref{lem.trD},
\begin{align}
&  1_{\{n \geq n_0\}}\prod_{i=1}^{n}\ \det \left[  I +
\frac{1}{6}A_{i} (b_n)\Delta_{i} s^{2} +
\eta_{i}(b_n) \right]  - \exp \left[{-\frac{1}{6}\sum_{i = 1}^{n} Scal(\phi_n(s_{i-1}))\frac{1}{n} } \right]\nonumber\\
&  = 1_{\{n \geq n_0\}}\exp \left[{ \frac{1}{6}\sum_{i=1}^{n}
\left( {\mathrm{{tr}}} \ A_{i} (b_n)\Delta_{i}
s^{2} + \psi_{i}^n \right) } \right] - \exp\left[{-\frac{1}{6}\sum_{i = 1}^{n} Scal(\phi_n(s_{i-1}))\frac{1}{n}} \right]\nonumber\\
&  = e^{-\frac{1}{6}\sum_{i = 1}^{n}
Scal(\phi_n(s_{i-1}))\frac{1}{n}} \left\{ 1_{\{n \geq
n_0\}}\exp\left[ \frac{1}{6}\sum_{i=1}^{n}\left( {\mathrm{{tr}}} \
A_{i} (b_n)\Delta_{i} s^{2}  + Scal(\phi_n(s_{i-1}))\frac{1}{n}
\right) + \sum_{i=1}^{n} \psi_{i}^n
\right] - 1 \right\} \nonumber\\
&  \longrightarrow0\ \nonumber
\end{align}
$\mu$-a.s. as $n \rightarrow\infty$. Finally,
\begin{align}
&  1_{\{n \geq n_0\}}\prod_{i=1}^{n}\ \det \left[  I +
\frac{1}{6}A_{i} (b_n)\Delta_{i} s^{2} +
\eta_{i}(b_n) \right]  - e^{-\frac{1}{6}\int_{0}^{1} Scal(\tilde{\phi}(s))\ ds}\nonumber\\
&  = 1_{\{n \geq n_0\}}\prod_{i=1}^{n}\ \det \left[  I +
\frac{1}{6}A_{i} (b_n)\Delta_{i} s^{2} +
\eta_{i}(b_n) \right]  - e^{-\frac{1}{6}\sum_{i = 1}^{n} Scal(\phi_n(s_{i-1}))\frac{1}{n}} \nonumber \\
& + e^{-\frac{1}{6}\sum_{i = 1}^{n} Scal(\phi_n(s_{i-1}))\frac{1}%
{n}} - e^{-\frac{1}{6}\int_{0}^{1} Scal(\tilde{\phi}(s))\ ds}\nonumber\\
&  \longrightarrow0\ \nonumber
\end{align}
$\mu$-a.s. as $n \rightarrow\infty$.
\end{proof}

\subsection{Convergence of $\det(\mathcal{I}^n+{\mathcal{U}^{n}}) \circ b_n$}

\label{s.cU}

Recall, from Equation (\ref{eq.u}), that $\mathcal{U}^{n}:=(\mathcal{S}^{n}%
)^{T}\mathcal{C}^{n}\mathcal{S}^{n}$ where%
\begin{equation}
{\mathcal{C}_{ij}^{n}}\circ
b_n=\delta_{ij}\left(\frac{1}{45}\left(A_{i}^{2}(b_{n})\Delta
_{i}s^{4}+A_{i+1}^{2}(b_{n})\Delta_{i+1}s^{4}\right)\right)+1_{\{|j-i|=1\}}\frac
{7}{360}A_{i\vee j}^{2}(b_{n})\Delta_{i\vee j}s^{4}\nonumber
\end{equation}
with $A_{n+1}^{2}(b_{n})\Delta_{n+1}s^{4}:=0$ as in Equation
(\ref{e.cij}).  In order to compute the limit of
$\det(\mathcal{I}^n+{\mathcal{U}^{n}})$ as $n\rightarrow
\infty,$ it will be necessary to compute%
\begin{equation}
{\mathrm{{Tr}}}\ \left(  [{\mathcal{U}^{n}}]^{k}\right)  ={\mathrm{{Tr}}%
}\ \left(
[(\mathcal{S}^{n})^{T}\mathcal{C}^{n}\mathcal{S}^{n}]^{k}\right)
={\mathrm{{Tr}}}\ [({\mathcal{S}^{n}})^{T}{\mathcal{C}^{n}}({\mathcal{B}^{n}%
}{\mathcal{C}^{n}})^{k-1}{\mathcal{S}^{n}}]={\mathrm{{Tr}}}\
\left( [{\mathcal{B}^{n}}{\mathcal{C}^{n}}]^{k}\right)  ,\nonumber
\end{equation}
where
${\mathcal{B}^{n}:}={\mathcal{S}^{n}}({\mathcal{S}^{n}})^{T}.$

\begin{lem}
\label{lem.cB}The matrix, ${\mathcal{B}^{n}:}={\mathcal{S}^{n}}({\mathcal{S}%
^{n}})^{T},$ is given by%
\begin{equation}
{\mathcal{B}^{n}_{lm}}=(\min lm)I\text{ for }l,m=1,2,\ldots,n.\label{eq.blm}%
\end{equation}
Moreover, $\mathcal{B}^{n}$ and $\mathcal{S}^{n}$ satisfy the norm estimates,%
\begin{equation}
\parallel{\mathcal{S}^{n}}\parallel=O(n)\text{ and }\parallel{\mathcal{B}^{n}%
}\parallel=O(n^{2}).\nonumber
\end{equation}

\end{lem}

\begin{proof}
By definition,
\begin{align}
{\mathcal{B}^{n}_{lm}} &  =\sum_{k=1}^{n}{\mathcal{S}^{n}_{lk}}({\mathcal{S}%
^{n}_{mk}})^{T}=\sum_{k=1}^{n}1_{l\geq k}1_{m\geq k}I\nonumber\\
&  =\sum_{k=1}^{n}1_{\{\min lm\geq k\}}I=(\min lm)I.\nonumber
\end{align}
Let $\lambda_{1}\geq\lambda_{2}\ldots\geq\lambda_{nd}$ be the
eigenvalues of ${\mathcal{B}^{n}}$. Since it is a positive
definite matrix, we have that
$\parallel{\mathcal{B}^{n}}\parallel=\lambda_{1}$. Therefore, we
have
\begin{equation}
\parallel{\mathcal{B}^{n}}\parallel=\lambda_{1}\leq\sum_{i=1}^{nd}\lambda
_{i}={\mathrm{{Tr}}}\ {\mathcal{B}^{n}}=\sum_{l=1}^{n}l\ {\mathrm{{tr}}%
}\ I=d\frac{n(n+1)}{2}=O(n^{2}).\nonumber
\end{equation}
Since
\begin{align}
\parallel{\mathcal{B}^{n}}\parallel &  =\sup_{\parallel v\parallel=1}%
\langle{\mathcal{B}^{n}}v,v\rangle=\sup_{\parallel v\parallel=1}%
\langle({\mathcal{S}^{n}})^{T}v,({\mathcal{S}^{n}})^{T}v\rangle\nonumber\\
&  =\left\Vert ({\mathcal{S}^{n}})^{T}\right\Vert ^{2}=\left\Vert
{\mathcal{S}^{n}}\right\Vert ^{2},\nonumber
\end{align}
it follows that $\parallel{\mathcal{S}^{n}}\parallel=O(n).$
\end{proof}

The following definition will be useful in describing the limiting
behavior of  $A_{m}^{2}(b_{n})\Delta_ms^{4}$ as
$n\rightarrow\infty.$

\begin{defn}\label{defn.Gammau}
Define $\Gamma:O(M)\longrightarrow\mathbb{R}^{d\times d}$ (the
$d\times d$
matrices) by%
\[
\Gamma(v)=\sum_{i,j=1}^{d}\Big(\Omega_{v}(e_{i},\Omega_{v}(e_{i},\cdot
)e_{j})e_{j}+\Omega_{v}(e_{i},\Omega_{v}(e_{j},\cdot)e_{i})e_{j}+\Omega
_{v}(e_{i},\Omega_{v}(e_{j},\cdot)e_{j})e_{i}\Big)
\]
where $\{e_{i}\}_{i=1,2,\ldots,d}$ is any orthonormal basis for
$T_{o}M$.
\end{defn}

\begin{notation}
For $a_{1},\ a_{2},\ a_{3},\ a_{4}\in\mathbb{R}^{d}$ and $1\leq
m\leq n,$ let $\widetilde{T}_{m}^{n}\left(  a_{1}\otimes
a_{2}\otimes a_{3}\otimes a_{4}\right)  :W\left(
\mathbb{R}^{d}\right)  \rightarrow\mathbb{R}^{d\times
d}$ be defined by%
\begin{equation}
\widetilde{T}_{m}^{n}\left(  a_{1}\otimes a_{2}\otimes
a_{3}\otimes
a_{4}\right)  :=\Omega_{u_{n}(s_{m-1})}(a_{1},\Omega_{u_{n}(s_{m-1})}%
(a_{2},\cdot)a_{3})a_{4}\in\mathbb{R}^{d\times d}.\nonumber
\end{equation}
If $\tau$ is a permutation of $\left\{  1,2,3,4\right\}  ,$ let%
\begin{equation}
(\tau\widetilde{T}_{m}^{n})(a_{1}\otimes a_{2}\otimes a_{3}\otimes
a_{4})=\widetilde{T}_{m}^{n}(a_{\tau(1)}\otimes a_{\tau(2)}\otimes
a_{\tau (3)}\otimes a_{\tau(4)})\nonumber
\end{equation}
and
\[
(\tau\widetilde{T}_{m,i}^{n})(a_{1},a_{2})=(\tau\widetilde{T}_{m}^{n}%
)(e_{i}\otimes e_{i}\otimes a_{1}\otimes a_{2})
\]
where $\{e_{i}\}_{i=1}^{d}$ is the standard orthonormal basis for
$\mathbb{R}^{d}.$
\end{notation}

Let $\underline{0}=0$ and for $t\in(s_{m-1},s_{m}],$ let
$\underline
{t}=s_{m-1},$ $\Delta:=s_{m}-s_{m-1}=1/n$,%
\begin{align*}
\Delta b(t) &  =b(t)-b(\underline{t}),\\
\Delta_{m}b &  =\Delta b\left(  s_{m}\right)  =b(s_{m})-b(s_{m-1}),\\
\Delta_{m}b^{i} &  =\Delta_{m}b\cdot e_{i},\text{ and}\\
\Delta b(t)^{\otimes^{3}}\otimes db(t) &  =\Delta
b(t)\otimes\Delta b(t)\otimes\Delta b(t)\otimes db(t).
\end{align*}

\begin{lem}
\label{l.ag-est}Let $\{e_{i}\}_{i=1,2,\ldots,d}$ be an orthonormal
basis for
$T_{o}M$, $\Delta = \frac{1}{n}$ and $m\in\left\{  1,2,\dots,n\right\}  .$ Then%
\begin{equation}
A_{m}^{2}(b_{n})\Delta^{4}-\Delta^2\Gamma(u_{n}(s_{m-1}))=\epsilon_{m}%
^{n}\label{eq.epsilonm}%
\end{equation}
where%
\begin{equation}
\epsilon_{m}^{n}:=\frac{1}{6}\int_{s_{m-1}}^{s_{m}}\sum_{\tau}%
(\tau\widetilde{T}_{m}^{n})(\Delta b(s)^{\otimes^{3}}\otimes
db(s))+\frac
{1}{2}\int_{s_{m-1}}^{s_{m}}(s_{m}-s)\sum_{i,\tau}(\tau\widetilde{T}_{m,i}%
^{n})(\Delta b(s)\otimes db(s)).\nonumber
\end{equation}
The above sums range over permutations, $\tau,$ of $\left\{
1,2,3,4\right\} $ and $i=1,2,\dots,d.$
\end{lem}

\begin{proof}
By definition, we have
\begin{align}
A_{m}^{2}(b_{n})\Delta^{4} &  =\Omega_{u_{n}(s_{m-1})}(\Delta_{m}%
b,\Omega_{u_{n}(s_{m-1})}(\Delta_{m}b,\cdot)\Delta_{m}b)\Delta_{m}b\nonumber\\
&  =\widetilde{T}_{m}^{n}\left(
\Delta_{m}b\otimes\Delta_{m}b\otimes
\Delta_{m}b\otimes\Delta_{m}b\right)  \nonumber\\
&  =\sum_{i,j,k,l}\widetilde{T}_{m}^{n}(e_{i}\otimes e_{j}\otimes
e_{k}\otimes
e_{l})\Delta_{m}b^{i}\Delta_{m}b^{j}\Delta_{m}b^{k}\Delta_{m}b^{l}.\nonumber
\end{align}
By Ito's formula,
\begin{align}
\widetilde{T}_{m}^{n}&(\Delta_{m}b\otimes \Delta_{m}b\otimes\Delta
_{m}b\otimes\Delta_{m}b)\nonumber\\
= &
\int_{s_{m-1}}^{s_{m}}\sum_{j=1}^{4}\widetilde{T}_{m}^{n}(\Delta
b(s)\otimes\ldots\underset{j^{\text{th}}\text{ - spot}}{\underbrace{db(s)}%
}\ldots\otimes\Delta b(s))+
\frac{1}{4}\int_{s_{m-1}}^{s_{m}}\sum_{i=1}^{d}\sum_{\tau}(\tau
\widetilde{T}_{m}^{n})(e_{i}\otimes e_{i}\otimes\Delta
b(s)\otimes\Delta
b(s))ds\nonumber\\
= &  \frac{1}{6}\int_{s_{m-1}}^{s_{m}}\sum_{\tau}(\tau\widetilde{T}_{m}%
^{n})(\Delta b(s)^{\otimes^{3}}\otimes db(s))+\frac{1}{4}\int_{s_{m-1}}%
^{s_{m}}\sum_{i,\tau}(\tau\widetilde{T}_{m,i}^{n})(\Delta
b(s)\otimes\Delta
b(s))\ ds.\label{eq.deltab^4}%
\end{align}
Now we claim that
\begin{equation}
\int_{s_{m-1}}^{s_{m}}\Delta b(s)\otimes\Delta b(s)\ ds=\int_{s_{m-1}}^{s_{m}%
}(s_{m}-t)\Delta b(t)\vee
db(t)+\frac{\Delta^{2}}{2}\sum_{i=1}^de_{i}\otimes
e_{i},\label{eq.db-db}%
\end{equation}
where $a_{1}\vee a_{2}=a_{1}\otimes a_{2}+a_{2}\otimes a_{1}.$ Let
\begin{align}
W_{s} &  =\int_{s_{m-1}}^{s}(s-t)\Delta b(t)\vee db(t)+\frac{(s-s_{m-1})^{2}%
}{2}\sum_{i=1}^{d}e_{i}\otimes e_{i}\label{eq.ws}\\
&  =s\int_{s_{m-1}}^{s}\Delta b(t)\vee
db(t)-\int_{s_{m-1}}^{s}t\Delta b(t)\vee
db(t)+\frac{(s-s_{m-1})^{2}}{2}\sum_{i=1}^{d}e_{i}\otimes
e_{i}\nonumber
\end{align}
and observe that $W_{s_{m}}$ is equal to the right side of Equation (\ref{eq.db-db}%
). Since,%
\begin{align}
dW_{s} &  =\left(  \int_{s_{m-1}}^{s}\Delta b(t)\vee db(t)\right)
ds+s\Delta
b(s)\vee db(s)-s\Delta b(s)\vee db(s) + \left(  (s-s_{m-1})\sum_{i=1}e_{i}\otimes e_{i}\right)  ds\nonumber\\
&  =\left(  \int_{s_{m-1}}^{s}\Delta b(t)\vee db(t)+(s-s_{m-1})\sum_{i=1}%
e_{i}\otimes e_{i}\right)  \ ds\nonumber\\
&  =(\Delta b(s)\otimes\Delta b(s))\ ds,\nonumber
\end{align}
it follows that
\[
W_{s_{m}}=\int_{s_{m-1}}^{s_{m}}\Delta b(s)\otimes\Delta b(s)\ ds
\]
which verifies Equation (\ref{eq.db-db}).

From Equation (\ref{eq.db-db}), we have%
\begin{align*}
&
\int_{s_{m-1}}^{s_{m}}\sum_{i,\tau}(\tau\widetilde{T}_{m,i}^{n})(\Delta
b(s)\otimes\Delta b(s))\ ds\\
&  =\int_{s_{m-1}}^{s_{m}}(s_{m}-s)\sum_{i,\tau}(\tau\widetilde{T}_{m,i}%
^{n})(\Delta b(s)\vee db(s))+\frac{\Delta^{2}}{2}\sum_{i,\tau}\sum_{j=1}%
^{d}(\tau\widetilde{T}_{m,i}^{n})(e_{j}\otimes e_{j})\\
&  =\int_{s_{m-1}}^{s_{m}}(s_{m}-s)\sum_{i,\tau}(\tau\widetilde{T}_{m,i}%
^{n})(\Delta b(s)\vee db(s))+\frac{\Delta^{2}}{2}\sum_{\tau}\sum_{i,j=1}%
^{d}(\tau\widetilde{T}_{m}^{n})(e_{i}\otimes e_{i}\otimes
e_{j}\otimes
e_{j})\\
&  =2\int_{s_{m-1}}^{s_{m}}(s_{m}-s)\sum_{i,\tau}(\tau\widetilde{T}_{m,i}%
^{n})(\Delta b(s)\otimes db(s))+\frac{\Delta^{2}}{2}8\Gamma
(u(s_{m-1}))\\
&  =2\int_{s_{m-1}}^{s_{m}}(s_{m}-s)\sum_{i,\tau}(\tau\widetilde{T}_{m,i}%
^{n})(\Delta b(s)\otimes db(s))+4\Delta^{2}\Gamma(u_{n}(s_{m-1})).
\end{align*}
When we sum over all permutations of $\{1,2,3,4\}$, we will have
$4! = 24$ terms. However, $\Gamma$ is a sum of 3 terms, hence we
will end up with $8$ copies of $\Gamma$, which explains the factor
8 in the second last equality. Combining this equation with
Equation (\ref{eq.deltab^4}) proves Equation (\ref{eq.epsilonm}).
\end{proof}

\begin{lem}
\label{lem.epsilon2}There exists $C\,<\infty$ such that for any
$\lambda>0,$ $n\in\mathbb{N},$ and $\left\{  c_{m}\right\}
_{m=1}^{n}\subset\left[ 0,\lambda\right]  ,$
\begin{equation}
{{\mathbb{E}}}\left\Vert
\sum_{m=1}^{n}c_{m}\epsilon_{m}^{n}\right\Vert ^{2}\leq
C\frac{\lambda^{2}}{n^{3}},\nonumber
\end{equation}
which as usual we abbreviate as
\[
{{\mathbb{E}}}\left\Vert
\sum_{m=1}^{n}c_{m}\epsilon_{m}^{n}\right\Vert ^{2}=O\left(
\frac{\lambda^{2}}{n^{3}}\right)  .
\]

\end{lem}

\begin{proof}
We may write
\[
\sum_{m=1}^{n}c_{m}\epsilon_{m}^{n}=\xi_{n}+\chi_{n}%
\]
where%
\begin{equation}
\xi_{n}=\frac{1}{6}\sum_{m=1}^{n}c_{m}\int_{s_{m-1}}^{s_{m}}\sum_{\tau}%
(\tau\widetilde{T}_{m}^{n})(\Delta b(s)^{\otimes^{3}}\otimes
db(s))\nonumber
\end{equation}
and
\begin{equation}
\chi_{n}=\frac{1}{2}\sum_{m=1}^{n}c_{m}\int_{s_{m-1}}^{s_{m}}(s_{m}%
-s)\sum_{i,\tau}(\tau\widetilde{T}_{m,i}^{n})(\Delta b(s)\otimes
db(s)).\nonumber
\end{equation}
Using the isometry property of the Ito integral, we have
\begin{align}
\mathbb{E}\left(  6\xi_{n}\right)  ^{2} &
={{\mathbb{E}}}\left\Vert
\sum_{m=1}^{n}c_{m}\int_{s_{m-1}}^{s_{m}}\sum_{\tau}(\tau\widetilde{T}_{m}%
^{n})(\Delta b(s)^{\otimes^{3}}\otimes db(s))\right\Vert ^{2}\nonumber\\
&  ={{\mathbb{E}}}\left\Vert \int_{0}^{1}\sum_{m=1}^{n}c_{m}1_{J_{m}}%
(s)\sum_{\tau}(\tau\widetilde{T}_{m}^{n})(\Delta
b(s)^{\otimes^{3}}\otimes
db(s))\right\Vert ^{2}\nonumber\\
&  =O\left(  {{\mathbb{E}}}\int_{0}^{1}\sum_{m=1}^{n}c_{m}^{2}1_{J_{m}%
}(s)\parallel\Delta b(s)\parallel^{6}\ ds\right)  \nonumber\\
&  =O\left(  \int_{0}^{1}\sum_{m=1}^{n}c_{m}^{2}1_{J_{m}}(s)(s-s_{m-1}%
)^{3}\ ds\right)  \nonumber\\
&  =O\left(  \sum_{m=1}^{n}c_{m}^{2}\frac{1}{n^{4}}\right)
=O\left( \frac{\lambda^{2}}{n^{3}}\right)  .\nonumber
\end{align}
Similarly,%
\begin{align}
\mathbb{E}\left(  2\chi_{n}\right)  ^{2} &=
{{\mathbb{E}}}\left\Vert
\sum_{m=1}^{n}c_{m}\int_{s_{m-1}}^{s_{m}}(s_{m}-s)\sum_{i,\tau}(\tau
\widetilde{T}_{m,i}^{n})(\Delta b(s)\otimes db(s))\Big)\right\Vert
^{2}\nonumber\\
&  ={{\mathbb{E}}}\left\Vert \int_{0}^{1}\sum_{m=1}^{n}c_{m}1_{J_{m}(s)}%
(s_{m}-s)\sum_{i=1}^{d}\sum_{\tau}(\tau\widetilde{T}_{m,i}^{n})(\Delta
b(s)\otimes db(s))\right\Vert ^{2}\nonumber\\
&  =O\left(  {{\mathbb{E}}}\int_{0}^{1}\sum_{m=1}^{n}c_{m}^{2}1_{J_{m}%
}(s)(s_{m}-s)^{2}\parallel\Delta b(s)\parallel^{2}\ ds\right)  \nonumber\\
&  =O\left(  \int_{0}^{1}\sum_{m=1}^{n}c_{m}^{2}1_{J_{m}}(s)(s_{m}%
-s)^{2}(s-s_{m-1})\ ds\right)  \nonumber\\
&  =O\left(  \sum_{m=1}^{n}c_{m}^{2}\frac{1}{n^{4}}\right)
=O\left( \frac{\lambda^{2}}{n^{3}}\right)  .\nonumber
\end{align}
Thus
\[
\left(  {{\mathbb{E}}}\left\Vert \sum_{m=1}^{n}c_{m}\epsilon_{m}%
^{n}\right\Vert ^{2}\right)  ^{\frac{1}{2}}\leq\left(  {{\mathbb{E}}%
}\left\Vert \xi_{n}\right\Vert ^{2}\right)  ^{\frac{1}{2}}+\left(
{{\mathbb{E}}}\left\Vert \chi_{n}\right\Vert ^{2}\right)  ^{\frac{1}{2}%
}=O\left(  \sqrt{\frac{\lambda^{2}}{n^{3}}}\right)
\]
and hence
\begin{equation}
{{\mathbb{E}}}\left\Vert
\sum_{m=1}^{n}c_{m}\epsilon_{m}^{n}\right\Vert ^{2}=O\left(
\frac{\lambda^{2}}{n^{3}}\right)  .\nonumber
\end{equation}

\end{proof}

\begin{defn}
\label{defn.lambda&phi} For $l,m=1,2,\ldots n,$ define
$\Lambda_{lm}^{n}$ and $\Phi_{lm}^{n}$ by
\begin{align}
\Lambda_{lm}^{n}  & =\Big(\min{l}{(m-1)}\Big)\frac{7}{360}\Gamma
(u_{n}(s_{m-1}))\frac{1}{n^{2}}+\Big(\min{l}{m}\Big)\frac{1}{45}\Gamma
(u_{n}(s_{m}))\frac{1}{n^{2}}1_{\{m<n\}}\nonumber\\
&  +\Big(\min{l}{m}\Big)\frac{1}{45}\Gamma(u_{n}(s_{m-1}))\frac{1}{n^{2}%
}+\Big(\min{l}{(m+1)}\Big)\frac{7}{360}\Gamma(u_{n}(s_{m}))\frac{1}{n^{2}%
}1_{\{m<n\}}\nonumber
\end{align}
and
\begin{align}
\Phi_{lm}^{n}  &  =\Big(\min{l}{(m-1)}\Big)\frac{7}{360}\epsilon_{m}%
^{n}+\Big(\min{l}{m}\Big)\frac{1}{45}\epsilon_{m+1}^{n}1_{\{m<n\}}
\nonumber \\%
&  +\Big(\min{l}{m}\Big)\frac{1}{45}\epsilon_{m}^{n}+\Big(\min{l}%
{(m+1)}\Big)\frac{7}{360}\epsilon_{m+1}^{n}1_{\{m<n\}}.\nonumber
\end{align}

\end{defn}

With this notation along with Equation (\ref{eq.blm}) and Equation (\ref{eq.epsilonm}), we have%
\begin{align}
({\mathcal{B}^{n}}{\mathcal{C}^{n}})_{lm}\circ b_n &  = \left(\sum_{k=1}^{n}{\mathcal{B}^{n}%
_{lk}}{\mathcal{C}^{n}_{km}}\right)\circ b_n=\left(\sum_{k=1}^{n}(\min lk){\mathcal{C}^{n}_{km}}%
\right)\circ b_n\nonumber\\
&  =\Lambda_{lm}^{n}+\Phi_{lm}^{n}.\label{e.BC}%
\end{align}

\begin{notation}
For any $k\in\mathbb{N}$ and $d\times d$ matrices
$\{M_{k}\}_{l=1}^{k}$, let
\begin{equation}
\Big[\prod_{l=1}^{k}\Big]M_{k}:=M_{1}M_{2}\cdots M_{k}.\nonumber
\end{equation}

\end{notation}

\begin{thm}
\label{thm.trf_k} For each $k=1,2,\dots$, define
$\gamma_{k}^{n}:W\left( \mathbb{R}^{d}\right)
\rightarrow\mathbb{R}$ by
\begin{equation}
\gamma_{k}^{n}:=\sum_{r_{k}=1}^{n}\cdots\sum_{r_{1}=1}^{n}{\mathrm{{tr}}%
}\ \Big[\prod_{l=1}^{k}\Big]\Lambda_{r_{i},r_{i+1}}^{n},\nonumber
\end{equation}
where $r_{k+1} = r_1.$ Then
\begin{equation}
\lim_{n\rightarrow\infty}\mathbb{E}\left\vert {\mathrm{{Tr}}}\ ([{\mathcal{U}%
^{n}}\circ b_n]^{k})-\gamma_{k}^{n}\right\vert =0\nonumber
\end{equation}
and in particular%
\[
{\mathrm{{Tr}}}\ \left([{\mathcal{U}^{n}}\circ
b_n]^{k}\right)-\gamma_{k}^{n}\longrightarrow0\text{ (in
}\mu\text{ -- measure) as }n\rightarrow\infty.
\]

\end{thm}

\begin{proof}
For a fixed $k\in\mathbb{N},$%
\begin{align}
{\mathrm{{Tr}}}\ \left( [{\mathcal{B}^{n}}{\mathcal{C}^{n}}\circ
b_n]^{k}\right)   &
=\left({\mathrm{{tr}}}\ \sum_{r_{k}=1}^{n}\cdots\sum_{r_{1}=1}^{n}({\mathcal{B}^{n}%
}{\mathcal{C}^{n}})_{r_{1},r_{2}}({\mathcal{B}^{n}}{\mathcal{C}^{n}}%
)_{r_{2},r_{3}}\cdots({\mathcal{B}^{n}}{\mathcal{C}^{n}})_{r_{k},r_{1}%
} \right) \circ b_n\nonumber\\
&  ={\mathrm{{tr}}}\
\sum_{r_{k}=1}^{n}\cdots\sum_{r_{1}=1}^{n}\Lambda
_{r_{1},r_{2}}^{n}\Lambda_{r_{2},r_{3}}^{n}\cdots\Lambda
_{r_{k},r_{1}}^{n}+\Xi^{n}\nonumber\\
&  =\gamma_{k}^{n}+\Xi^{n},\nonumber
\end{align}
where $\Xi^{n}$ consists of a finite sum of terms of the form
\begin{equation}
{\mathrm{{tr}}}\ \sum_{r_{k}=1}^{n}\cdots\sum_{r_{1}=1}^{n}c_{r_{1},r_{2}%
}c_{r_{2},r_{3}}\ldots
c_{r_{k},r_{1}}\theta_{r_{1}}\theta_{r_{2}}\ldots
\theta_{r_{k}},\nonumber
\end{equation}
with
\[
\theta_{r_{i}}\in\left\{
\frac{1}{n^{2}}\Gamma(u_{n}(s_{r_{i}-1})),\frac
{1}{n^{2}}\Gamma(u_{n}(s_{r_{i}})),\epsilon_{r_{i}}^{n},\epsilon_{r_{i}%
+1}^{n}\right\}  ,
\]
$1\leq c_{r_{i},r_{i+1}}\leq n$ for $i\leq n$ and for at least one
$r_{i},$ $\theta_{r_{i}}=\epsilon_{r_{i}}^{n}$. Since trace is
invariant under
cyclic permutation, we can assume that $\theta_{r_{k}}=\epsilon_{r_{k}}%
^{n}$. To finish the proof it suffices to show,
$\lim_{n\rightarrow\infty
}{{\mathbb{E}}}\parallel\Xi^{n}\parallel=0$
and for this it suffices to show,%
\begin{equation}
\lim_{n\rightarrow\infty}{{\mathbb{E}}}\left\Vert
\sum_{r_{k}=1}^{n}\cdots
\sum_{r_{1}=1}^{n}c_{r_{1},r_{2}}c_{r_{2},r_{3}}\ldots c_{r_{k},r_{1}}%
\theta_{r_{1}}\theta_{r_{2}}\ldots\theta_{r_{k}}\right\Vert
=0.\nonumber
\end{equation}
Let $r=(r_{1},\ldots r_{k})$ and write
\begin{equation}
c(r)=\prod_{i=1}^{k}c_{r_{i},r_{i+1}},\nonumber
\end{equation}
where by convention, $r_{k+1}:=r_{1}$. Then
\begin{align}
T_{n} &  :=\left(  {{\mathbb{E}}}\left\Vert
\sum_{r_{k}=1}^{n}\cdots
\sum_{r_{1}=1}^{n}c(r)\theta_{r_{1}}\theta_{r_{2}}\ldots\theta_{r_{k-1}%
}\epsilon_{r_{k}}\right\Vert \right)  ^{2}\nonumber\\
&  =\left(  {{\mathbb{E}}}\left\Vert \sum_{r_{k-1}=1}^{n}\cdots\sum_{r_{1}%
=1}^{n}\theta_{r_{1}}\theta_{r_{2}}\ldots\theta_{r_{k-1}}\sum_{r_{k}=1}%
^{n}c(r)\epsilon_{r_{k}}\right\Vert \right)  ^{2}\nonumber\\
&  \leq\left(  {{\mathbb{E}}}\sum_{r_{k-1}=1}^{n}\cdots\sum_{r_{1}=1}%
^{n}\left(  \left\Vert \theta_{r_{1}}\theta_{r_{2}}\ldots\theta_{r_{k-1}%
}\right\Vert \left\Vert
\sum_{r_{k}=1}^{n}c(r)\epsilon_{r_{k}}\right\Vert
\right)  \right)  ^{2}\nonumber\\
&  \leq{{\mathbb{E}}}\left(  \sum_{r_{1},\ldots,r_{k-1}}\prod_{i=1}%
^{k-1}\left\Vert \theta_{r_{i}}\right\Vert ^{2}\right)
{{\mathbb{E}}}\left(
\sum_{r_{1},\ldots,r_{k-1}}\left\Vert \sum_{r_{k}=1}^{n}c(r)\epsilon_{r_{k}%
}\right\Vert ^{2}\right)  ,\label{eq.EquationSb}%
\end{align}
where the last inequality is a consequence of the Cauchy Schwartz
inequality. Let
$C_{1}:=\sup_{v\in O(M)}\left\Vert \Gamma(v)\right\Vert <\infty,$ then%
\begin{equation}
\left\Vert \frac{1}{n^{2}}\Gamma(u(s_{i-1}))\right\Vert
^{2(k-1)}\leq C_{1}^{2(k-1)}\left(  \frac{1}{n}\right)
^{4(k-1)}=O\left(  \frac {1}{n^{4(k-1)}}\right)  .\nonumber
\end{equation}
Using Equation (\ref{eq.epsilonm}), we observe that there exists a
constant $C_{2}$ such that
\begin{equation}
{{\mathbb{E}}}\ \parallel\epsilon_{m}^{n}\parallel^{2(k-1)}\leq
{{\mathbb{E}}}\ \left[  C_{2}\left(\
\parallel\Delta_{m}b\parallel^{4}+\frac {1}{n^{2}}\right)\right]
^{2(k-1)}=O\left(  \frac{1}{n^{4(k-1)}}\right) .\nonumber
\end{equation}
Thus we can find a constant $C(k)$ such that
\begin{equation}
{{\mathbb{E}}}\ \parallel\theta_{r_{i}}\parallel^{2(k-1)}\leq
C\left(\frac {1}{n}\right)^{4(k-1)}.\nonumber
\end{equation}
By Holder's Inequality,
\begin{align}
&  {{\mathbb{E}}}\left(  \sum_{r_{k-1}=1}^{n}\cdots\sum_{r_{1}=1}^{n}%
\prod_{i=1}^{k-1}\left\Vert \theta_{r_{i}}\right\Vert ^{2}\right)
\leq
\sum_{r_{k-1}=1}^{n}\cdots\sum_{r_{1}=1}^{n}\prod_{i=1}^{k-1}\left(
{{\mathbb{E}}}\parallel\theta_{r_{i}}\parallel^{2(k-1)}\right)
^{\frac
{1}{k-1}}\nonumber\\
&  =\sum_{r_{k-1}=1}^{n}\cdots\sum_{r_{1}=1}^{n}C\Big(\frac{1}{n}%
\Big)^{4(k-1)}=C\left(\frac{1}{n}\right)^{3(k-1)}.\label{eq.equationTa}%
\end{align}
Note that $\sup_{r}|c(r)|\leq n^{k}$. Using Lemma
\ref{lem.epsilon2}, with $\lambda=n^{k}$, we have
\begin{align}
&  {{\mathbb{E}}}\left(  \sum_{r_{k-1}=1}^{n}\cdots\sum_{r_{1}=1}%
^{n}\left\Vert \sum_{r_{k}=1}^{n}c(r)\epsilon_{r_{k}}\right\Vert
^{2}\right)  =\sum_{r_{k-1}=1}^{n}\cdots\sum_{r_{1}=1}^{n}{{\mathbb{E}}%
}\left\Vert \sum_{r_{k}=1}^{n}c(r)\epsilon_{r_{k}}\right\Vert
^{2}\nonumber\\
&  =\sum_{r_{k-1}=1}^{n}\cdots\sum_{r_{1}=1}^{n}O\left(  \frac{n^{2k}}{n^{3}%
}\right)  =n^{k-1}O\left(  \frac{n^{2k}}{n^{3}}\right)
=n^{3(k-1)}O\left(
\frac{1}{n}\right)  .\label{eq.equationTb}%
\end{align}
Hence using Equations (\ref{eq.equationTa}) and
(\ref{eq.equationTb}), we have from Equation
(\ref{eq.EquationSb}),
\[
T_{n}\leq C\frac{1}{n^{3(k-1)}}n^{3(k-1)}\frac{1}{n}=O\left(  \frac{1}%
{n}\right)  .
\]
Thus putting all together, we have
\begin{equation}
{{\mathbb{E}}}\parallel\Xi^{n}\parallel\leq
d\cdot\sqrt{T_{n}}=O\left( \frac{1}{\sqrt{n}}\right)  .\nonumber
\end{equation}

\end{proof}

\begin{defn}
Let
\[
h(s,t)=(\min{s}{t})\Gamma(u(t))
\]
and for $k\in\mathbb{N},$ let%
\begin{equation}
\gamma_{k}=\Big(\frac{1}{12}\Big)^{k}\int_{s_{k}=0}^{1}\cdots\int_{s_{1}%
=0}^{1}{\mathrm{{tr}}}\ \Big[\prod_{l=1}^{k}\Big]h(s_{l},s_{l+1}%
)\ ds_{1}\ldots ds_{k},\nonumber
\end{equation}
where again by convention, $s_{k+1}:=s_{1}$.
\end{defn}

\begin{prop}
\label{p.riem}Continuing the notation in the above definition,
\begin{equation}
\gamma_{k}^{n}\longrightarrow\gamma_{k}\nonumber
\end{equation}
$\mu$-a.s. as $n\rightarrow\infty$.
\end{prop}

\begin{proof}
To begin with, let us write
\begin{align*}
s  & =(s_{1},s_{2},\ldots s_{k})\in\lbrack0,1)^{k},\ \\
r  & =(r_{1},r_{2},\ldots,r_{k})\in W_{n}^{k}:=\{1,2,\ldots
n\}^{k},\ \\
ds &= ds_1 \ldots ds_k%
\end{align*}
and
\begin{equation}
H(s)=\Big[\prod_{l=1}^{k}\Big]h(s_{l},s_{l+1})\nonumber
\end{equation}
where $s_{k+1} = s_1$. Denote $r_{l}^{-}:=r_{l}-1$ and
\begin{equation}
V_{n}(r)=V_{n}(r_{1},r_{2},\ldots,r_{k}):=J_{r_{1}}\times
J_{r_{2}}\cdots
J_{r_{k}}=\left(\frac{r_{1}^{-}}{n},\frac{r_{1}}{n}\right]\times\left(\frac
{r_{2}^{-}}{n},\frac{r_{2}}{n}\right]\cdots\left(\frac{r_{k}^{-}}{n},\frac{r_{k}%
}{n}\right].\nonumber
\end{equation}
With this new notation and $r_{k+1} = r_1$,
\begin{align}
  \gamma_{k}^{n}-\gamma_{k}
&  =\sum_{r\in W_{n}^{k}}{\mathrm{{tr}}}\ \Big[\prod_{l=1}%
^{k}\Big]\Lambda_{r_{i},r_{i+1}}^{n}-\left(\frac{1}{12}\right)^{k}%
\int_{s\in\lbrack0,1)^{k}}{\mathrm{{tr}}}\ H(s)\ ds\nonumber\\
&  =\sum_{r\in W_{n}^{k}}\left({\mathrm{{tr}}}\ \Big[\prod_{l=1}%
^{k}\Big]\Lambda_{r_{i},r_{i+1}}^{n}-\left(\frac{1}{12}\right)^{k}%
\int_{s\in V_{n}(r)}{\mathrm{{tr}}}\ H(s)\ ds%
\right)\nonumber\\
&  =\sum_{r\in W_{n}^{k}}\int_{s\in V_{n}(r)}\left(n^{k}{\mathrm{{tr}}%
}\
\Big[\prod_{l=1}^{k}\Big]\Lambda_{r_{i},r_{i+1}}^{n}-\left(\frac
{1}{12}\right)^{k}{\mathrm{{tr}}}\ H(s)\right)\ ds\nonumber\\
&  =\int_{s\in\lbrack0,1)^{k}}\sum_{r\in W_{n}^{k}}1_{V_{n}(r)}%
(s)\left({\mathrm{{tr}}}\ \Big[\prod_{l=1}^{k}\Big]n\Lambda_{r_{i},r_{i+1}%
}^{n}-\left(\frac{1}{12}\right)^{k}{\mathrm{{tr}}}\ H(s)\right)\ ds%
.\nonumber
\end{align}
Now $u_n \rightarrow u$ in the sup norm $\mu$-a.s. and
$\Gamma(\cdot)$ is continuous, thus
\begin{align}
\sum_{l, m = 1}^n1_{J_{l}\times J_{m}}(\tau,
t)\left(\min{\frac{l}{n}}{\frac{m-1}{n}}\right)\Gamma
\left(u_{n}\left(\frac{m-1}{n}\right)\right) - h(\tau,t) &
\longrightarrow 0\ \text{and} \nonumber \\
\sum_{l, m = 1}^n1_{J_{l}\times J_{m}}(\tau,
t)\left(\min{\frac{l}{n}}{\frac{m}{n}}\right)\Gamma
\left(u_{n}\left(\frac{m-1}{n}\right)\right) - h(\tau,t) &
\longrightarrow 0 \nonumber
\end{align}
$\mu$-a.s. as $n \rightarrow \infty$. For $(\tau, t) \in [0,1)^2$,
by using Definition \ref{defn.lambda&phi},
\begin{align}
&\sum_{l, m = 1}^n  1_{J_{l}\times J_{m}}(\tau,
t)n\Lambda_{lm}^{n}
-\frac{1}{12}h(\tau, t) \nonumber \\
&=\sum_{l, m = 1}^n  1_{J_{l}\times J_{m}}(\tau, t) \left\{
\frac{7}{360}\left(\min{\frac{l}{n}}{\frac{m-1}{n}}\right)\Gamma
\left(u_{n}\left(\frac{m-1}{n}\right)\right)+\frac{1}{45}\left(\min{\frac{l}{n}}{\frac{m}{n}}\right)\Gamma
\left(u_{n}\left(\frac{m}{n}\right)\right)1_{\{m<n\}}\atop
+\frac{1}{45}\left(\min{\frac{l}{n}}{\frac{m}{n}}\right)\Gamma
\left(u_{n}\left(\frac{m-1}{n}\right)\right)+\frac{7}{360}\left(\min{\frac{l}{n}}{\frac{m+1}{n}}\right)\Gamma
\left(u_{n}\left(\frac{m}{n}\right)\right)1_{\{m<n\}} \right\} \nonumber \\
 &  \hspace{1cm}- \frac{1}{12}h(\tau, t) \nonumber \\
& \longrightarrow 0 \nonumber
\end{align}
$\mu$-a.s. as $n \rightarrow \infty$. And since taking trace and
products are continuous operations, for $s \in [0,1)^k$, we hence
have
\begin{align}
g_{n}(s)&:=\sum_{r\in W_{n}^{k}}1_{V_{n}(r)}(s)\left({\mathrm{{tr}}%
}\
\Big[\prod_{l=1}^{k}\Big]n\Lambda_{r_{i},r_{i+1}}^{n}-\Big(\frac
{1}{12}\Big)^{k}{\mathrm{{tr}}}\ H(s) \right) \nonumber \\
&= \sum_{r\in W_{n}^{k}}1_{V_{n}(r)}(s)\left({\mathrm{{tr}}%
}\
\Big[\prod_{l=1}^{k}\Big]n\Lambda_{r_{i},r_{i+1}}^{n}-{\mathrm{{tr}}}\
\Big[\prod_{l=1}^{k}\Big]\frac {1}{12}h(s_{l},s_{l+1})
\right)\longrightarrow 0 \nonumber
\end{align}
$\mu$-a.s. as $n \rightarrow \infty$. Thus we can now apply
dominated convergence theorem and this gives us
\begin{align}
\gamma_{k}^{n}-\gamma_{k} &  =\int_{s \in [0,1)^k} g_{n}(s)\ ds
\longrightarrow0\nonumber
\end{align}
as $n\rightarrow\infty$.
\end{proof}

\begin{rem}
\label{rem.fk} If $S_{v}$ is the sectional curvature of the
manifold and $\sup_{v \in O(M)} \parallel \nolinebreak
S_{v}\parallel < \frac{3}{17d}$, then for any orthonormal frame
$\{e_{i}\}_{i=1}^{d} \subseteq T_oM$,
\begin{equation}
\sup_{v \in O(M)} \parallel\Omega_{v}(e_{i}, \cdot)e_{j}
\parallel< 2/d,\nonumber
\end{equation}
and hence
\begin{equation}
\sup_{v \in O(M)} \parallel\Gamma(v) \parallel< 3d^{2}\cdot4/d^{2}
= 12.\nonumber
\end{equation}
See Definition \ref{def.K} and Proposition \ref{p.curve}. Let
$\kappa= \sup_{v \in O(M)} \parallel\Gamma(v) \parallel/ 12$. Then
$\kappa< 1$ and hence
\begin{equation}
|\gamma_{k}^{n}| \leq\sum_{r_{k}=1}^{n} \cdots\sum_{r_{1}=1}^{n}
d\Big(\frac{1}{12n}\Big)^{k} \sup_{v \in O(M)} \parallel\Gamma(v)
\parallel^{k} = d\kappa^{k}.\nonumber
\end{equation}

\end{rem}

A standard result in probability states that $x_n \rightarrow x$
in $\mu$-measure iff for every subsequence of $\{x_n,\ n \geq 1\}$
has itself a subsequence converging $\mu$-a.s. to $x$. Hence for
any continuous function $f$, $f(x_n)$ converges to $f(x)$ in
$\mu$-measure. In the next proof, we will be using this result
without further comment.

\begin{thm}
\label{thm.fk} If $\ \sup_{v \in O(M)}\parallel S_{v} \parallel <
\frac{3}{17d}$, then
\begin{equation}
\det[ \left({\mathcal{I}^{n}} + {\mathcal{U}^{n}}\right)\circ b_n]
\longrightarrow e^{\gamma }\ \nonumber
\end{equation}
in $\mu$-measure as $n \rightarrow\infty$, where $\gamma$ is
defined as
\begin{equation}
\gamma := \sum_{k=1}^{\infty}(-1)^{k+1}\frac{1}{k}
\gamma_{k}.\nonumber
\end{equation}

\end{thm}

\begin{proof}
To ease the notation, in this proof only, we will write
$\mathcal{U}^n = \mathcal{U}^n \circ b_n$ and $\mathcal{I}^n =
\mathcal{I}^n \circ b_n$. By Remark \ref{rem.fk},
$|\gamma_{k}^{n}| \leq d\kappa^{k}$, $\kappa= \sup_{v \in O(M)}
\parallel \Gamma(v)
\parallel/ 12 < 1$. Observe that ${\mathcal{U}^{n}}$ satisfies the
hypothesis in Lemma \ref{lem.posdetformula}. Thus we can apply the
formula in Equation (\ref{eq.detformula2}). Hence
\begin{equation}
\det( {\mathcal{I}^{n}} + {\mathcal{U}^{n}}) = \exp( \Psi_{r} +
R_{r+1} ),\nonumber
\end{equation}
where $\Psi_{r} = \sum_{k=1}^{r} (-1)^{k+1}\frac{1}{k}
{\mathrm{{Tr}}} \ ([ {\mathcal{U}^{n}}]^{k})$ and $|R_{r+1}|
\leq\frac{1}{r+1} {\mathrm{{Tr}}} \ ([ {\mathcal{U}^{n}}]^{r+1})$.
Therefore
\begin{align}
 & \left|  \det( {\mathcal{I}^{n}} + {\mathcal{U}^{n}}) -
\exp\left( \sum_{k=1}^{r} (-1)^{k+1}\frac{1}{k}
\gamma_{k}^{n}\right)  \right|&
\nonumber\\
&  = \exp\left(  \sum_{k=1}^{r} \frac{(-1)^{k+1}}{k}
\gamma_{k}^{n} \right) \left|\exp\left (\sum_{k=1}^{r}
\frac{(-1)^{k+1}}{k}\left( {\mathrm{{Tr}}} \ ([
{\mathcal{U}^{n}}]^{k}) - \gamma_{k}^{n} \right) + R_{r+1}\right)
- 1
\right|\nonumber\\
&  \leq\exp( d\kappa/(1-\kappa) ) \left|\exp\left (\sum_{k=1}^{r}
(-1)^{k+1}\frac{1}{k}\left( {\mathrm{{Tr}}} \ ([
{\mathcal{U}^{n}}]^{k}) - \gamma_{k}^{n} \right) + R_{r+1}\right)
- 1 \right|.\nonumber
\end{align}
The last inequality follows from
\begin{equation}
\exp\left(  \sum_{k=1}^{r} (-1)^{k+1}\frac{1}{k} \gamma_{k}^{n}
\right) \leq\exp\left(  \sum_{k=1}^{r} \frac{1}{k} d\kappa^{k}
\right)  < \exp (d\kappa/(1-\kappa)).\nonumber
\end{equation}
Now by Theorem \ref{thm.trf_k}, we have
\begin{equation}
{\mathrm{{Tr}}} \ ([ {\mathcal{U}^{n}}]^{k}) - \gamma_{k}^{n}
\rightarrow 0\ \nonumber
\end{equation}
in $\mu$-measure as $n \rightarrow\infty$. Together with
$\parallel\gamma _{k}^{n} \parallel\leq d\kappa^{k}$, we will have
\begin{equation}
\limsup_{n \rightarrow\infty} |R_{r+1}| \leq
\frac{1}{r+1}\limsup_{n \rightarrow \infty} \mathrm{Tr}\ (\left[
\mathcal{U}^n \right]^{r+1}) \leq \frac{1}{(r+1)
}d\kappa^{r+1}\leq d\kappa^{r+1},\nonumber
\end{equation}
in $\mu$-measure. Therefore
\begin{align}
&  \limsup_{n \rightarrow\infty} \left|  \det[ {\mathcal{I}^{n}} +
{\mathcal{U}^{n}}] - \exp\left(  \sum_{k=1}^{r}
(-1)^{k+1}\frac{1}{k} \gamma^{n}_{k}\right)  \right|
\leq\exp(d\kappa/(1-\kappa)) \left|\exp\left( d\kappa^{r+1}\right)
- 1
\right| \label{eq.limsup1}%
\end{align}
in $\mu$-measure. From Proposition \ref{p.riem}, we know that
$\gamma_{k}^{n} \rightarrow\gamma_{k}$ and hence $\gamma_{k} \leq
d\kappa^{k}$ , which implies $p_{r} = \sum_{k=1}^{r}
\frac{(-1)^{k+1}}{k} \gamma_{k}$ converges to $\gamma$ and
$|p_{r}| \leq\frac{d\kappa}{1-\kappa}$. Thus
\begin{align}
&  \lim_{n \rightarrow\infty}\left|  e^{\gamma} - \exp\left(
\sum_{k=1}^{r} \frac{(-1)^{k+1}}{k} \gamma_{k}^{n}\right)  \right|
\leq\left|  e^{\gamma} - e^{p_{r}} \right|  + \lim_{n
\rightarrow\infty} \left|  e^{p_{r}} - \exp\left(  \sum_{k=1}^{r}
\frac{(-1)^{k+1}}{k} \gamma_{k}^{n}\right)  \right|
\nonumber\\
&  \leq e^{p_{r}} \left|  \exp\left( \sum_{k=r+1}^{\infty}
(-1)^{k+1}\frac {1}{k} \gamma_{k}\right) - 1\right|  \leq
e^{p_{r}}\left| \exp\left( \sum
_{k=r+1}^{\infty} \frac{1}{k} d\kappa^{k} \right) - 1\right| \nonumber\\
&  \leq e^{d\kappa/(1-\kappa)}\left|  \exp\left(
\frac{d\kappa^{r+1}}{1-\kappa}
\right) - 1\right|  \label{eq.limsup2}%
\end{align}
$\mu$-a.s.. Therefore using Equations (\ref{eq.limsup1}) and
(\ref{eq.limsup2}),
\begin{align}
  \limsup_{n \rightarrow\infty} &\left|  \det[ {\mathcal{I}^{n}} +
{\mathcal{U}^{n}}] - e^{\gamma} \right|  \nonumber\\
  \leq& \limsup_{n \rightarrow\infty}\left\{  \left|  e^{\gamma} -
\exp\left( \sum_{k=1}^{r} \frac{(-1)^{k+1}}{k}
\gamma_{k}^{n}\right) \right| \atop  + \left| \det[
{\mathcal{I}^{n}} + {\mathcal{U}^{n}}] - \exp\left(
\sum_{k=1}^{r}
\frac{(-1)^{k+1}}{k} \gamma_{k}^{n}\right)  \right|  \right\} \nonumber\\
  \leq& \exp(d\kappa/(1-\kappa))\left|  \exp\left(  \frac{d\kappa^{r+1}%
}{1-\kappa} \right)  - 1\right|  + \exp(d\kappa/(1-\kappa))
\left|\exp( d\kappa^{r+1}) - 1 \right|\nonumber
\end{align}
in $\mu$-measure. Since the inequality holds for every $r$ and as
$r \rightarrow\infty$,
\begin{equation}
\left|\exp( a\kappa^{r+1} ) - 1 \right| \rightarrow0,\nonumber
\end{equation}
for any constant $a$, we thus have
\begin{equation}
\limsup_{n \rightarrow\infty} \left|\ \det[ \left({\mathcal{I}^{n}} + {\mathcal{U}%
^{n}}\right)\circ b_n] - e^{\gamma}\ \right| = 0\ \nonumber
\end{equation}
in $\mu$-measure.
\end{proof}

\begin{defn}\label{defn.Ku}
Define an integral operator $K_{u} : L^{2}([0,1]
\rightarrow\mathbb{R}^{d}) \longrightarrow L^{2}([0,1]
\rightarrow\mathbb{R}^{d})$ by
\begin{align}
(K_{u} f)(s) &:= \int_{0}^{1} h(s,t) f(t)\ dt
\nonumber \\
&= \int_{0}^{1} (\min{s}{t})\ \Gamma(u(t)) f(t)\ dt,\nonumber
\end{align}
where $\Gamma$ was defined in Definition \ref{defn.Gammau}.

\end{defn}

\begin{prop}
$K_{u}$ is trace class.
\end{prop}

\begin{proof}
Note that $L^{2}([0,1] \rightarrow\mathbb{R}^{d})$ is a separable
Hilbert space. Let $(Af)(s) = \int_{0}^{1} (\min{s}{t})f(t)\ dt$
and $(Bf)(s) = \Gamma(u(s))f(s)$. Thus $K_{u} = AB$. By
Proposition \ref{prop.trace1}, $A$ is trace class. Since $B$ is
bounded, by Proposition \ref{prop.trace2}, ${\mathrm{{tr}}} \ |AB|
< \infty$ and hence $K_{u}$ is trace class.
\end{proof}

\begin{prop}
\label{prop.traceK_u}
\begin{equation}
e^{\gamma} = \det\left(I + \frac{1}{12}K_{u} \right).\nonumber
\end{equation}

\end{prop}

\begin{proof}
We will use Equation (\ref{eq.detformula}) to prove the statement.
Note that by Remark \ref{rem.fk}, for any $f \in L^{2}([0,1]
\rightarrow\mathbb{R}^{d})$,
\begin{eqnarray}
\Big \Vert \int_0^1 \frac{1}{12}(K_uf)^T(s)f(s)\ ds \Big \Vert &
\leq & \left( \int_0^1 \frac{1}{12}\parallel (K_uf)(s)
\parallel^2\ ds \right)^{\frac{1}{2}} \left( \int_0^1 \frac{1}{12}
\parallel f(s)
\parallel^2\ ds \right)^{\frac{1}{2}} \nonumber \\
& \leq & \frac{1}{12}\parallel \Gamma_u \parallel \int_0^1
\parallel f(s) \parallel^2\ ds < \int_0^1
\parallel f(s) \parallel^2\ ds \nonumber
\end{eqnarray}
and hence $\parallel \frac{1}{12}K_u \parallel < 1$ . Thus
\begin{equation}
\det\left(I + \frac{1}{12}K_{u}\right) = \exp
\left[{\sum_{k=1}^{\infty}(-1)^{k+1}\frac{1}{k} {\mathrm{{tr}}} \
\left(\frac{1}{12}K_{u}\right)^{k}} \right].\nonumber
\end{equation}
Therefore it suffices to show that for $k \in \mathbb{N}$,
\begin{align}
{\mathrm{{tr}}} \ (K_{u}^{k} )  &  = {\mathrm{{tr}}} \ \int_{s_k =
0}^{1} \ldots \int_{s_1 = 0}^{1} \Big[\prod_{l=1}^{k}\Big]
h(s_{l}, s_{l+1})\ ds_{1} \ldots ds_{k}  = \gamma_{k}(u),\nonumber
\end{align}
where by convention, $s_{k+1} = s_1$. Let
\begin{equation}
(K_{u}^{k} f)(s) = \int_{0}^{1} p_k(s,s_{1})f(s_{1})\
ds_{1},\nonumber
\end{equation}
where
\begin{equation}
p_1(s, s_1) =  h(s, s_1) \nonumber
\end{equation}
and for $k \geq 2$, \begin{equation} p_k(s, s_1) = \int_{s_k =
0}^{1} \cdots\int_{s_2 = 0}^{1} h(s,
s_{2})\Big[\prod_{l=2}^{k}\Big] h(s_{l}, s_{l+1})\ ds_{2} \cdots
ds_{k}. \nonumber
\end{equation}
Let $\Omega\subseteq L^{2}([0,1])$ be an orthonormal basis and
hence
\begin{equation}
\{\varpi e_{i}\ |\ \varpi\in\Omega,\ i = 1,2, \ldots d \}\nonumber
\end{equation}
is an orthonormal basis for $L^{2}([0,1], \mathbb{R}^{d})$. Let
$(, )$ denote the inner product on \\$L^{2}([0,1],
\mathbb{R}^{d})$ and $\langle, \rangle$ be the inner product on
$\mathbb{R}^{d}$. Then
\begin{align}
{\mathrm{{tr}}} \ (K_{u}^{k})  &  = \sum_{i,\varpi \in \Omega}
(K_{u}^{k} \varpi e_{i},\ \varpi
e_{i} )  = \sum_{i,\varpi \in \Omega} \int_{0}^{1} \int_{0}^{1} \langle p_k(s,t) \varpi(t)e_{i}%
,\ \varpi(s)e_{i} \rangle\ dt\ ds\nonumber\\
&  = \sum_{\varpi \in \Omega} \int_{0}^{1} \int_{0}^{1}
\varpi(t)\varpi(s) {\mathrm{{tr}}} \ p_k(s,t) \ dt\ ds  \nonumber
\\
& = \int_{0}^{1}  \sum_{\varpi \in \Omega}\left(\int_{0}^{1}
\varpi(t) {\mathrm{{tr}}}
\ p_k(s,t) \ dt \right) \varpi(s)\ ds \nonumber \\
& = \int_{0}^{1} {\mathrm{{tr}}} \ p_k(s,s)\ ds  = {\mathrm{{tr}}}
\ \int_{s_k=0}^{1} \ldots\int_{s_1=0}^{1}
\Big[\prod_{l=1}^{k}\Big] h(s_{l}, s_{l+1})\ ds_{1} \ldots
ds_{k}.\nonumber
\end{align}
Since $K_{u}$ is trace class, all the interchanging between the
integrals and the sums are valid.
\end{proof}

\subsection{Convergence of $\widetilde{ \det(\mathcal{I}^n + \mathcal{X}^n)\circ b_n }$} \label{s.X}

\begin{lem}
\label{lem.1limit}
\begin{equation}
\widetilde { \det ( {\mathcal{I}^{n}}+\mathcal{X}^n )\circ b_n } =
1_{\{n \geq n_0\}}\det ( {\mathcal{I}^{n}}+\mathcal{X}^n )\circ
b_n \longrightarrow1 \nonumber
\end{equation}
$\mu-a.s.$ as $n\rightarrow\infty$.
\end{lem}

\begin{proof}
By Lemma \ref{lem.walpha}, it suffices to consider $
W_\alpha(\mathbb{R}^d)$ with $\frac{7}{15} < \alpha <
\frac{1}{2}$. Recall
\begin{equation}
{\mathcal{X}^{n}}=( {\mathcal{I}^{n}}+{\mathcal{U}^{n}})^{-1}%
({\mathcal{S}^{n}})^{T}{\mathcal{E}^{n}}{\mathcal{S}^{n}}.\nonumber
\end{equation}
Now observe that ${\mathcal{I}^{n}}+{\mathcal{U}^{n}}$ is a
positive definite matrix with eigenvalues greater than or equal to
1. Thus since the norm of a symmetric matrix is equal to the
maximum eigenvalue and the eigenvalues of $(
{\mathcal{I}^{n}}+{\mathcal{U}^{n}})^{-1}$ is the reciprocal of
the eigenvalues of ${\mathcal{I}^{n}} +{\mathcal{U}^{n}}$, we have
\begin{equation}
\left \Vert ( {\mathcal{I}^{n}}+{\mathcal{U}^{n}})^{-1}\right
\Vert \leq 1.\nonumber
\end{equation}
Note that $\parallel \mathcal{E}^n \circ b_n
\parallel=O\left(\bigvee_{i=1,\ldots,n}\parallel\Delta
_{i}b\parallel^{5}\right)$ for $n \geq n_0(\cdot)$. Thus by Lemma
\ref{lem.cB}, on $W_\alpha(\mathbb{R}^d)$ and for $n \geq
n_0(\cdot)$, we have
\begin{equation}
\parallel{\mathcal{X}^{n}}\circ b_n\parallel=O \left(n^{2}\bigvee_{i=1,\ldots
,n}\parallel\Delta_{i}b\parallel^{5}\right).\nonumber
\end{equation}
Let $\parallel{\mathcal{X}^{n}}\circ b_n\parallel\leq
C[n^{2}\bigvee_{i=1,\ldots ,n}\parallel\Delta_{i}b\parallel^{5}]$,
$C$ is a fixed constant independent of $\omega$ and $n$.
Therefore \beq \parallel{\mathcal{X}^{n}%
}\circ b_n\parallel \leq
Cn^2\bigvee_{i=1,\ldots,n}\parallel\Delta_{i}b\parallel^5 < C n^{2
- 5\alpha}<Cn^{-\frac{1}{3}} \nonumber \eeq and hence
$\parallel{\mathcal{X}^{n}}\parallel<1$. Thus using Equation
(\ref{eq.detformula}), we have
\begin{align}
\det\left[  {\mathcal{I}^{n}}+{\mathcal{X}^{n}}\right]
&=\exp\left[\sum
_{k=1}^{\infty}\frac{1}{k}(-1)^{k+1}{\mathrm{{Tr}}}\ \left[({\mathcal{X}^{n}}%
)^{k}\right]\right] \nonumber \\
&=\exp\left[{\mathrm{{Tr}}%
}\ {\mathcal{X}^{n}}-\frac{1}{2}{\mathrm{{Tr}}}\
\left[({\mathcal{X}^{n}})^{2}\right]+ R_3\left(
{\mathcal{X}^{n}}\right) \right],\nonumber
\end{align}
where $R_3({\mathcal{X}^{n}})=\sum_{k=3}^{\infty}\frac{1}{k}(-1)^{k+1}%
{\mathrm{{Tr}}}\ \left[({\mathcal{X}^{n}})^{k}\right]$. Note that
$\mathcal{X}^{n}$ is a $nd \times nd$ matrix and hence, \beq
|R_3({\mathcal{X}^{n}\circ b_n})|\leq
\frac{ndC^{3}n^{3(2 - 5\alpha)}%
}{1-Cn^{2 - 5\alpha}}=\frac{dC^{3}n^{3(2 - 5\alpha)+1}%
}{1-Cn^{2 - 5\alpha}}. \nonumber \eeq By the choice of $\alpha$,
$\tau := 3(2 - 5\alpha)+1 < 0$ and hence \beq |
R_3({\mathcal{X}^{n}\circ b_n})| = O\left(n^{\tau}\right)
\longrightarrow 0 \nonumber \eeq as $n \rightarrow \infty$. Let
$P:=( {\mathcal{I}^{n}}+{\mathcal{U}^{n}} \circ b_n)^{-1},$
$S=\mathcal{S}^{n}$ and $E=\mathcal{E}^{n} \circ b_n$ to simplify
notation. Then, using Proposition \ref{prop.linalgebra2} twice,
\begin{align*}
\left\vert \mathrm{Tr}\ {\mathcal{X}^{n} \circ b_n}\right\vert  &
=\left\vert \operatorname{Tr}\left(  PS^{T}ES\right)  \right\vert
=\left\vert
\operatorname{Tr}\left(  SPS^{T}E\right)  \right\vert \\
&  \leq\left\Vert E\right\Vert \operatorname{Tr}\left(
SPS^{T}\right) =\left\Vert E\right\Vert \operatorname{Tr}\left(
S^{T}SP\right) \leq\left\Vert E\right\Vert \left\Vert P\right\Vert
\operatorname{Tr}\left( S^{T}S\right)
\end{align*}
where $\operatorname{Tr}\left(  S^{T}S\right)  =O\left(
n^{2}\right) $. This is because from Lemma \ref{lem.cB}, we have
\beq \mathrm{Tr}\ \left( S^{T}S\right) = \mathrm{Tr}\ \cB =
\sum_{m=1}^n m\mathrm{tr}\ I = \frac{dn(n+1)}{2}. \nonumber \eeq
Thus
\begin{align}
| {\mathrm{{Tr}}} \ \mathcal{X}^{n} \circ b_n|  &  \leq O(n^{2})
\parallel {E}
\parallel
= O\left( n^2 (n^{-\alpha})^5 \right) \nonumber \\
&= O\left( n^2 (n^{-7/15})^5 \right) = O( n^{-1/3} )
\nonumber \\
& \longrightarrow0 \nonumber
\end{align}
as $n \rightarrow\infty$. Similarly, using Lemma
\ref{prop.linalgebra2} twice and the fact that $P$ is positive
definite so that $SPS^{T}$ is positive definite, we have
\begin{align*}
\left\vert {\mathrm{{Tr}}}\left[ ({\mathcal{X}^{n}}\circ
b_n)^{2}\right] \right\vert &  =\left\vert \operatorname{Tr}\left(
PS^{T}ESPS^{T}ES\right)  \right\vert
=\left\vert \operatorname{Tr}\left(  SPS^{T}ESPS^{T}E\right)  \right\vert \\
&  \leq\operatorname{Tr}\left(  SPS^{T}\right)  \left\Vert ESPS^{T}%
E\right\Vert =\operatorname{Tr}\left(  S^{T}SP\right)  \left\Vert
ESPS^{T}E\right\Vert \\
&  \leq\operatorname{Tr}\left(  S^{T}S\right)  \cdot\left\Vert
P\right\Vert
\left\Vert ESPS^{T}E\right\Vert \\
&  \leq O(n^{2})\left\Vert ESPS^{T}E\right\Vert  \leq O(n^{2})\left\Vert E\right\Vert ^{2}\left\Vert S\right\Vert ^{2}\\
&  = O(n^{4})\left\Vert E\right\Vert ^{2}  =
 O\left( n^4 (n^{-\alpha})^{10} \right) \\
&= O\left( n^4 (n^{-7/15})^{10} \right) = O\left( n^4 n^{-14/3} \right) \\
& = O\left( n^{\frac{12-14}{3}} \right) = O(n^{-2/3}) \\
& \longrightarrow0
\end{align*}
as $n \rightarrow\infty$. Therefore,
\begin{align}
\widetilde{ \det ( {\mathcal{I}^{n}}+{\mathcal{X}^{n}})\circ b_n }
& =  1_{\{n \geq n_0\}}\det (
{\mathcal{I}^{n}}+{\mathcal{X}^{n}})\circ b_n
\nonumber \\
& = \exp
\left[{\mathrm{{Tr}} }\ {\mathcal{X}^{n}}-\frac{1}{2}{\mathrm{{Tr}}%
}\ \left[({\mathcal{X}^{n}})^{2}\right] + R_3\left(  {\mathcal{X}^{n}%
}\right)  \right] \circ b_n \nonumber\\
&  \longrightarrow 1 \nonumber
\end{align}
$\mu$-a.s. as $n\rightarrow\infty$.
\end{proof}

As a random variable on $(W(\mathbb{R}^{d}),\ \mu)$, the next
theorem shows that $\rho_{n} \circ\phi\circ b_n$ converges in
$\mu$-measure.

\begin{thm}
\label{thm.rho1} If $\ \sup_{v \in O(M)}\parallel S_v \parallel <
\frac{3}{17d}$, then
\begin{equation}
\rho_{n} \circ\phi\circ b_n \longrightarrow
e^{-\frac{1}{6}\int_{0}^{1} Scal(\tilde{\phi}(s))\ ds
}\sqrt{\det\left(I + \frac{1}{12}K_{u} \right)} \ .\nonumber
\end{equation}
in $\mu$-measure as $n \rightarrow\infty$.
\end{thm}

\begin{proof}
It suffices to consider on $W_\alpha(\mathbb{R}^d)$, with
$\frac{7}{15} < \alpha < \frac{1}{2}$. From Theorem \ref{t.ket},
\begin{equation}
\det [\langle (\cF)^T \cF \rangle] = [\det (\mathcal{V}^n)]^2
\det( {\mathcal{I}^{n}}+{\mathcal{U}^{n}})%
\det\Big[ {\mathcal{I}^{n}}+\Big( {\mathcal{I}^{n}}+{\mathcal{U}^{n}%
}\Big)^{-1}({\mathcal{S}^{n}})^{T}{\mathcal{E}^{n}}{\mathcal{S}^{n}%
}\Big].\nonumber
\end{equation}
Using Lemma \ref{lem.prodscal}, Theorem \ref{thm.fk}, Proposition
\ref{prop.traceK_u} and Lemma \ref{lem.1limit}, we have,
\begin{align}
1_{\{n \geq n_o\}} \det [\langle (\cF)^T \cF \rangle \circ b_n] &
\longrightarrow e^{\frac{1}{3}\int_{0}^{1} Scal(\tilde{\phi
}(s))\ ds + \gamma}\nonumber\\
&  = e^{\frac{1}{3}\int_{0}^{1} Scal(\tilde{\phi}(s))\
ds}\det\left(I + \frac{1}{12}K_{u} \right)\nonumber
\end{align}
in $\mu$-measure as $n \rightarrow\infty$. But by Theorem
\ref{thm.rho}, for $n \geq n_0(\cdot)$,
\begin{equation}
\left(\rho_{n} \circ\phi \circ b_n \right)^{2} = \det [\langle
(\cF)^T \cF \rangle]\circ b_n  .\nonumber
\end{equation}
Thus taking the square root completes the proof.
\end{proof}

\section{$L^{1}$ Convergence of $\{\rho_{
n} \circ \phi \circ b_n\}_{n=1}^\infty$}\label{s.L1}

\begin{defn} \label{d.mup}
Let $\mu_{G_{ {\mathcal{P}}}^{1}}$ be defined as in Theorem
\ref{t.m_G} on $H_{\mathcal{P}}(\mathbb{R}^d)$ by the density
\begin{equation}
 \frac{1}{Z_{ {\mathcal{P}}}^{1}}%
e^{-\frac{1}{2}E_{\mathbb{R}^{d}}}Vol_{ G_{
{\mathcal{P}}}^{1}}\nonumber
\end{equation}
where $Z_{ {\mathcal{P}}}^{1} := (2\pi)^{\frac{dn}{2}}$ and $E_{\mathbb{R}%
^{d}}(\omega) := \int_{0}^{1} \parallel
\omega^{\prime}(s)\parallel^{2}\ ds$.
\end{defn}

The following 2 results can be found in \cite{MR1698956} and hence
we will omit the proofs.

\begin{thm}
\label{thm.pullbackmu} Let $\mu_{G_{ {\mathcal{P}}}^{1}}$ be
defined as in Definition \ref{d.mup}. Write $\phi_{
{\mathcal{P}}}= \phi|_{H_{ {\mathcal{P}}}(\mathbb{R}^d)}$. Then
$\mu_{G_{ {\mathcal{P}}}^{1}}$ is the pullback of $\nu_{G_{
{\mathcal{P}}}^{1}}$ by $\phi_{ {\mathcal{P}}}$. i.e.
\begin{equation}
\mu_{G_{ {\mathcal{P}}}^{1}} = (\phi_{ {\mathcal{P}}})^{*}
\nu_{G_{ {\mathcal{P}}}^{1}}.\nonumber
\end{equation}

\end{thm}

Let $\pi_{ {\mathcal{P}}}: W(\mathbb{R}^{d})
\rightarrow(\mathbb{R}^{d})^{n}$ be given by
\begin{equation}
\pi_{ {\mathcal{P}}}(\omega) := \Big( \omega(s_{1}),\
\omega(s_{2}),\ \ldots, \omega(s_{n}) \Big).\nonumber
\end{equation}
Note that $\pi_{ {\mathcal{P}}}: H_{
{\mathcal{P}}}(\mathbb{R}^{d}) \rightarrow(\mathbb{R}^{d})^{n}$ is
a linear isomorphism of finite dimensional vector spaces. Let $i_{
{\mathcal{P}}}: (\mathbb{R}^{d})^{n} \rightarrow H_{
{\mathcal{P}}}(\mathbb{R}^{d})$ denote the inverse of ${\pi_{ {\mathcal{P}}}%
}|_{H_{ {\mathcal{P}}}(\mathbb{R}^{d})}$.

\begin{lem}
\label{lem.pullbackmu} Let $(y_1, \ldots y_n)$ denote the standard
cartesian coordinates on $(\mathbb{R}^{d})^{n}$ and $y_{0} := 0$.
Then
\begin{equation}
i_{ {\mathcal{P}}}^{*}\mu_{G_{ {\mathcal{P}}}^{1}} = \frac{1}{Z_{
{\mathcal{P}}}^{1}} \left( \prod_{i=1}^{n}
(\Delta_{i}s)^{-\frac{d}{2}} \exp\left(-\frac{1}{2 \Delta_{i}s}
\parallel y_{i} - y_{i-1} \parallel^{2}\right)
\right)\ dy_{1} dy_{2} \ldots dy_{n}. \label{eq.pullbackmu}%
\end{equation}
This equation can also be written as
\begin{equation}
i_{ {\mathcal{P}}}^{*}\mu_{G_{ {\mathcal{P}}}^{1}} = \left(
\prod_{i=1}^{n} p_{\Delta_{i}s}(y_{i-1},\ y_{i}) \right)\ dy_{1}
dy_{2} \ldots dy_{n}\nonumber
\end{equation}
where $p_{s}(x,\ y) :=
(\frac{1}{2\pi})^{-\frac{d}{2}}\exp(-\frac{\parallel x - y
\parallel}{2s})$ is the heat kernel on $\mathbb{R}^{d}$.
\end{lem}

We are now ready to prove the main result.

\begin{thm}
\label{thm.nu_converge} Let $M$ be a compact Riemannian manifold,
$f : W(M) \rightarrow\mathbb{R}$ be a bounded continuous function
and ${\mathcal{P}}= \{0 < \frac{1}{n} < \cdots < \frac{n}{n}=1\}$
is an equally spaced partition. Let $b$ be brownian motion, $u$
solves Equation (\ref{e.par.trans}), $\tilde{\phi} = \pi \circ u$
and $\widetilde{//}$ is stochastic parallel translation defined in
Definition \ref{d.s//trans}. Suppose  $ \sup_{v \in O(M)}\parallel
S_v \parallel <\frac{3}{17d}$, where $S_v$ is sectional curvature
and $\rho_{ {\mathcal{P}}}\circ\phi \circ b_\cP$ is uniformly
integrable. Then
\begin{align}
&  \int_{H_{ {\mathcal{P}}}(M)} f(\sigma)\ d\nu_{
{\mathcal{P}}}(\sigma) = \int_{H_{ {\mathcal{P}}}(M)} f(\sigma)\
\rho_{ {\mathcal{P}}}(\sigma)
\ d\nu_{G_{ {\mathcal{P}}}^{1}}(\sigma)\nonumber\\
&  {%
\genfrac{}{}{0pt}{}{ \longrightarrow}{{|\cP| \rightarrow 0} }%
} \int_{W(M)} f(\sigma) e^{-\frac{1}{6}\int_{0}^{1} Scal(\sigma(s))\ ds }%
\sqrt{\det\left(I + \frac{1}{12}K_{\widetilde{//}(\sigma)}
\right)}\ d\nu(\sigma)\nonumber
\end{align}
where $\nu$ is Wiener measure on $M$. $Scal$ is the scalar
curvature of the manifold and $K_{\widetilde{//}(\sigma)} :
L^{2}([0,1] \rightarrow\mathbb{R}^{d}) \longrightarrow L^{2}([0,1]
\rightarrow\mathbb{R}^{d})$ is given by
\begin{equation}
\left(K_{\widetilde{//}(\sigma)} v \right)(s) = \int_{0}^{1}
(\min{s}{t})\ \Gamma(\widetilde{//}_t(\sigma)) v(t)\ dt\nonumber
\end{equation}
where
\begin{eqnarray}
\Gamma(\widetilde{//}_t(\sigma)) & = & \sum_{i,j = 1}^{d} \left\{
{\Omega_{\widetilde{//}_t(\sigma)} (e_{i},
\Omega_{\widetilde{//}_t(\sigma)}(e_{i}, \cdot )e_{j} )e_{j} +
  \Omega_{\widetilde{//}_t(\sigma)} (e_{i},
\Omega_{\widetilde{//}_t(\sigma)}(e_{j}, \cdot)e_{i} )e_{j} \atop
+ \Omega_{\widetilde{//}_t(\sigma)} (e_{i},
\Omega_{\widetilde{//}_t(\sigma)}(e_{j}, \cdot )e_{j} )e_{i} }
\right\}\nonumber
\end{eqnarray}
for any orthonormal basis $\{e_{i}\}_{i=1}^{d} \subseteq T_oM$,
and \beq \Omega_{\widetilde{//}_t(\sigma)}(a, c) :=
\widetilde{//}_t^{-1}(\sigma) R\left(\widetilde{//}_t(\sigma)a,
\widetilde{//}_t(\sigma)c \right)\widetilde{//}_t(\sigma)
\nonumber \eeq for any vectors $a, c \in T_oM$, $R$ being the
curvature tensor.
\end{thm}

\begin{proof}
By Theorem \ref{thm.wongkai}, $f \circ\phi\circ b_{
{\mathcal{P}}}:= f_{ {\mathcal{P}}}$ converges to $f
\circ\tilde{\phi}$ $\mu-a.s.$ as $| {\mathcal{P}}| \rightarrow0$.
Thus we have
\begin{align}
 \int_{H_{ {\mathcal{P}}}(M)} f\ d\nu_{
{\mathcal{P}}}
& = \int_{H_{ {\mathcal{P}}}(M)} f\ \rho_{ {\mathcal{P}}}\ d\nu_{G_{ {\mathcal{P}}}^{1}} \nonumber\\
&  = \int_{H_{ {\mathcal{P}}}(\mathbb{R}^{d})} (f \rho_{
{\mathcal{P}}})\circ \phi\ d\mu_{G_{ {\mathcal{P}}}^{1}}%
\ (By\ Theorem\ \ref{thm.pullbackmu})\nonumber\\
&  = \int_{W(\mathbb{R}^{d})}   \ (f\rho_{
{\mathcal{P}}})\circ\phi\circ b_\cP \ d\mu
\ (By\ Lemma\ \ref{lem.pullbackmu})\nonumber\\
&  = \int_{W(\mathbb{R}^{d})}  \ f_{ {\mathcal{P}}}\cdot \rho_{ {\mathcal{P}}%
}\circ\phi\circ b_{ {\mathcal{P}}} \ d\mu.\nonumber
\end{align}
By the assumption on $\rho_{ {\mathcal{P}}}\circ\phi\circ b_{
{\mathcal{P}}}$, since $f$ is bounded and continuous, we have $f_{
{\mathcal{P}}}\cdot\ \rho_{ {\mathcal{P}}}\circ\phi\circ b_{
{\mathcal{P}}}$ is uniformly integrable. Therefore by Theorem
\ref{thm.rho1},
\begin{align}
&  \int_{H_{ {\mathcal{P}}}(M)} f(\sigma)\ d\nu_{
{\mathcal{P}}}(\sigma) =
\int_{W(\mathbb{R}^{d})} \left(f_{ {\mathcal{P}}} \cdot \rho_{ {\mathcal{P}}}%
\circ\phi\circ b_{ {\mathcal{P}}}\right)(\omega) \ d\mu(\omega)\nonumber\\
&  \longrightarrow\int_{W(\mathbb{R}^{d})} \left(f
\circ\tilde{\phi} \cdot e^{-\frac {1}{6}\int_{0}^{1}
Scal(\tilde{\phi}(s))\ ds }\sqrt{\det\left(I + \frac
{1}{12}K_{u} \right)}\right)(\omega)\ d\mu(\omega) \nonumber\\
&  = \int_{W(M)} f(\sigma)e^{-\frac{1}{6}\int_{0}^{1}
Scal(\sigma(s))\ ds }\sqrt{\det\left(I +
\frac{1}{12}K_{\widetilde{//}(\sigma)} \right)}\
d\nu(\sigma)\nonumber
\end{align}
as $|\cP| \rightarrow 0$, where $K_u$ was defined in Definition
\ref{defn.Ku}. Note that $\nu = \mu\tilde{\phi}^{-1}$ and $u =
\widetilde{//}(\tilde{\phi})$ from (\ref{rem.s//trans.1}) and
(\ref{rem.s//trans.2}) of Remark \ref{rem.s//trans} respectively.
\end{proof}

\begin{cor}
\label{cor.nu_converge} Let $M$ be a compact Riemannian manifold,
$f: W(M) \rightarrow\mathbb{R}$ be a bounded continuous function
and ${\mathcal{P}}$ is an equally spaced partition. Suppose that
$\forall v \in O(M)$, $0 \leq S_v <\frac{3}{17d}$. Then
\begin{align}
&  \int_{H_{ {\mathcal{P}}}(M)} f(\sigma)\ d\nu_{
{\mathcal{P}}}(\sigma)\nonumber \\%
&  {%
\genfrac{}{}{0pt}{}{ \longrightarrow}{{|\cP| \rightarrow 0} }%
} \int_{W(M)} f(\sigma) e^{-\frac{1}{6}\int_{0}^{1} Scal(\sigma(s))\ ds }%
\sqrt{\det\left(I + \frac{1}{12}K_{\widetilde{//}(\sigma)}
\right)}\ d\nu(\sigma).\nonumber
\end{align}

\end{cor}

\begin{proof}
By Theorem \ref{thm.uniformint}, we know that $\rho_{
{\mathcal{P}}}\circ \phi\circ b_{ {\mathcal{P}}}$ is uniformly
integrable under the assumptions on the sectional curvature. Hence
the corollary now follows from Theorem \ref{thm.nu_converge}.
\end{proof}

\section{Further Questions}\label{s.fq}

Two questions immediately arise from Corollary
\ref{cor.nu_converge}.

\begin{enumerate}
\item Can we remove the upper bound restriction on the sectional
curvature? In other words, does the result hold true for any
compact manifold with non-negative sectional curvature?

\item Can the result be extended to an arbitrary compact manifold?

\end{enumerate}

The restriction of the result on a non-negative manifold arises
when we are trying to prove the uniform integrability of
$\rho_{{\mathcal{P}}}$. Thus one way to improve the result is to
obtain a better upper estimate for $\rho_{{\mathcal{P}}}$, and
show that without any restrictions on the sectional curvature of
the manifold, $\rho_{ {\mathcal{P}}}$ is still uniformly
integrable.
\par One possible research problem will be to consider a different metric on
$H(M)$. For example, by considering a $L^2$-metric $G^0$ on
$TH(M)$ where
\begin{align}
G^{0}(X,\ X)  &  := \int_{0}^{1} g\left(X(s),\ X(s)\right)\ ds, \nonumber%
\end{align}
for $X \in T_{\sigma}H(M)$. Instead of considering Wiener space
$W(M)$, one can consider the space of pinned paths, equipped with
pinned Wiener measure and carry out a similar analysis.

\section{Acknowledgements}
The author would like to thank Professor Bruce Driver for all his
help and input into this paper.

\appendix \label{s.appendix}

\section{Trace Class Operators}\label{s.tclass}

A background knowledge on trace class operators can be found in
\cite {MR1009163} and \cite{MR541149}.

\par In general, it is difficult to determine if an integral
operator is trace class. However, the following theorem, taken
from Theorem 2.12 of \cite{MR541149}, gives a condition for which
an integral operator is trace class. See also Section XI.4 of
\cite{MR529429}.

\begin{thm}
\label{thm.trace_class} Let $\mu$ be a Baire measure on a locally
compact space $X$. Let $K$ be a function on $X \times X$ which is
continuous and Hermitian positive, that is
\begin{equation}
\sum_{i,j=1}^{N} \overline{z_{i}}z_{j}K(x_{i}, x_{j})
\geq0\nonumber
\end{equation}
for any $x_{1}, \dots, x_{N} \in X$, $z_{1}, \ldots, z_{N}
\in\mathbb{C}^{N}$ and for any $N$. Then $K(x,x) \geq0$ for all
$x$. Suppose that in addition,
\begin{equation}
\int K(x,x)\ d\mu(x) < \infty.\nonumber
\end{equation}
Then there exists a unique trace class integral operator $A$ such
that
\begin{equation}
(Af)(x) = \int K(x,y)f(y)\ d\mu(y)\nonumber
\end{equation}
and
\begin{equation}
\parallel A \parallel_{1} = \int K(x,x)\ d\mu.\nonumber
\end{equation}

\end{thm}

\begin{prop}
\label{prop.trace1} Let $(Af)(s) = \int_{0}^{1} (\min{s}{t}) f(t)\
dt$. Then $A$ is a trace class operator.
\end{prop}

\begin{proof}
Let $X = [0,1]$ and $\mu$ be Lebesgue measure. Using Theorem
\ref{thm.trace_class}, it suffices to show that $\min{s}{t}$ is
Hermitian positive. Let $z_{1}, z_{2}, \ldots,z_{N}$ be any
complex numbers and let $x_{1}, \ldots,x_{N} \in[0,1]$. The proof
is by induction. Clearly when $N = 1$, it is trivial. Suppose it
is true for all values from $k = 1,2,
\ldots,N-1$. By relabelling, we can assume that $x_{1} \leq x_{k}%
$, $k = 2,\ldots N$. Hence $\min{x_{1}}{x_{k}} = x_{1}$ for any
$k$. Thus
\begin{align}
\sum_{j=1}^{N} \overline{z_{1}}z_{j}(\min{x_{1}}{x_{j}}) +
\sum_{j=1}^{N} \overline{z_{j}}z_{1}(\min{x_{j}}{x_{1}})  &  =
x_{1}\sum_{j=1}^{N}
\overline{z_{1}}z_{j} + x_{1}\sum_{j=1}^{N} \overline{z_{j}}z_{1}\nonumber\\
&  =
x_{1}\left(\sum_{j=1}^{N}\overline{z_{j}}\right)\left(\sum_{j=1}^{N}
z_{j}\right) - x_{1}\sum_{i,j=2}^{N}
\overline{z_{i}}z_{j}.\nonumber
\end{align}
Therefore,
\begin{align}
\sum_{i,j=1}^{N} \overline{z_{i}}z_{j}(\min{x_{i}}{x_{j}})  &  =
x_{1}\left(\sum_{j=1}^{N}\overline{z_{j}}\right)\left(\sum_{j=1}^{N}z_{j}\right)
- x_{1}\sum_{i,j=2}^{N} \overline{z_{i}}z_{j} + \sum_{i,j=2}^{N}
\overline
{z_{i}}z_{j}(\min{x_{i}}{x_{j}})\nonumber\\
&  =
x_{1}\left(\sum_{j=1}^{N}\overline{z_{j}}\right)\left(\sum_{j=1}^{N}
z_{j}\right) + \sum_{i,j=2}^{N}
\overline{z_{i}}z_{j}(\min{c_{i}}{c_{j}}),\nonumber
\end{align}
where $c_{i} = x_{i} - x_{1} \geq 0,\ i = 2, \ldots,N$. Thus by
induction hypothesis,
\begin{equation}
\sum_{i,j=2}^{N} \overline{z_{i}}z_{j}(\min{c_{i}}{c_{j}})
\geq0\nonumber
\end{equation}
and hence
\begin{equation}
\sum_{i,j=1}^{N} \overline{z_{i}}z_{j}(\min{x_{i}}{x_{j}})
\geq0.\nonumber
\end{equation}

\end{proof}

\begin{prop}
\label{prop.trace2} Let ${\mathcal{H}}$ be a separable Hilbert
space. Suppose $A$ is trace class and $B$ is a bounded operator.
Then
\begin{equation}
{\mathrm{{tr}}} \ |AB| \leq\parallel B \parallel{\mathrm{{tr}}} \
|A|\nonumber
\end{equation}
where $\parallel\cdot\parallel$ is the operator norm.
\end{prop}

For a proof, see Theorem VI.25 in \cite{MR751959}.

\section{Perturbation Formulas}\label{s.pf}

\begin{lem}
\label{lem.posdetformula} Let $U$ be a $d\times d$ matrix and for
$r \in \mathbb{N}$, let \beq
\Psi_{r}(U):=\sum_{k=1}^{r}(-1)^{k+1}\frac{1}{k}{\mathrm{{tr}}}\
U^{k}. \nonumber \eeq If $U$ is positive definite, then there
exists $R_{r+1}(U)$ such that
\begin{equation}
\det(I+U)=\exp\left(  \Psi_{r}(U)+R_{r+1}(U)\right)  , \label{eq.detformula2}%
\end{equation}
where
\[
|R_{r+1}(U)|\leq\frac{1}{r+1}{\mathrm{{tr}}}\ U^{r+1}.
\]
If $U$ is any $d \times d$ matrix (not necessarily positive) such
that $\parallel U \parallel< 1$, then
\begin{equation}
\det(I + U) = \exp \left( \Psi_r(U)+R_{r+1}(U) \right), \label{eq.detformula}%
\end{equation}
with
\begin{equation}
R_{r+1}(U) = \sum_{k=r+1}^{\infty}(-1)^{k+1}\frac{1}{k}
{\mathrm{{tr}}} \ U^{k} \nonumber
\end{equation}
such that
\begin{equation}
\left| R_{r+1}(U) \right| \leq \frac{d\parallel U
\parallel^{r+1}}{1 - \parallel U \parallel}. \label{e.bound}
\end{equation}
\end{lem}

\begin{proof}
Let $\{\lambda_{1}, \lambda_{2}, \ldots, \lambda_{d}\}$ be the set
of eigenvalues of $U$. Then
\begin{align}
\det(I + U)  &  = \prod_{i=1}^{d} (1 + \lambda_{i})\nonumber\\
&  = \exp\left(\sum_{i=1}^{d} \ln(1 + \lambda_{i}) \right), \label{eq.det1}%
\end{align}
since by assumption, $\lambda_{i} \geq0$ for all $i$. But by
applying Taylor's Theorem to $\ln(1 + x)$, we have
\begin{equation}
\ln(1 + x) = \sum_{k =1}^{r} (-1)^{k+1}\frac{1}{k}x^{k} +
R_{r+1}(x),\nonumber
\end{equation}
where
\begin{equation}
R_{r+1}(x) = x^{k+1}\frac{(-1)^{r}}{r+1}
\frac{1}{(1+c)^{r+1}}\nonumber
\end{equation}
for some $c \in(0,x)$. Therefore,
\begin{equation}
|R_{r+1}(x)| \leq\frac{1}{r+1}x^{r+1}.\nonumber
\end{equation}
Since
\begin{equation}
{\mathrm{{tr}}} \ U^{k} = \sum_{i=1}^{d} \lambda_{i}^{k},\nonumber
\end{equation}
we have
\begin{align}
\sum_{i=1}^{d} \ln(1 + \lambda_{i})  &  = \sum_{i=1}^{d} \left(
\sum_{k =1}^{r}
(-1)^{k+1}\frac{1}{k}\lambda_{i}^{k} + R_{r+1}(\lambda_{i}) \right)\nonumber\\
& = \sum_{k=1}^{r}
(-1)^{k+1}\frac{1}{k}\tr\ U^k + \sum_{i=1}^{d}R_{r+1}(\lambda_{i})\nonumber \\
&  = \Psi_{r}(U) + R_{r+1}(U),\nonumber
\end{align}
where $R_{r+1}(U) := \sum_{i=1}^{d}R_{r+1}(\lambda_i)$ and
$|R_{r+1}(U)| \leq \sum_{i=1}^d \lambda_i^{r+1} =\frac{1}{r+1}
{\mathrm{{tr}}} \ U^{r+1}$. To prove Equation
(\ref{eq.detformula}), we assume that $\parallel U
\parallel< 1$. Thus the eigenvalues of $U$, $\lambda_{i} < 1$ for
all $i$. Therefore Equation (\ref{eq.det1}) holds. But $\ln(1 +
\lambda_{i}) = \sum_{k=1}^{\infty
}(-1)^{k+1}\frac{1}{k}\lambda_{i}^{k}$ and this sum converges
absolutely. Thus
\begin{align}
\sum_{i=1}^{d}\sum_{k=1}^{\infty}(-1)^{k+1}\frac{1}{k}\lambda_{i}^{k}
&  =
\sum_{k=1}^{\infty}(-1)^{k+1}\frac{1}{k} \sum_{i=1}^{d} \lambda_{i}%
^{k}\nonumber\\
&  = \sum_{k=1}^{\infty}(-1)^{k+1}\frac{1}{k} {\mathrm{{tr}}} \ U^{k}%
.\nonumber
\end{align}
Hence this proves Equation (\ref{eq.detformula}). In this case,
\begin{align}
\left| R_{r+1}(U) \right| &= \left|
\sum_{k=r+1}^{\infty}(-1)^{k+1}\frac{1}{k} {\mathrm{{tr}}} \ U^{k}
\right|  \leq \sum_{k=r+1}^{\infty}\frac{1}{k} \left|
{\mathrm{{tr}}} \ U^{k} \right| \nonumber \\
& \leq \sum_{k=r+1}^{\infty} d\parallel U \parallel^k =
\frac{d\parallel U
\parallel^{r+1}}{1 - \parallel U \parallel}.
\end{align}
Alternatively, Equation (\ref{eq.detformula}) can be proved by a
rewriting of the standard formula,
\begin{align}
\frac{d}{ds}\log(\det(I - sU))  &  = - {\mathrm{{tr}}} ((I - sU)^{-1}%
U)\nonumber\\
&  = - {\mathrm{{tr}}} \ \Big( \sum_{k=0}^{\infty}s^{k}U^{k}U \Big)\nonumber\\
&  = -\sum_{k=0}^{\infty}s^{k} {\mathrm{{tr}}} \ (U^{k+1}
).\nonumber
\end{align}

\end{proof}

Let $\mathcal{J}_1$ be the space of trace class operators,
equipped with the norm $\parallel \cdot \parallel_1$, define by
\beq
\parallel A \parallel_1 = \tr\ |A|. \nonumber \eeq
Let $A \in \mathcal{J}_1$. Suppose $-z^{-1}$ lies in the resolvent
set of $A$, i.e. $-z^{-1} \notin \sigma(A)$. Then the mapping $A
\rightarrow \ln(I + zA)$ is defined and is continuous on $A$. (See
Lemma 15 in Section XI.9.22 of \cite{MR1009163}.) Thus, one can
define a function $\det(I + \cdot)$ on $\mathcal{J}_1$ in the
following way.

\begin{defn}
Define for $-z^{-1} \notin \sigma(A)$, \beq \det (I + zA) =
\exp(\tr\ \ln(I + zA)). \nonumber \eeq
\end{defn}

This function $\det (I + zA)$ is analytic in $z$ and has only
removable singularities at the points $z$ such that $-z^{-1} \in
\sigma(T)$. We can now show that Equation (\ref{eq.detformula})
holds even for trace class operator.

\begin{lem}
Let $A$ be a trace class operator on a Hilbert space
${\mathcal{H}}$, with $\parallel A \parallel< 1$ where
$\parallel\cdot\parallel$ is the operator norm. Then Equation
(\ref{eq.detformula}) holds with $U$ replaced by $A$.
\end{lem}

\begin{proof}
Since $\tr$ is a continuous function on $\mathcal{J}_1$ (See
Theorem 19 in Section XI.9.20 of \cite{MR1009163}.), we have
\begin{eqnarray}
\tr\ \ln(I + zA) & = &\tr\ \lim_{N \rightarrow \infty}\sum_{k=1}^N
\frac{(-1)^{k+1}}{k}(zA)^k = \lim_{N \rightarrow \infty}\tr\
\sum_{k=1}^N
\frac{(-1)^{k+1}}{k}(zA)^k \nonumber \\
& = & \lim_{N \rightarrow \infty}\sum_{k=1}^N
\frac{(-1)^{k+1}}{k}\tr\ (zA)^k = \sum_{k=1}^\infty
\frac{(-1)^{k+1}}{k}\tr\ (zA)^k. \nonumber
\end{eqnarray}
Because $\parallel A \parallel < 1$, setting $z=1$ gives us our
desired result.
\end{proof}

\section{Matrix Inequalities\label{s.matrix}}

\begin{prop}
\label{prop.linalgebra2} If $A,B$ are two $N\times N$ matrices
with $B$ being positive semi definite, then
\begin{equation}
\left\vert {\operatorname{tr}}(AB)\right\vert \leq\left\Vert
A\right\Vert
{\operatorname{tr}}\ B \label{eq.a}%
\end{equation}
and in particular by taking $B=I,$%
\begin{equation}
\left\vert {\operatorname{tr}}\ A\right\vert \leq N\left\Vert
A\right\Vert {.}
\label{eq.b}%
\end{equation}

\end{prop}

\begin{proof}
Let $\{e_{i}\}_{i=1,2,\ldots N}$ be an orthonormal basis of
eigenvectors of $B$ with corresponding eigenvalues,
$\{\lambda_{i}\geq0\}_{i=1,2,\ldots N}.$ Then
\begin{align}
\left\vert \mathrm{{\operatorname{tr}}}\left(  AB\right)
\right\vert & =\left\vert \sum_{i=1}^{N}\langle
ABe_{i},e_{i}\rangle\right\vert =\left|\sum
_{i=1}^{N}\lambda_{i}\langle Ae_{i},e_{i}\rangle\right|\nonumber\\
&  \leq\sum_{i=1}^{N}\lambda_{i}|\langle
Ae_{i},e_{i}\rangle|\leq\sum _{i=1}^{N}\lambda_{i}\left\Vert
A\right\Vert =\left\Vert A\right\Vert \mathrm{{tr}}\ B.\nonumber
\end{align}

\end{proof}

\begin{prop}
\label{prop.linalgebra1}Suppose that $M$ is a positive definite
$N\times N$ matrix and $\alpha>0.$ Then
\begin{equation}
\det\left(  M\right)  \leq\left(  \frac{\operatorname{tr}\left(
M\right) }{N}\right)  ^{N}\leq\alpha^{N}e^{\operatorname{tr}\left(
\alpha
^{-1}M-I\right)  }. \label{eq.1}%
\end{equation}
Moreover if $\alpha\geq1,$ then%
\begin{equation}
\det\left(  M\right)
\leq\alpha^{N}e^{\alpha^{-1}\operatorname{tr}\left(
M-I\right)  }. \label{eq.2}%
\end{equation}

\end{prop}

\begin{proof}
Let $\{\lambda_{i}\}_{i=1}^N$ be the eigenvalues of $M,$ then
\[
\det M=\lambda_{1}\cdots\lambda_{N}=\left(
\lambda_{1}^{1/N}\cdots\lambda _{N}^{1/N}\right)  ^{N}\leq\left(
\sum_{i=1}^{N}\frac{1}{N}\lambda
_{i}\right)  ^{N}=\left(  \frac{\operatorname{tr}\left(  M\right)  }%
{N}\right)  ^{N}.
\]
Now suppose that $\alpha>0.$ Then%
\begin{align*}
\det M  &  =\alpha^{N}\det\left(  \alpha^{-1}M\right)
\leq\alpha^{N}\left(
\frac{\operatorname{tr}\left(  \alpha^{-1}M\right)  }{N}\right)  ^{N}\\
&  \leq\alpha^{N}\left(  1+\frac{\operatorname{tr}\left(  \alpha
^{-1}M-I\right)  }{N}\right)
^{N}\leq\alpha^{N}e^{\operatorname{tr}\left( \alpha^{-1}M-I\right)
},
\end{align*}
where the last inequality follows from the inequality%
\[
\left(  1+\frac{x}{N}\right)  ^{N}\leq e^{x}\text{ for all }N\in
\mathbb{N}\text{ and }x\geq-N.
\]
Indeed, elementary calculus shows that $f\left(  x\right)
:=e^{-x}\left( 1+\frac{x}{N}\right)  ^{N}$ for $x\geq-N$ has a
global maximum of one at $x=0.$

Alternatively,%
\begin{align*}
\det\left(  M\right)   &  =\alpha^{N}\det\left(
\alpha^{-1}M\right)
=\alpha^{N}\det\left(  \alpha^{-1}M-I+I\right) \\
&  =\alpha^{N}\prod_{i=1}^{N}\left(  1+\left(  \alpha^{-1}\lambda
_{i}-1\right)  \right)  \leq\alpha^{N}\prod_{i=1}^{N}e^{\left(
\alpha ^{-1}\lambda_{i}-1\right)
}=\alpha^{N}e^{\operatorname{tr}\left(  \alpha
^{-1}M-I\right)  }%
\end{align*}
wherein we have used the inequality, $1+x\leq e^{x},$ which
results from the convexity of the exponential function. Finally,
by optimizing this inequality over $\alpha>0,$ (take
$\alpha=N^{-1}\operatorname{tr}M),$ the previous inequality
implies $\det M\leq\left(  N^{-1}\operatorname{tr}M\right)  ^{N}.$
\end{proof}

\section{Bounds on Curvature} \label{s.bds}

\begin{defn} \label{def.K}
Let $M$ be a manifold, $u \in O(M)$ and $\Omega_u$ as defined in
Definition \ref{defn.Omegau}. Define for any vectors $x, y, w, z
\in T_oM$, \beq R_u(x, y, w, z) = \langle \Omega_u(x, y)w, z
\rangle, \nonumber \eeq and for any linearly independent vectors
$x, y$, \beq S_u(x, y) = \frac{R_u(x, y, x, y)}{\parallel x
\parallel^2
\parallel y
\parallel^2 - \langle x,y \rangle^2}. \nonumber \eeq
\end{defn}

The curvature tensor, $R_u$ satisfies the following symmetries,
\begin{align}
R_u(x, y, w, z) & = -R_u(y, x, w, z), \label{e.sym1} \\
R_u(x, y, w, z)  & = -R_u(x, y, z, w), \label{e.sym2} \\
R_u(x, y, w, z)  & = R_u(w, z, x, y), \label{e.sym3} \\
 R_u(x, y, w, z)  &+ R_u(w, x, y, z) + R_u(y, w, x, z) = 0.
\label{e.sym4}
\end{align}
The last equality is the Bianchi Identity.
\par $S_u$ is called the sectional curvature and if we define
$S_u(V) := S_u(x, y)$ where $V$ is the plane spanned by any
linearly independent vectors $x, y$, then it can be shown that
$S_u(V)$ is independent of the choice of $x$ and $y$. We will
write \beq
\parallel S_u \parallel = \sup\{|S_u(V)|\ |\ V\ {\rm is\ a\
plane\ in}\ T_0M\}, \nonumber \eeq and \beq K = \sup_{u \in O(M)}
\parallel S_u
\parallel. \nonumber \eeq
Since we are considering a compact manifold $M$, $K$ is finite.
Immediate from this definition, we have \beq |R_u(x, y, x, y)|
\leq
\parallel S_u\parallel\parallel x\parallel^2
\parallel y\parallel^2 \label{e.Kbd} \eeq for any $x, y \in T_oM$.
\par For any vector $v$, Equation (\ref{e.sym3}) shows that
$\Omega_u(v, \cdot)v$ is a symmetric matrix. It follows that
\begin{align}
\parallel \Omega_u(v, \cdot)v \parallel &= \sup_{\parallel w \parallel = 1}
\parallel \Omega_u(v, w)v \parallel = \sup_{\parallel w \parallel = 1}
|\langle \Omega_u(v, w)v, w\rangle| \leq \parallel S_u \parallel
\parallel v \parallel^2. \nonumber
\end{align}
The last inequality follows from Equation (\ref{e.Kbd}). As a
consequence, if we take the supremum over $O(M)$, then \beq
\sup_{u \in O(M)}
\parallel \Omega_u(v, \cdot)v \parallel \leq K\parallel v
\parallel^2. \nonumber \eeq

\begin{prop} \label{p.curve}
If $\sup_{u \in O(M)}\parallel S_u \parallel < \frac{3}{17d}$,
then for any orthonormal frame $\{e_i\}_{i=1}^d \subseteq T_oM$,
\beq \sup_{u \in O(M)}\parallel \Omega_u(e_i, \cdot)e_j \parallel
< \frac{2}{d},\nonumber \eeq where $d$ is the dimension of the
manifold.
\end{prop}

\begin{proof}
Let $x, y, w, z \in T_oM$. Then using Equations (\ref{e.sym1}),
(\ref{e.sym2}) and (\ref{e.sym3}),
\begin{align}
A & := R_u(x+w, y+z, x+w, y+z) - R_u(x+w, y, x+w, y) - R_u(x+w, z,
x+w, z)
\nonumber \\
& - R_u(x, y+z, x, y+z) - R_u(w, y+z, w, y+z) + R_u(w, y, w, y) +
R_u(x, z,
x, z) \nonumber \\
& = R_u(x+w, y, x+w, z) + R_u(x+w, z, x+w, y) - R_u(x, y, x, z) - R_u(x, z, x, y) - R_u(x, y, x, y) \nonumber \\
& - R_u(w, y, w, z) - R_u(w, z, w, y) - R_u(w, z, w, z) \nonumber \\
& = R_u(x, y, w, z) + R_u(w, y, x, z) + R_u(x, z, w, y) + R_u(w,
z, x, y) - R_u(x, y, x, y) - R_u(w, z, w, z) \nonumber \\
& = 2R_u(x, y, w, z) + 2R_u(w, y, x, z) - R_u(x, y, x, y) - R_u(w,
z, w, z). \nonumber
\end{align}
Similarly,
\begin{align}
B & := R_u(x+z, y+w, x+z, y+w) - R_u(x+z, y, x+z, y) - R_u(x+z, w,
x+z,
w) \nonumber \\
&  - R_u(x, y+w, x, y+w) - R_u(z, y+w, z, y+w) + R_u(x, w, x, w) +
R_u(z,
y, z, y) \nonumber \\
& = R_u(x+z, y, x+z, w) + R_u(x+z, w, x+z, y) - R_u(x, y, x, w) - R_u(x, w, x, y) - R_u(x, y, x, y) \nonumber \\
& - R_u(z, y, z, w) - R_u(z, w, z, y) - R_u(z, w, z, w) \nonumber \\
& = R_u(x, y, z, w) + R_u(z, y, x, w) + R_u(x, w, z, y) + R_u(z,
w, x, y) - R_u(x, y, x, y) - R_u(z, w, z, w) \nonumber \\
& = 2R_u(x, y, z, w) + 2R_u(z, y, x, w) - R_u(x, y, x, y) - R_u(w,
z, w, z). \nonumber
\end{align}
Thus
\begin{align}
A - B & = 2R_u(x, y, w, z) + 2R_u(w, y, x, z) - 2R_u(x, y, z, w) - 2R_u(z, y, x, w)\nonumber \\
& = 2R_u(x, y, w, z) + 2R_u(x, y, w, z) + 2R_u(w, y, x, z) +
2R_u(y, z, x,
w) \nonumber \\
& = 4R_u(x, y, w, z) + 2R_u(w, y, x, z) + 2R_u(x, w, y, z) \nonumber \\
& = 4R_u(x, y, w, z) - 2R_u(y, x, w, z) \nonumber \\
& = 6R_u(x, y, w, z). \nonumber
\end{align}
The second last equality follows from the Bianchi Identity. Hence
for any unit vectors $x, y, w, z$,
\begin{align}
| A |  & \leq \parallel S_u \parallel \left[\parallel x+w
\parallel^2 \parallel y+z
\parallel^2 + 2\parallel x+w
\parallel^2
+ 2\parallel y+z \parallel^2 + 2 \right]
\nonumber \\
& \leq \parallel S_u \parallel [4^2 + 2\cdot8 + 2] = 34\parallel
S_u
\parallel. \nonumber
\end{align}
Similarly, $| B | \leq 34 \parallel S_u \parallel$. Therefore, for
any unit vectors $x, y, w, z$, \beq | R_u(x, y, w, z) | =
\frac{1}{6}| A - B | \leq \frac{34}{3}\parallel S_u
\parallel. \nonumber \eeq and hence \beq \sup_{\parallel x \parallel
= \parallel y \parallel =
\parallel w \parallel = \parallel z \parallel = 1}| R_u(x,
z, w, z) | \leq \frac{34}{3}\parallel S_u \parallel. \nonumber
\eeq Now for any $u \in O(M)$ and orthonormal frame
$\{e_i\}_{i=1}^d$, let $\Omega_{ij}(y) := \Omega_u(e_i, y)e_j$.
Then
\begin{align}
\sup_{\parallel y \parallel = 1} \parallel \Omega_{ij}(y)
\parallel & = \sup_{\parallel y \parallel = 1} \left\langle \Omega_{ij}(y),
\frac{\Omega_{ij}(y)}{\parallel \Omega_{ij}(y)
\parallel}\right\rangle
 = \sup_{\parallel y \parallel = 1} \left|R_u\left(e_i, y, e_j,
\frac{\Omega_{ij}(y)}{\parallel \Omega_{ij}(y)\parallel}\right)\right| \nonumber \\
& \leq \sup_{\parallel x \parallel =  \parallel y \parallel =
\parallel w \parallel = \parallel z \parallel = 1}| R_u(x,
y, w, z) | \leq \frac{34}{3}\parallel S_u \parallel. \nonumber
\end{align}
So, for any $u \in O(M)$ and orthonormal basis $\{e_i\}_{i=1}^d
\subseteq T_oM$, $\parallel \Omega_{ij} \parallel \leq
\frac{34}{3}\parallel S_u
\parallel$. If we choose $\sup_{u \in O(M)}\parallel S_u \parallel < 3/(17d)$, then
\beq \sup_{u \in O(M)}\parallel \Omega_u(e_i, \cdot)e_j \parallel
< \frac{34}{3}\cdot \frac{3}{17d}  = \frac{2}{d}. \nonumber \eeq
\end{proof}

\nocite{*}
\bibliographystyle{amsalpha}
\bibliography{Wiener}



\end{document}